\documentclass[12pt]{amsart}
\newcounter{dummy}
\usepackage{amsmath,amsthm,amscd,amsfonts,wasysym,amssymb,epic,eepic,bbm,tikz-cd,mathtools,multicol,enumitem,mathrsfs}
\usepackage[pagebackref,colorlinks=true,linkcolor=blue,citecolor=blue]{hyperref}
\usepackage[noabbrev]{cleveref}
\usepackage{MnSymbol}
\makeatletter
\newcommand\myitem[1][]{\item[#1]\refstepcounter{dummy}\def\@currentlabel{#1}}
\makeatother



\allowdisplaybreaks
\newtheorem{thm}{Theorem}[section]
\newtheorem{cor}[thm]{Corollary}
\newtheorem{conj}[thm]{Conjecture}
\newtheorem{lem}[thm]{Lemma}

\newtheorem{prop}[thm]{Proposition}
\theoremstyle{remark}
\newtheorem*{rem}{Remark}

\newcounter{remarkscounter}

\newcommand{\ga}{\gamma}

\newcommand{\Lam}{\Lambda}

\numberwithin{equation}{section}
\newcommand{\A}{\mathbb{A}}
\newcommand{\GL}{\mathrm{GL}}

\newcommand{\SL}{\mathrm{SL}}

\newcommand{\ZZ}{\mathbb{Z}}

\newcommand{\QQ}{\mathbb{Q}}

\newcommand{\lto}{\longrightarrow}
\newcommand{\lra}{\longrightarrow}

\newcommand{\OO}{\mathcal{O}}
\newcommand{\CC}{\mathbb{C}}
\newcommand{\RR}{\mathbb{R}}
\newcommand{\GG}{\mathbb{G}}

\newcommand{\Hom}{\mathrm{Hom}}

\newcommand{\quash}[1]{}

\theoremstyle{definition}

\newenvironment{psmatrix}
  {\left(\begin{smallmatrix}}
  {\end{smallmatrix}\right)}

\renewcommand{\bar}{\overline}
\numberwithin{equation}{subsection}

\newcommand{\one}{\mathbbm{1}}

\newcommand{\norm}[1]{\lVert#1\rVert}

\allowdisplaybreaks

\makeindex

\begin{document}

\title[Triple product $L$-functions]{Triple product $L$-functions \\ and the fiber bundle method}

\author{Jayce R. Getz}
\address{Department of Mathematics\\
Duke University\\
Durham, NC 27708}
\email{jgetz@math.duke.edu}

\author{Miao (Pam) Gu}
\address{Department of Mathematics\\
Aarhus University\\
DK-8000 Aarhus, Denmark}
\email{pmgu@umath.au.dk}

\author{Chun-Hsien Hsu}
\address{Department of Mathematics\\
University of Chicago\\
Chicago, IL 60637}
\email{chunhsien@uchicago.edu}

\author{Spencer Leslie}
\address{Department of Mathematics\\
Boston College \\
Chestnut Hill, MA 02467}
\email{spencer.leslie@bc.edu}

\thanks{The first author is thankful for partial support provided by NSF grant DMS-1901883 and NSF RTG grant DMS-2231514. The second author appreciates partial support provided by an AMS-Simons Travel Grant. The fourth author was partially supported by NSF grants DMS-2200852 and DMS-2440741. Any opinions, findings, and conclusions or recommendations expressed in this material are those of the authors and do not necessarily reflect the views of the National Science Foundation. 
}

\begin{abstract}
We introduce multivariable zeta integrals that unfold to Euler products representing the triple product $L$-function times a product of $L$-functions with known analytic properties.  Motivated by this integral representation, we formulate a conjectural Poisson summation formula and show it implies the analytic properties of triple product $L$-functions. Finally, we propose a strategy, the fiber bundle method, to reduce the conjectural formula to two known Poisson summation formulae along with certain local compatibility statements. 
\end{abstract}

\maketitle

\setcounter{tocdepth}{2}
\begin{multicols}{2}
\tableofcontents
\end{multicols}

\section{Introduction}

Let $F$ be a number field with ring of adeles $\A_F.$  If $G$ is a reductive $F$-group, we denote by ${}^LG$ its Langlands dual.  
Consider the tensor product
\begin{align*}
    \otimes^2:{}^L(\GL_{r_1} \times \GL_{r_2}) \lto {}^L\GL_{r_1r_2}.
\end{align*}
Langlands functoriality predicts that there is a corresponding functorial transfer of automorphic representations of $\GL_{r_1}(\A_F) \times \GL_{r_2}(\A_F)$ to $\GL_{r_1r_2}(\A_F).$  We refer to this transfer as the \textbf{automorphic tensor product}.   The tensor product is a cornerstone of representation theory. Likewise, constructing the automorphic tensor product is a fundamental open problem.

One might hope to approach this problem via converse theory.  
Let $r_1,r_2,r_3$ be a triple of positive integers and let $\pi=\otimes_{i=1}^3\pi_i$ where the $\pi_i$ are all cuspidal automorphic representations of $\GL_{r_i}(\A_F)$.  
Thus $\pi$ is a cuspidal automorphic representation of $\prod_{i=1}^3\GL_{r_i}(\A_F)$.
We denote by 
$$
L(s,\pi, \otimes^3)=L(s, \pi_1\times \pi_2 \times \pi_3)
$$
the corresponding triple product $L$-function.  
It is the Langlands $L$-function defined by the tensor product
representation
\begin{align*}
\otimes^3:{}^L(\GL_{r_1} \times \GL_{r_2} \times \GL_{r_3}) \lto \GL_{r_1r_2r_3}(\CC).
\end{align*}
If we assume the expected analytic properties of triple product $L$-functions, then the converse theorem of Cogdell and Piatetski-Shapiro \cite{Cogdell:PS:ConverseII} implies the existence of the automorphic tensor product.  This is explained in detail in \cite{Getz:Hahn:Yao}, where the authors also showed that knowledge of the analytic properties of triple product $L$-functions implies the Ramanujan conjecture.

\subsection{A novel integral representation}

In this paper we introduce certain multivariable zeta integrals.  We prove that they converge absolutely in an appropriate cone and unfold to an Euler product representing the triple product $L$-function times a product of $L$-functions with known analytic properties.   The construction is inspired by the integral representation for triple product $L$-functions obtained by Garrett \cite{GarrettTripleAnnals} (an adelic treatment was given in \cite{PS:Rallis:triple}).  The construction is subtle.  It does not fit neatly into established paradigms in the study of integral representations of $L$-functions; we view this as a strength of the construction.
It would be very interesting to see if it can be generalized.

To construct the zeta integrals, we introduce the notion of an affine $\Psi$-bundle.  
  Affine $\Psi$-bundles allow us to formulate a precise conjectural Poisson summation formula.  It is similar to, but distinct from, the related conjectures appearing in  \cite{BK-lifting,LafforgueJJM,NgoSums,Ngo:Hankel,SakellaridisSph}.
  We prove that this conjectural summation formula implies the expected analytic properties of triple product $L$-functions for $r_i>2.$  Finally, we end the paper by outlining a strategy to reduce this Poisson summation conjecture to local questions.

\subsection{Affine $\Psi$-bundles}
 \label{ssec:period}
 We now introduce the groups and schemes  needed to define our zeta integrals. 
 To orient the reader, if $X$ is a scheme of finite type over $F,$ we will denote by $X^\circ \subset X$ an open dense orbit under a suitable group action.

 For $r_1,r_2,r_3 \in \ZZ_{>0}$, let 
\begin{align*}
\GL_{\underline{r}}:=\GL_{r_1} \times \GL_{r_2} \times \GL_{r_3}.
\end{align*}
We will often be treating spaces that come in triples in the current work, and we will use notation such as $\underline{r-2}$ to denote $(r_1-2,r_2-2,r_3-2)$, $\underline{1}$ to denote $(1,1,1),$ etc.  Let $I_{\underline{r}}$ be the identity of $\GL_{\underline{r}}(F)$.  For an $F$-algebra  $R$, consider the group scheme 
\begin{align*} 
H(R):=\{(h_1,h_2,h_3) \in \GL_{\underline{2}}(R): \det h_1=\det h_2=\det h_3\} \index{$H$}
\end{align*}
and let 
$
\nu:H \to \GG_m 
$ 
denote the 
character $(h_1,h_2,h_3) \mapsto \det h_1.$ For technical reasons explained in \S \ref{sec: z-ext}, we work with a certain extension $H^e \to H.$
Let
$$
G:=\GL_{\underline{r}} \times \GL_{\underline{r-2}} \times H^e.
$$
We write $g=(g^{(1)},g^{(2)},g^{(3)})$ for the projections to the three factors.
Let
$$
X^\circ_{\underline{r}}:=N_{\underline{r-2},\underline{2}} \backslash \GL_{\underline{r}},
$$
where $N_{\underline{r-2},\underline{2}}$ is the unipotent radical of the standard parabolic $P_{\underline{r-2},\underline{2}}$ of type $(\underline{r-2},\underline{2})$ (see \eqref{Prdef}).  We denote by $X_{\underline{r}}:=\overline{X^\circ_{\underline{r}}}^{\mathrm{aff}}$ the affine closure of $X^\circ_{\underline{r}};$ this is a singular spherical 
$G$-variety with open orbit $X^\circ_{\underline{r}}$.

Let $M_{r,s}(R)$ be the space of $r \times s$ matrices with coefficients in $R.$ 
Assume $r_i \geq 3$ for all $i$ and let\index{$\mathcal{M}$} \index{$\mathcal{M}_i$}
\begin{align*}
\begin{split}
    \mathcal{M}_{i}(R)&=\{(u_1,u_2) \in M_{r_i-2,2}(R): u_1 \wedge u_2=0\},\\
    \mathcal{M}^{ \circ}_i:&=\mathcal{M}_{i}-\{0\}, \quad
    \mathcal{M}:=\prod_{i=1}^3\mathcal{M}_{i},  \quad 
    \mathcal{M}^{ \circ}:=\prod_{i=1}^3\mathcal{M}_{ i}^\circ.  \end{split}
\end{align*} 
Thus $\mathcal{M}_i^\circ(F) \subset M_{r_i-2,2}(F)$ is the subset of matrices of rank $1.$ The group $\GL_{\underline{r-2}}\times H^e$ acts on  $\mathcal{M}$ as in \eqref{R:def} and $\mathcal{M}^\circ$ is the unique open orbit. We thus obtain a $G$-scheme $X_{\underline{r}}\times \mathcal{M}$ with open orbit $X_{\underline{r}}^\circ\times \mathcal{M}^\circ$.

Let \index{$V_3$}
\begin{align*}
V_3:=\GG_a^2 \otimes \GG_a^2 \otimes \GG_a^2
\end{align*} 
denote the $8$-dimensional representation of $H$ given on pure tensors by  \index{$v.h$}
\begin{align*}
(v_1 \otimes v_2 \otimes v_3).(h_1,h_2,h_3)= (\det h_1)^{-1} (v_1h_1 \otimes v_2h_2 \otimes v_3h_3).
\end{align*}
This action extends to $H^e$ so that $\GL_{\underline{r-2}}\times H^e$ acts on $\mathcal{M}\times V_3$; see \eqref{R:def} for details. 
In \S\ref{sec:Y:def}, we construct a (singular) closed subscheme $Y \subset \mathcal{M} \times V_3$ that is preserved under this action such that $Y_{\mathcal{M}^{\circ}}:=Y\times_{\mathcal{M}}\mathcal{M^\circ}$ is a vector bundle of rank $4$ over $\mathcal{M}^{\circ}$ (see Lemma \ref{lem:Y:vb}). 
We denote by $Y^\circ\subset Y_{\mathcal{M}^{\circ}}\subset Y$ the open $(\GL_{\underline{r-2}} \times H^e)$-orbit and set
$$
X:= X_{\underline{r}} \times Y.
$$

 Our zeta integrals are built from an equivariant torsor over 
 $$
 X^\circ_{\underline{r}} \times Y.
 $$
 More precisely, consider the canonical $N_{\underline{r-2},\underline{2}}$-torsor
$$
\GL_{\underline{r}} \times Y \lto X^\circ_{\underline{r}} \times Y.
$$
It is $G$-equivariant. In \eqref{Psi}, we define a morphism
\begin{align*}  \begin{split}
    \Psi:\GL_{\underline{r}} \times Y \times  N_{\underline{r-2},\underline{2}}(R)  &\lto \GG_a\end{split}.
\end{align*}
In \S \ref{ssec:gen whit} we abstract the key properties of the pair 
$(\GL_{\underline{r}} \times Y\to  X^\circ_{\underline{r}} \times Y,\Psi)$ to define the notion of an \textbf{affine $\Psi$-bundle}.  Although we will not mention it again in the introduction, this notion isolates the geometry that underlies  our work.

\subsection{The integral representation}

Using $\Psi$ and a nontrivial additive character $\psi:F \backslash \A_F \to \CC^\times,$ in \S \ref{ssec:adelic:sec} we construct a Hermitian line bundle $\mathcal{H}$ over $X^\circ(\A_F).$   The bundle is equipped with an action $\mathcal{R}$ of $G(\A_F).$

For suitable sections $f$ of $\mathcal{H}$ over $X^{\circ}(\A_F)$ and $g\in G(\A_F),$ set \index{$\Theta_f(g,g',h)$}
\begin{align*}\begin{split}
    \Theta_f(g):&=\sum_{x \in X^{\circ}(F)}\mathcal{R}_u(g)f(x) \end{split}
\end{align*}
whenever this expression is absolutely convergent.  Here $\mathcal{R}_u$ is the unitarily normalized regular action (see \eqref{Ru}).

Let $\pi$ and $\pi'$ be cuspidal automorphic representations of $A_{\GL_{\underline{r}}} \backslash \GL_{\underline{r}}(\A_F)$ and $ A_{\GL_{\underline{r-2}}} \backslash \GL_{\underline{r-2}}(\A_F)$ with central characters $\omega$ and $\omega'$ satisfying $\omega\omega'=1.$ In \S\ref{sec:local zeta}, we define an unramified quasi-character
\begin{align*}
\eta_{\underline{s},s} :\GL_{\underline{r-2}}(\A_F)\times H^e(\A_F)\lra \mathbb{C}^\times,
\end{align*}
where  $(\underline{s},s) \in  \CC^{\underline{1}} \times \CC.$
We continue to denote by $\eta_{\underline{s},s}:G(\A_F) \to \CC^\times$ the pull-back of $\eta_{\underline{s},s}$ along the canonical projection to $G(\A_F) \to \GL_{\underline{r-2}}(\A_F) \times H^e(\A_F).$

For smooth $\varphi\in \pi$ and $\varphi'\in \pi'$, we define the zeta integral \index{$Z(\varphi,\varphi',f,\underline{s},s)$}
\begin{align} \label{glob:zeta}
&Z(\varphi,\varphi',f,\underline{s},s):=\int_{ [G]' }\Theta_f(g)\varphi(g^{(1)})\varphi'(g^{(2)})\eta_{\underline{s},s}(g)dg,
\end{align}
where $[G]':=[\GL_{\underline{r}}\times  \GL_{\underline{r-2}}(\A_F)]\times H^e(F)\backslash H^e(\A_F)'$ is defined as in \eqref{He1}.

Our first main theorem, Theorem \ref{thm:intro1}, states that under appropriate assumptions on $f,$ the integral above converges absolutely in a suitable cone and unfolds to an Eulerian integral.  Let $U_{r_i}$ be the unipotent radical of the Borel subgroup of $\GL_{r_i}$ consisting of upper triangular matrices. By abuse of notation, let $\psi$ denote the usual generic character of $U_{r_i}(\A_F)$ induced by $\psi$. 
We also denote by $\psi$ the character 
\begin{align*}
U_{\underline{r}}(\A_F) &\lto \CC^\times\\
(u_1,u_2,u_3) &\longmapsto \psi(u_1)\psi(u_2)\psi(u_3).
\end{align*}
Define the usual Whittaker functions \index{$W^{\varphi}_{\psi}$} \index{$W^{\varphi'}_{\overline{\psi}}$}
\begin{align*}
W^{\varphi}_{\psi}(g):=\int_{[U_{\underline{r}}]}\varphi(ug)\overline{\psi}(u)du \quad  \textrm{ and } \quad
W^{\varphi'}_{\overline{\psi}}(g):=\int_{[U_{\underline{r-2}}]}\varphi'(ug)\psi(u)du.
\end{align*}
Let $N_0 \leq \SL_{\underline{2}}$ be the unipotent group defined in \eqref{N0:def} and let $y_0:=(m_0,v_0) \in Y(F)$ be the representative for the open $\GL_{\underline{r-2}} \times H^e$ orbit of \eqref{Ycirc}.

\begin{thm} \label{thm:intro1}
Assume $f \in C^\infty(X^\circ(\A_F),\mathcal{H})$ satisfies the assumptions in \S \ref{sec:conv}.
For $(\mathrm{Re}(\underline{s}),\mathrm{Re}(s))$ in a suitable nonempty cone in $\RR^{\underline{1}}_{>0} \times \RR_{>0},$ the integral $Z(\varphi,\varphi',f,\underline{s},s)$ converges absolutely and is equal to $\mathrm{meas}(Z_{\GL_{\underline{2}}}(F) \backslash Z_{\GL_{\underline{2}}}(\A_F)^1)$ times
\begin{align*}
    &\int_{Z_{\GL_{\underline{2}}}(\A_F) N_0(\A_F)\backslash H^e(\A_F)}\int_{U_{\underline{r-2}}(\A_F) \backslash \GL_{\underline{r-2}}(\A_F)} \int W^{\varphi}_{\psi}(g)W^{\varphi'}_{\overline{\psi}}(g)\mathcal{R}_u\left(g,g',h \right)f(x_0)\eta_{\underline{s},s}(g',h)dgdg'dh,
\end{align*}
where the inner integral is over 
$N_{\underline{r-2},\underline{2}}(\A_F) \backslash \GL_{\underline{r}}(\A_F).$
\end{thm}
\noindent 
Theorem \ref{thm:intro1} is inspired by the integral representation for triple product $L$-functions for $\GL_{\underline{2}}$ due to Garrett \cite{GarrettTripleAnnals,PS:Rallis:triple}.   To be clear, it is not at all obvious how to generalize the construction in loc.~cit. to higher rank.  Indeed, Theorem \ref{thm:intro1} is not a generalization in the strict sense; the assumption $r_i \geq 3$ is built into our integral representation.

\subsection{Local integrals}
Let $v$ be a place of $F.$  We omit $v$ from notation, writing $F:=F_v,$ etc.  Let $\pi$ and $\pi'$ be generic irreducible admissible representations of $\GL_{\underline{r}}(F)$ and $\GL_{\underline{r-2}}(F),$ and let $\mathcal{W}(\pi,\psi)$ and $\mathcal{W}(\pi',\overline{\psi})$ be the usual Whittaker models.

Let $(W,W') \in \mathcal{W}(\pi,\psi) \times \mathcal{W}(\pi',\overline{\psi})$ and let $f \in C^\infty(X^\circ(F) ,\mathcal{H}).$  The integral in Theorem \ref{thm:intro1} is visibly Eulerian, with local factors of the form \index{$Z(W,W',f,\underline{s},s)$}
\begin{align*}
\begin{split}
&Z(W,W',f,\underline{s},s)\\&:=\int_{Z_{\GL_{\underline{2}}}(F) N_0(F)\backslash H^e(F)}
\int_{U_{\underline{r-2}}(F) \backslash \GL_{\underline{r-2}}(F)}\int W(g)W'(g') \mathcal{R}_u(g,g',h)f(x_0)\eta_{\underline{s},s}(g',h)dg dg'dh.\end{split}
\end{align*}
Here the inner integral is over $N_{\underline{r-2},\underline{2}}(F) \backslash \GL_{\underline{r}}(F).$
   We prove in Proposition \ref{prop:ram:conv} that the integral $Z(W,W',f,\underline{s},s)$ converges for $f$ satisfying appropriate assumptions if  $(\mathrm{Re}(\underline{s}),\mathrm{Re}(s))$ lies in a suitable cone in $\RR^{\underline{1}}_{>0} \times \RR_{>0}$. 

To connect this with the triple product $L$-function, we consider the case in which $F$ is unramified over $\QQ_p$ with $p\geq 5$, and $f=b^{\mathrm{nai}}$ is the ``na\"ive basic function'' defined in \eqref{naive}.  This is essentially the characteristic function of the integral points, and is an approximation of the true basic function discussed in \S\ref{ssec:desi:loc:ss} and \S \ref{sec:add:des} (see \S\ref{ssec:unr}) .

\begin{thm} \label{thm:loc:comp:intro} Assume that $\pi,\pi',\psi$ are unramified and let $(W,W') \in \mathcal{W}(\pi,\psi)^{\GL_{\underline{r}}(\OO_F)}\times \mathcal{W}(\pi',\overline{\psi})^{\GL_{\underline{r-2}}(\OO_F)}$ be the unique vectors satisfying $W(I_{\underline{r}})=W'(I_{\underline{r-2}})=1.$ 
    One has that
    \begin{align*}
Z&(W,W',b^{\mathrm{nai}},\underline{s},s)\\=&\sum_{\substack{\underline{k} \in \ZZ_{\geq 0}}}\sum_{\ell=0}^\infty q^{-\ell (s+1/2)}\prod_{i=1}^3 L(C_i(\underline{s},s)+\tfrac{1}{2}
, \pi_i \times \pi_i') 
   \sum_{n_i=0}^{r_i-2}\frac{(-1)^{n_i} \mathbb{S}_{0,\dots,0,n_i-\ell-k_i}(\alpha_i)\mathrm{Tr}(\wedge^{n_i} \alpha'_i)}{q^{n_i(C_i(\underline{s},s)+1/2)+k_i(s_i+1/2)}}.
    \end{align*}
\end{thm}
\noindent Here $q$ is the order of the residue field of $F,$ $C_i(\underline{s},s)$ is defined in \eqref{CA}, and $\alpha_{i} \in \GL_{r_i}(\CC),$ $\alpha_{i}' \in \GL_{r_i-2}(\CC)$ are the Langlands class of $\pi_i,$ $\pi_i',$ respectively.

\begin{cor} \label{cor:loc:comp:intro}
Assume additionally $\pi$ and $\pi'$ are tempered.  
Let $\varepsilon>0$ and assume that 
$\mathrm{Re}(s) >\varepsilon$ and that for each $i$, $\mathrm{Re}(s_i) > \varepsilon$ and $C_i(\underline{s},s)>\varepsilon$.  Then
$$
\frac{Z(W,W',b^{\mathrm{nai}},\underline{s},s)}{
L(s+\tfrac{1}{2},\pi^\vee, \otimes^3)\prod_{i=1}^3L(C_i(\underline{s},s)+\tfrac{1}{2},\pi_i \times \pi_i')L(s_i+\tfrac{1}{2},\pi_i^\vee) }=1+O_\varepsilon(q^{-1-\varepsilon}).
$$
Moreover, even if $\pi$ and $\pi'$ are not tempered, the above quotient
    is a rational function of $q^{\pm (r_i-2)^{-1}s},q^{\pm (r_i-2)^{-1}s_i}.$  
\end{cor}
\noindent  These results are restated as Theorem \ref{thm:loc:comp} and Corollaries \ref{cor:loc:comp} and \ref{cor:rational} below. We have not dwelt on the problem of explicitly computing $Z(W,W',b^{\mathrm{nai}},\underline{s},s)$ as a rational function of $q^{\pm (r_i-2)^{-1} s
},q^{\pm (r_i-2)^{-1}s_i}$ since $b^{\mathrm{nai}}$ is merely an approximation to the actual basic function.  One must replace $b^{\mathrm{nai}}$ with the true basic function to obtain the meromorphic continuation of the corresponding global $L$-functions.  

\begin{rem}
As pointed out to the second author by P.~Sarnak,
it might be possible to study triple product $L$-functions at the edge of the critical strip using the work in this paper.  In principle, the na\"ive basic function and our integral representation already yield an expression for the poles of triple product $L$-functions at $s=1,$ even without the conjectural Poisson summation formula.  
\end{rem}

\subsection{The fiber bundle method} \label{ssec:Fb:method}

Of course, our main motivation for studying the zeta integrals $Z(\varphi,\varphi',f,\underline{s},s)$ is to prove the analytic properties of triple product $L$-functions.  To this end, we formulate a version of the Poisson summation conjecture for $X$ in Conjecture \ref{PS:conj} below.  We note that this version of the Poisson summation conjecture does not reduce to others in the literature as the relevant spaces of local invariant linear functionals are not one-dimensional. They are only one-dimensional after integrating over $H^{\mathrm{der}}(F_v).$

In Theorem \ref{thm:func:eqn}, we prove that Conjecture \ref{PS:conj} implies the analytic continuation and functional equation of $Z(\varphi,\varphi',f,\underline{s},s)$ for $f$ satisfying suitable local assumptions.    The following corollary is restated and proved as Corollary \ref{cor:mero} below:

\begin{cor} \label{cor:mero:intro}
 Assume Conjecture \ref{PS:conj} and the local nonvanishing hypothesis \ref{nonvanish}.  Let $\pi$ be a cuspidal automorphic representation of $A_{\GL_{\underline{r}}} \backslash \GL_{\underline{r}}(\A_F).$  Then $L(s,\pi,\otimes^3)$ admits a meromorphic continuation to the plane.  
\end{cor}

We expect known techniques will imply the local nonvanishing hypothesis; see the sentence below \ref{nonvanish} for details. 
Thus one is left with proving Conjecture \ref{PS:conj}.  While we do not yet have a proof of this conjecture, we present a nontrivial reduction to a simpler situation via the \textbf{fiber bundle method}.  This is discussed in more detail in \S \ref{sec:fb} below, so we content ourselves with the main idea in the introduction.  
In \S \ref{sec:fb} we construct a diagram
\begin{equation*}
\begin{tikzcd}
X_{\underline{r}} \times \mathcal{M}& \arrow[l,"\phi_1",swap] X \arrow[r,"\phi_2"] & V_3.
\end{tikzcd}
\end{equation*}
 General fibers of $\phi_1$ are vector spaces and general fibers of $\phi_2$ are simpler affine $\Psi$-bundles. The Poisson summation formula for the fibers of $\phi_1$ is thus classical.  The Poisson summation formula for the fibers of $\phi_2$ is nontrivial, but it is proved in the recent paper \cite{GGH} under suitable assumptions that should be no loss of generality for our purposes.  What remains is to combine these two Poisson summation formulae.  This is the subject of future work, but in any cases the techniques involved are orthogonal to those of this paper.  Crucially, it appears to be a  local problem.

\subsection{Outline} \label{ssec:outline} We begin with function-theoretic preliminaries in \S \ref{sec:FT:prelim}.  After this, we pass to  geometric preliminaries in \S \ref{sec:groups:orbits}.
In \S \ref{sec:SS} we discuss conjectures regarding Schwartz spaces and make our conjectural Poisson summation formula precise.
We bound the local integrals in \S \ref{sec:local:int}.  We turn to the computation of the unramified local integrals in \S \ref{sec:unramified}.  Motivated by this computation, we introduce additional local desiderata for the Schwartz space in \S \ref{sec:add:des}.  We prove the convergence of the global integrals in \S \ref{sec:conv}.  We then prove that they unfold in \S \ref{sec:unfold}.

In \S \ref{sec:apply} we prove that our conjectural Poisson summation formula implies the meromorphic continuation of triple product $L$-functions.  Finally, in \S \ref{sec:fb} we propose a reduction of the Poisson summation conjecture to local assertions.  

To aid the reader we have appended an index of notation.

\subsection*{Acknowledgments}
 The first author began this paper while at the Institute for Advanced Study in the Spring of 2018, and he thanks the Charles Simonyi endowment for their support.  He also thanks the Postech Mathematics Institute and Y.~Choie for hospitality during a sabbatical in 2020-2021 and again for a month in 2024.

The first author thanks B.~C.~Ng\^o and F.~Shahidi for many useful conversations and P.~Garrett, D.~Jiang, and B.~Liu for encouragement.  He learned of the fiber bundle method from D.~Kazhdan during a visit to the Einstein Institute of Mathematics in 2018; he thanks Kazhdan for the invitation and his insight.  He thanks  A.~Slipper for the proof of Lemma \ref{lem:Pl:Iso}. W-T. Gan also deserves acknowledgment for questions that led to a simplification of the definition of the period integral studied in this paper.  Finally, he thanks H.~Hahn for her help with editing and her constant encouragement.

\section{Function-theoretic preliminaries}
 \label{sec:FT:prelim}
\subsection{Quasi-characters and norms} \label{ssec:quasi} 

Let $F$ be a local field of characteristic zero and let $|\cdot|$ be the normalized norm on $F$. It is the usual Euclidean norm if $F=\RR$ and the square of the usual norm if $F=\CC$. If $F$ is non-Archimedean with ring of integers $\mathcal{O}_F$ and uniformizer $\varpi$, then $q:=|\varpi|^{-1}$ is the cardinality of the residue field of $\mathcal{O}_F$. For a quasi-character $\eta:F^\times\to \CC^\times$ and $s\in \CC$, let $\eta_s:=\eta|\cdot|^s$. We define $\mathrm{Re}(\eta)$ to be the unique real number such that $\eta_{-\mathrm{Re}(\eta)}$ is unitary.  If $X$ is an integral separated scheme of finite type over $F,$ we make use of the notion of log norms
\begin{align*}
\sigma_X:X(F) \lto \RR_{>0}
\end{align*}
(\cite[\S 1.2]{BP:Ast} is a convenient reference).  
If $X$ is an affine algebraic group, we drop it from notation, writing $\sigma:=\sigma_X.$

Let $X_1$ and $X_2$ be integral separated schemes of finite type over $F.$ For $(x_1,x_2) \in X_1(F) \times X_2(F),$ 
\begin{align} \label{prod:log}
\sigma_{X_1}(x_1)^{1/2}\sigma_{X_2}(x_2)^{1/2} \ll \sigma_{X_1 \times X_2}(x_1,x_2) \ll \sigma_{X_1}(x_1)\sigma_{X_2}(x_2).
\end{align}
By \cite[Lemma 1.2.1]{BP:Ast}, if $\phi:X_1 \to X_2$ is a morphism, then
\begin{align} \label{pullback}
\sigma_{X_1}(x_1) \gg \sigma_{X_2}(\phi(x_1)).
\end{align}
If, moreover, $\phi$ is finite, then
\begin{align} \label{iso:inv}
\sigma_{X_1}(x_1) \asymp \sigma_{X_2}(\phi(x_1)).
\end{align}

Log norms are only defined up to a notion of equivalence discussed in \cite[\S 1.4]{BP:Ast}.  When applying local bounds to the adelic situation, it is sometimes important to fix a particular representative in the equivalence class.  For our purposes, fixing the following two representatives is enough:
\begin{align*}
\sigma_{\GG_a^n}(x_1,\dots,x_n):&=1+\max(1,|x_1|,\dots,|x_n|),\\
\sigma_{\GL_n}(g):&=1+\max(\norm{g}_{\mathrm{op}},\norm{g^{-1}}_{\mathrm{op}}). 
\end{align*}
Here $\norm{g}_{\mathrm{op}}$ is the operator norm with respect to a choice of metric on $F^n.$  In the non-Archimedean case, we take the metric to be the max norm, so $\norm{g}_{\mathrm{op}}$ is the maximum of the norms of the entries of $g.$

\subsection{Measures}\label{ssec:measure} 
Let $Q$ be a locally compact Hausdorff group and let $Z$ be a Hausdorff $Q$-homogeneous space. Then the stabilizer of any point of $Z$ is closed in $Q$. 
Let $\chi:Q \to \RR_{>0}$ be a continuous quasi-character. By a \textbf{$\chi$-eigenmeasure} on $Z$ we mean a positive Radon measure $dz$ such that 
$d(z.q)=\chi(q)dz.
$
An \textbf{eigenmeasure} is a $\chi$-eigenmeasure for some $\chi.$   
Below we will make use of the fact that $\chi$-eigenmeasures are uniquely determined up to scalars in the sense of \cite[Proposition B.1.6]{BHV}.  
We refer to this result as uniqueness of eigenmeasures.
Let $z \in Z.$  A normalization of the $\chi$-eigenmeasure may be specified by choosing a right Haar measure on $Q$ and a right Haar measure on $Q_z.$  
By loc.~cit., the quasi-character $\chi$ above is an extension of the modular quasi-character of $Q_z$ times the inverse of the modular quasi-character of $Q.$  In particular, if no such extension exists, then no eigenmeasure exists.  If 
$Q_0$ is an affine algebraic group over $F$ we denote by $\delta_{Q_0}$ the modular quasi-character of $Q_0(F).$

Let $\psi:F\to \CC^\times$ be a nontrivial additive character. Let $dx$ be the unique self-dual Haar measure on $F$ with respect to $\psi$, and let
$
d^\times x:=\zeta(1)|x|^{-1}dx,
$
where $\zeta$ is the usual local zeta function.

\subsection{Function spaces} \label{ssec:fspaces}

If $X$ is a smooth scheme of finite type over $F$ with $X(F) \neq \emptyset$, then $X(F)$ is an analytic manifold. 
For any Hermitian vector bundle $\mathcal{H}$ over $X(F)$ (for example, the trivial bundle $\CC$) we let
\begin{align*}
C_c^\infty(X(F),\mathcal{H})\leq C^\infty(X(F),\mathcal{H})
\end{align*}
be the space of  compactly supported smooth sections and  smooth sections, respectively.  Here and below, when $\mathcal{H}=\CC$ is the trivial bundle we drop it from notation.
When $F$ is non-Archimedean, we set $\mathcal{S}(X(F),\mathcal{H}):=C_c^\infty(X(F),\mathcal{H}).$  When $F$ is Archimedean, we define $\mathcal{S}(X(F),\mathcal{H})$ as in \cite{AG:Nash}.

Let $X$ be a quasi-affine scheme of finite type over a local field $F.$ Now assume that $F$ is Archimedean.  Elazar and Shaviv define a Schwartz space $\mathcal{S}_{\mathrm{ES}}(X(F)),$ whether or not $X$ is smooth.  Briefly one views $X(F)=\mathrm{Res}_{F/\RR}X(\RR)$ as a real algebraic variety, chooses a closed embedding $X(F) \to \RR^n$ in the category of real algebraic varieties and defines $\mathcal{S}_{\mathrm{ES}}(X(F))$ to be the space of restrictions of Schwartz functions on $\mathcal{S}(\RR^n).$  One then endows $\mathcal{S}_{\mathrm{ES}}(X(F))$ with the nuclear Fr\'echet space structure induced by the quotient map $\mathcal{S}(\RR^n) \to \mathcal{S}_{\mathrm{ES}}(X(F))$ given by restriction of functions.  This space and its Fr\'echet space structure are independent of the choice of closed embedding. In the special case that $X$ is smooth, one has $\mathcal{S}(X(F))=\mathcal{S}_{\mathrm{ES}}(X(F)).$ For all of this, we refer to \cite{Elazar:Shaviv}. 

Assume that $X$ admits a right action by an affine algebraic group $G$ over $F$.  Let $\mathcal{U}(\mathfrak{g})$ be the universal enveloping algebra of the complexification of the real algebra $\mathfrak{g}:=\mathrm{Lie}\,\mathrm{Res}_{F/\RR}G(\RR).$

\begin{lem} \label{lem:D:closed}
 For all $D \in \mathcal{U}(\mathfrak{g})$ and $f \in \mathcal{S}_{\mathrm{ES}}(X(F)),$ one has $D.f \in \mathcal{S}_{\mathrm{ES}}(X(F)).$
\end{lem}
\begin{proof}
We can choose a representation $V$ of $G$ and an equivariant immersion $X \to V$ \cite[p. 217, Lemma 2]{Rosenlicht}.  Thus we may assume $X$ is a vector space, in which case the lemma is clear.
\end{proof}

For unity of notation, when $F$ is non-Archimedean, we let $\mathcal{S}_{\mathrm{ES}}(X(F)):=C_c^\infty(X(F))$ be the space of locally  constant compactly supported functions on $X(F).$

\subsection{Global setting} \label{ssec:global}

Let $F$ be a number field. For a place $v$, let $F_v$ be the completion of $F$ at $v$. We often use the symbol $\infty$ to denote the collection of Archimedean places of $F$, and set $F_\infty := \prod_{v\in \infty}F_v$. The idelic norm on $\mathbb{A}_F^\times$ is given by $|\cdot|:=\prod_{v} |\cdot|_v$. For an idelic quasi-character $\eta$, we define $\eta_s$ for $s\in \CC$ and $\mathrm{Re}(\eta)$ as in the local setting. We fix a nontrivial additive character $\psi=\prod_v \psi_v:F\backslash \mathbb{A}_F\to \CC^\times$.  
Using $\psi_v$ to define local measures as in \S \ref{ssec:measure}, we obtain Haar measures $dx:=\prod_v dx_v$ and $d^\times x:=\prod_v d^\times x_v$ on $\mathbb{A}_F$ and $\mathbb{A}_F^\times$ respectively. 

Let $G$ be an affine algebraic $F$-group.  We let $A_G$ be the neutral component of the real points of the maximal $\QQ$-split torus in the center of $\mathrm{Res}_{F/\QQ}G.$  Let $X^*(G)$ be the space of ($F$-rational) characters of $G$, and let 
$$
G(\A_F)^1:=\bigcap_{\chi \in X^*(G)}(|\cdot| \circ \chi:G(\A_F) \to \RR_{>0}).  
$$
If $G$ has reductive neutral component, then  $A_GG(\A_F)^1=G(\A_F).$  We write
\begin{align*}
    [G]:=G(F) \backslash G(\A_F) \quad \textrm{ and } \quad [G]^1:=G(F) \backslash G(\A_F)^1.
\end{align*}

 If $X$ is a quasi-affine scheme of finite type over $F$, the topological space $X(\A_F)$ is defined and has good properties \cite{Conrad_adelic_points}.  If $\mathcal{H}$ is a Hermitian vector bundle over $X(\A_F),$ we form the spaces 
\begin{align*}
    C_c^\infty(X(\A_F),\mathcal{H}) \leq C^\infty(X(\A_F),\mathcal{H}).
\end{align*}
We do not know how to define the adelic Schwartz space in this level of generality.  We will however define
\begin{align*}
    \mathcal{S}_{\mathrm{ES}}(X(\A_F)):=\left(\widehat{\otimes}_{v|\infty}\mathcal{S}_{\mathrm{ES}}(X(F_v)) \right) \otimes \mathcal{S}_{\mathrm{ES}}(X(\A_F^\infty)).
\end{align*}
Here the hat denotes the completed (projective) tensor product; this is unambiguous because $\mathcal{S}_{\mathrm{ES}}(X(F_v))$ is nuclear for $v|\infty.$
Let $S$ be a finite set of places of $F$ including the infinite places and let $\mathcal{X}$ be a flat $\OO_{F}^S$-scheme of finite type with $\mathcal{X}_{F}=X.$  Then $\mathcal{S}_{\mathrm{ES}}(X(\A_F^\infty))=\otimes_{v \nmid \infty}'\mathcal{S}_{\mathrm{ES}}(X(F_v)),$ where the restricted tensor product is with respect to $\one_{\mathcal{X}(\OO_{F_v})}$ for almost all places $v.$

\subsection{Affine $\Psi$-bundles}\label{ssec:gen whit}
 Let $G$ be an affine algebraic group over a characteristic zero field $F$ and let $X$ be an affine scheme of finite type over $F$ equipped with a right $G$-action.  We assume that $X$ admits a unique open $G$-orbit $X^{\circ}.$

Let $p:\mathcal{V}^\circ \to X^\circ$ be a $\GG_a^n$-torsor.   We assume that 
$\mathcal{V}^\circ$ is also equipped with an action of a semidirect product $\GG_a^n \rtimes G$ extending the action of $\GG_a^n$ and that $p$ is $\GG_a^n \rtimes G$-equivariant.  Here $\GG_a^n \rtimes G$ acts via its quotient $G$ on $X^\circ.$  
  We point out that we have an action
  \begin{align} \label{action:on:char:dom} \begin{split}
      \mathcal{V}^{\circ}(R) \times R^n \times G(R) &\lto \mathcal{V}^{\circ}(R) \times R^n\\
      (v,z,g) &\longmapsto (vg,(0\rtimes g)^{-1}z(0 \rtimes g)).\end{split}
  \end{align}
We assume we are given a morphism $
\Psi:\mathcal{V}^{\circ} \times \GG_a^n \lto \GG_a$ such that 
\begin{enumerate}
    \myitem[($\Psi1$)] \label{W1} $\Psi$ factors through the map $p \times \mathrm{Id}:\mathcal{V}^{\circ} \times \GG_a^n \to X^\circ \times \GG_a^n.$
    \myitem[($\Psi2$)] \label{W2} $\Psi$ is $G$-invariant with respect to \eqref{action:on:char:dom}.
    \myitem[($\Psi3$)] \label{W3} For all $w \in  \mathcal{V}^{\circ}(R)$, the map $\Psi(w,\cdot):\GG_{a,R}^n \to \GG_{a,R}$ is a homomorphism of group schemes over $R.$
\end{enumerate}
An \textbf{affine $\Psi$-bundle} is a pair $(p:\mathcal{V}^\circ\to X^\circ,\Psi)$ as above.  We point out that the action of $\GG_a^n \rtimes G$ on $\mathcal{V}^\circ$ is part of the data.  The notion of an affine $\Psi$-bundle is an axiomatization of the constructions in \S \ref{sec:Whittaker:ind} and \S \ref{sec:fiber} below. 
\begin{rem}
A Whittaker induction in the sense of \cite[\S 2.6]{SV} is a pair $(X^L\wedge^{P}G,\Lambda)$ where $P=LU$ is a parabolic subgroup of $G$ with Levi subgroup $L$ and unipotent radical $U,$ $X^L$ is a homogeneous spherical $L$-variety, and
    \[
    \Lam: X^L \lra \Hom(U,\mathbb{G}_a)
    \]
    is an $L$-equivariant map.
     From a Whittaker induction, one may construct an affine $\Psi$-bundle equipped with an action of $[U,U]\backslash U \rtimes (G\times L)$ as follows: set $\mathcal{V}^\circ = [U,U]\backslash G\times X^L$ and $X^\circ = U\backslash G\times X^L$, and let $p:\mathcal{V}^{\circ} \to X^\circ$ be the natural projection. Then $[U,U]\backslash U\simeq \mathbb{G}_a^n$ for some $n.$  We define $\Psi: [U,U]\backslash G\times X^L\times [U,U]\backslash U\to \mathbb{G}_a$ by
    \[
    \Psi(g,x,u) = \Lam(x)(u).
    \]
    The constructions in \S \ref{sec:Whittaker:ind} are essentially of this form, while those of \S \ref{sec:fiber} are more general.
\end{rem}

Let $w\in \mathcal{V}^\circ(F)$ and $x:=p(w).$  The morphism $\Psi$ induces characters
$$
\Psi_w:G_{x} \lto \GG_a
$$
as we now explain.  
The choice of $w\in \mathcal{V}_x^\circ(F)$ induces an isomorphism
\begin{align*}
    R^n &\tilde{\lto} \mathcal{V}_x^\circ(R)\\
    z &\longmapsto wz.
\end{align*}
For $g \in G_x(R),$ let $z_w(g)\in R^n$ be the unique element such that 
\begin{align}\label{zw}
wz_w(g)=wg.
\end{align}
We define
\begin{align*}
\Psi_w(g):=\Psi(w,z_w(g)).    
\end{align*}
\begin{lem} \label{lem:Psi}
    The map $\Psi_w:G_x \to \GG_a$ is independent of the choice of $w \in \mathcal{V}^{\circ}_x(F).$  It is  a character.  
\end{lem}

\begin{proof}
Let $w' \in \mathcal{V}^\circ_x(F).$  Then $w=w'z$ for a unique $z \in F^n.$  It follows that 
$$
z_{w'}(g)=z+z_w(g)-(0 \rtimes g)^{-1}z(0\rtimes g)
$$
and hence, using \ref{W1}
$$
\Psi_{w'}(g)=\Psi(w,z_{w'}(g))=\Psi(w,z+z_w(g)-(0 \rtimes g)^{-1}z(0\rtimes g)).
$$
Using \ref{W2} and \ref{W3}, we deduce that this is $\Psi(w,z_w(g))=\Psi_{w'}(g),$ proving independence.  On the other hand for $g,h \in G_x(R)$ one has 
$$
wz_w(gh)=wgh=wgz_{wg}(h)=
w(z_w(g)+z_{wg}(h))
$$
so 
\begin{align*}
\Psi_w(gh)&=\Psi(w,z_w(gh))=\Psi(w,z_w(g)+z_{wg}(h))=\Psi(w,z_w(g))+\Psi(w,z_{wg}(h)),
\end{align*}
where we have used that $\Psi(w,\cdot)$ is linear by  \ref{W3}. 
Using \ref{W1} and the fact that $\Psi_w$ is independent of the choice of $w \in \mathcal{V}_x^\circ(F)$,
\[
\Psi(w,z_{wg}(h))=\Psi(wg,z_{wg}(h))=\Psi(w,z_{w}(h))=\Psi_w(h),
\] so we deduce
\begin{equation*}
\Psi_w(gh)=\Psi_w(g)+\Psi_w(h).\qedhere
\end{equation*}
\end{proof}

\subsection{Hermitian line bundles} \label{ssec:basic} Let $(p:\mathcal{V}^\circ \to X^\circ,\Psi)$ be an affine $\Psi$-bundle. For a sufficiently large finite set of places $S \supseteq \infty$, we choose models of $\mathcal{V}^\circ,X^\circ,X$ and $G$ over $\OO_F^S;$ we denote them by the same symbols by abuse of notation.  We assume that the models of $\mathcal{V}^{\circ},$ $X^{\circ},$ and $G$ are smooth, and that all the obvious smooth $F$-morphisms involving these schemes (and $\GG_a^n$) extend to the integral models and are again smooth. 

Assume $\psi=\prod_v\psi_v$ is unramified outside of $S.$  There is a unique Hermitian line bundle
$\mathcal{H}$ over $X^\circ(\A_F)$ whose space of smooth sections is
\begin{align} \label{sections:V}
\{f \in C^\infty(\mathcal{V}^\circ(\A_F)): \mathcal{R}(a \rtimes I)f(w)=\psi(\Psi(w,a))f(w) \textrm{ for } a\in \A_F^n\}.
\end{align}
The action of $G(\A_F)$ on $C^\infty(\mathcal{V}^{\circ}(\A_F))$ preserves this subspace by \ref{W2} above.

We can also define Hermitian line bundles $\mathcal{H}_{v}$ over $X(F_v)$ for every place $v$ of $F$ and $\mathcal{H}_{\infty}$ over $X(F_{\infty})$ in the obvious manner.  Then
$$
\mathcal{H}:=\mathcal{H}_{\infty} \otimes \otimes_{v\nmid \infty}'\mathcal{H}_{v},
$$
where the restricted tensor product is with respect to 
$$
\one_{X^\circ(\OO_{F_v}),\psi_v},
$$
the unique sections taking the value $1$ on $\mathcal{V}^\circ(\OO_{F_v})$ for all $v \not \in S.$

Assume that $X^\circ$ admits a nowhere-vanishing section $\Omega$ of the canonical bundle such that $\Omega(xg)=\chi(g)\Omega(x)$ for some $\chi:G \to \GG_m.$   Then $X^\circ(F_v)$ admits a $G(F_v)$-eigenmeasure $dx_v$ for all $v.$ This allows us to form the $L^2$-spaces
$$
L^2(X^{\circ}(F_v),\mathcal{H}_v)
$$
for all $v$ and $L^2(X^\circ(F_\infty),\mathcal{H}_\infty).$

\subsection{Bounding sums by integrals} A standard method from calculus is to estimate the sum over $\ZZ$ of  a smooth function using an integral.  Finis and Lapid \cite{FinisLapid2011} have provided an elegant generalization of this to the adelic context that will prove useful below.  We recall it for the convenience of the reader.  Let $G$ be an affine algebraic group over the number field $F$.  Let $\mathcal{U}(\mathfrak{g}_{F_\infty})$ be the universal enveloping algebra of the complexificiation of $\mathfrak{g}_{F_\infty}:=\mathrm{Lie}\,G(F_\infty)$ (viewed as a real Lie algebra)  and let $K\leq G(\A_F^\infty)$ be a compact open subgroup.  Let $\mathcal{V}_G$ be a basis for the set of elements of $\mathcal{U}(\mathfrak{g}_{F_\infty})$ of degree less than or equal to $[F:\QQ]\dim_F G$ with respect to the usual grading.  Let $dg'=dg'_\infty \otimes dg'^\infty$ be a left Haar measure on $G(\A_F)$ such that $dg'^\infty(K)=1.$  Thus the space $L^1(G(\A_F)/K)$ is defined.  Set
\begin{align*}\mathcal{C}_1(G(\A_F),K):&=\{f:G(\A_F)/K \to \CC : \norm{D.f}_{L^1(G(\A_F)/K)}<\infty \textrm{ for all }D \in \mathcal{U}(\mathfrak{g}_{F_\infty})\}.
\end{align*}
For $f \in \mathcal{C}_1(G(\A_F),K),$ let
\begin{align*}
    \mu(f):&=\sum_{D \in \mathcal{V}_G}\norm{D.f}_{L^1(G(\A_F)/K)}.
\end{align*}
\begin{lem}[Finis and Lapid] \label{lem:FL}
Let $C$ be a compact neighborhood of the identity element of $G(\A_F)$. Then there
exists $c \in \RR_{>0}$ such that for any $f \in  \mathcal{C}_1(G(\A_F), K)$ and any $g \in G(\A_F)$ one has
\begin{align}
    |f(g)| \leq c\sum_{D \in \mathcal{V}_G}\int_{C}|D.f(gg')|dg'.
\end{align}
Moreover,
$$
\sum_{\gamma \in G(F)}|f(g \gamma )| \leq c\mu(f).
$$ 
\end{lem}
\begin{proof}
When $F=\QQ$ this is \cite[Lemma 3.3]{FinisLapid2011}.  Using restriction of scalars, it then follows for general $F.$  
\end{proof}

At various times we will infer a pointwise bound on smooth functions on $\RR_{>0}^n$ using bounds on integrals of the function and their derivatives.  This is possible due to the following identity extracted from the proof of \cite[Lemma 3.3]{FinisLapid2011}:
\begin{lem} \label{lem:Sobolev}
    For smooth $f:\RR^n \to \CC$ one has 
    \begin{align*}
        |f(0)| \leq \sum_{ I \subseteq \{1,\dots,n\}}\int_{[0,1]^n}\left|\frac{\partial^I f}{\prod_{i \in I}\partial x_i}(x) \right|dx.
    \end{align*} \qed
\end{lem}

\section{Groups and affine $\Psi$-bundles} \label{sec:groups:orbits}

Let $F$ be a field of characteristic zero. In this section, we define various groups and schemes that occur in the Poisson summation formulae conjectured in \S \ref{ssec:adelic:S} below.  In particular, we define certain affine $\Psi$-bundles in \S\ref{sec:Whittaker:ind} and \S\ref{sec:fiber} below.
\subsection{Groups and $z$-extensions}\label{sec: z-ext}
We recall our notation from \S \ref{ssec:period}. For $r_1,r_2,r_3 \in \ZZ_{>0}$, let 
\begin{align*}
\GL_{\underline{r}}:=\GL_{r_1} \times \GL_{r_2} \times \GL_{r_3}.
\end{align*}
We will often be treating spaces that come in triples in the current work, and we will use notation such as $\underline{r-2}$ to denote $(r_1-2,r_2-2,r_3-2)$, $\underline{1}$ to denote $(1,1,1),$ etc.  Let $I_{\underline{r}}$ be the identity of $\GL_{\underline{r}}(F)$.   

We denote by  $R$ an $F$-algebra, used to define the points of schemes. 
Let  
\begin{align} \label{H:def}
H(R):=\{(h_1,h_2,h_3) \in \GL_{\underline{2}}(R): \det h_1=\det h_2=\det h_3\} \index{$H$}
\end{align}
and let
\begin{align*}
\nu:H \lto \GG_m \index{$\nu$}
\end{align*}
be given on points by $(h_1,h_2,h_3) \mapsto \det h_1.$

Let $\mathrm{GSp}_{6}$ denote the symplectic similitude group on a $6$-dimensional vector space, and let $P \leq \mathrm{GSp}_{6}$, $M_P \leq P$ denote the usual Siegel parabolic and Levi subgroup.  More specifically, 
\begin{align}\begin{split}\label{eq:levi}
\mathrm{GSp}_{6}(R):&=\left\{g \in \GL_{6}(R): g^t\begin{psmatrix} & I_3 \\ -I_3 & \end{psmatrix}g=\lambda \begin{psmatrix} & I_3 \\ -I_3 & \end{psmatrix} \textrm{ for some } \lambda \in R^\times \right\},\\
M_P(R):&=\{\begin{psmatrix} \lambda A & \\ & A^{-t} \end{psmatrix}: A \in\GL_3(R), \, \lambda \in R^\times\},\\
N_P(R):&=\{\begin{psmatrix} I_3 &  Z\\ & I_3 \end{psmatrix} : Z \in M_{3,3}(R), \, Z^t=Z\}, \end{split}
\end{align}
and $P=M_PN_P.$ 

Let $\mathcal{E} \leq H$ be the finite \'etale group scheme over $F$ whose points in an $F$-algebra $R$ are \index{$\mathcal{E}$}
\begin{align}\label{eqn: E group}
\mathcal{E}(R)=\left\{\left(\begin{psmatrix} \epsilon_1 & 0 \\ 0 &  \epsilon_1 \end{psmatrix}, \begin{psmatrix} \epsilon_2 & 0 \\ 0 & \epsilon_2 \end{psmatrix},\begin{psmatrix} \epsilon_3 & 0 \\ 0&  \epsilon_3 \end{psmatrix}\right) \in Z_{\GL_{\underline{2}}}(R):\epsilon_i^2=1 ,\epsilon_1\epsilon_2\epsilon_3=1 \right\} .
\end{align}

Let $Z_G$ denote the center of an affine algebraic $F$-group $G.$  
If $G \leq \GL_{\underline{2}}$ is any subgroup containing $\mathcal{E},$ we set \index{$G^e$}
\begin{align} \label{Ge}
    G^e:=Z_{\GL_{\underline{2}}} \wedge^{\mathcal{E}} G=(Z_{\GL_{\underline{2}}} \times G)/\mathcal{E}
\end{align} 
Here $\mathcal{E}$ is embedded via $x \mapsto (x,x)=(x,x^{-1}).$ 
Thus $H^e$ is a $z$-extension of  $H /\mathcal{E} $ in the sense of \cite[\S 1]{Kottwitz:rational}. The groups $Z_{\GL_{\underline{2}}}$ and $H$ may be canonically identified with subgroups of $H^e.$
The character $\nu$ induces a character of $H/\mathcal{E}$ and hence $H^e.$  
We observe that $H^e$ admits maps 
\begin{align} \label{pi} \begin{split}
p_1:H^e \lto \GL_{\underline{2}} \index{$p_1$} \quad \textrm{ and }\quad
p_2:H^e \lto H/\mathcal{E}. \index{$p_2$} \end{split} 
\end{align}
Here $p_1$ is induced by the product map $Z_{\GL_{\underline{2}}} \times \GL_{\underline{2}} \to \GL_{\underline{2}}$ and $p_2$ is the canonical map.  We point out that $p_1$ restricted to either $Z_{\GL_{\underline{2}}}$ or $H$ is the identity.

Let
$$
G:=\GL_{\underline{r}} \times \GL_{\underline{r-2}} \times H^e.
$$
We write $g=(g^{(1)},g^{(2)},g^{(3)})$ for the projections to the three factors.

\subsection{A quotient of a symplectic group}  \label{ssec:geo:setup}

For an $F$-scheme $Y,$ we let
$\overline{Y}^{\mathrm{aff}}:=\mathrm{Spec}\, \Gamma(Y,\OO_Y).$
This is the affine closure of $Y.$  Let us record the following useful lemma:

\begin{lem} \label{lem:normal} 
If $Y$ is  locally noetherian, normal and connected, then $\overline{Y}^{\mathrm{aff}}$ is normal.
\end{lem}
\begin{proof}
Since $Y$ is normal, $\Gamma(Y,\OO_Y)$ is normal which implies that  $ \overline{Y}^{\mathrm{aff}}$ is normal \cite[Lemma 6.38]{Gortz_Wedhorn}.
\end{proof}

A quasi-affine scheme $Y$ of finite type over $F$ is said to be \textbf{strongly quasi-affine} if the algebra of global sections $\Gamma(Y,\OO_Y)$ is finitely generated and the natural map 
$
Y \to \overline{Y}^{\mathrm{aff}}
$ is an open immersion.

We recall a homogeneous space for a particular similitude group. Let
$$
X_P^{\circ}:=P^{\mathrm{der}}\backslash \mathrm{Sp}_6.
$$
It is known that $X_P^{\circ}$ is strongly quasi-affine \cite[Theorem 1.1.2]{Braverman:Gaitsgory}.  By Lemma \ref{lem:normal} $X_P:=\overline{X_P^{\circ}}^{\mathrm{aff}}$ \index{$X_P$} is normal. Consider the Pl\"ucker map
\begin{align*} 
\mathrm{Pl}_P:P^{\mathrm{der}} \backslash \mathrm{Sp}_6 \lto \displaystyle\wedge^3 \GG_a^6
\end{align*}
induced from the map $\mathrm{Sp}_6(R) \to \wedge^3 R^6$ given by
\begin{align*}
    \begin{psmatrix} * \\ v_1 \\ v_2 \\ v_3 \end{psmatrix} \longmapsto v_1 \wedge v_2 \wedge v_3
\end{align*}
where $v_1,v_2,v_3$ are the last three rows of an element of $\mathrm{Sp}_6(R).$ The Pl\"ucker embedding defines a closed immersion of $X_P$ into $\wedge^3 \GG_a^6,$ and we have
\begin{align*}
\mathrm{Pl}_P(X_P)=\{0\} \cup \mathrm{Pl}_P(X_P^{\circ})
\end{align*}
\cite[Theorem 1 and 2]{Popov:Vinberg}.
Therefore, we write $\{0\}$ for the point that is $X_P-X_P^{\circ}$. It is a $\mathrm{Sp}_6$-orbit. 

We define a homomorphism
\begin{align*}
\begin{split}
m:\mathrm{GSp}_6(R) &\lto M_P(R)\\
    g &\longmapsto \begin{psmatrix}1 & & & \\
                  & \nu(g)^{-1}I_3 && \\  & & & I_2
        \end{psmatrix},  \end{split}
\end{align*}
and let $\mathrm{GSp}_6$ act on $X_P^{\circ}$ via
\begin{align} \label{GSp:act}
\begin{split}
X_P^{\circ}(R) \times \mathrm{GSp}_6(R) &\lto X_P^{\circ}(R)\\
    (x,g) &\longmapsto m(g) xg.
\end{split}
\end{align}
It is easy to check from the description of $X_P$ in \cite[Theorem 1.1.2]{Braverman:Gaitsgory} that this extends to an action on $X_P$ sending $\{0\}$ to itself.
We denote this action using a dot:
\begin{align*} 
x.g:=m(g)xg.
\end{align*}
We endow $\wedge^3\GG_a^6$ with the action of $\mathrm{GSp}_6$ given on pure tensors by 
$$
(v_1 \wedge v_2 \wedge v_3,g) \mapsto \nu(g)^{-1} v_1g \wedge v_2g \wedge v_3g.
$$
Then $\mathrm{Pl}_P$ is $\mathrm{GSp}_6$-equivariant.

We now link the discussion above to the group $H$ defined in \eqref{H:def}.  One has an embedding
\begin{align*}
\begin{split}
    H(R) &\lto \mathrm{GSp}_6(R)\\
    \left(\begin{psmatrix} a_i & b_i \\ c_i & d_i \end{psmatrix} \right) &\longmapsto \begin{psmatrix} a_1 & & &  b_1 & & &\\
    & a_2 & & & b_2 &\\
    && a_3 & & & b_3 \\ c_1 & & &  d_1&& \\
    &c_2 & & & d_2 & \\ & & c_3 & & & d_3\end{psmatrix}.
    \end{split}
\end{align*}
Let
\begin{align*} \begin{split}
\gamma_0:=\left(\begin{smallmatrix} 0 & 0 & 0 & -1 & 0& 0\\ 0 & 1 & 0 & 0 & 0 &0\\
0 & 0 & 1 & 0& 0 & 0\\
1 & 1 & 1 & 0 & 0 & 0\\ 0 & 0 & 0 &-1 & 1 &  0\\
0 & 0 & 0 & -1 & 0 & 1  \end{smallmatrix}\right). \end{split}
\end{align*}
We identify $\gamma_0$ with its image under $\mathrm{Sp}_6(F) \to X_P^{\circ}(F)$ when convenient.
 For every $\gamma \in X_P^{\circ}(F)$, let $H_\gamma$ denote the stabilizer of $\gamma$ under the $H$-action, and let $O(\gamma)$ be the orbit of $\gamma$.   We use analogous notation for other group actions.
One has
\begin{align} \label{commute}
    \begin{psmatrix} 1 & & \\ &\lambda^{-1}I_3 & \\ & &I_2 \end{psmatrix}\gamma_0 =\gamma_0 \begin{psmatrix} \lambda^{-1}I_3 & \\ & I_3 \end{psmatrix}
\end{align}
for any $\lambda \in R^\times.$  This is why we defined \eqref{GSp:act} as above. 

Let
\begin{align*}
    \Delta:\GG_a \lto \GG_a^{\underline{1}}
\end{align*}
denote the diagonal embedding and let \index{$T_H$} \index{$N_0$}
\begin{align} \begin{split} 
T_H(R):&=\left\{\begin{psmatrix} \Delta(\lambda) & \\ & 1 \end{psmatrix}: \lambda \in R^\times  \right\},\\
\label{N0:def} N_0(R):&=\left\{\left(\begin{psmatrix} 1 & t_1 \\ & 1 \end{psmatrix},\begin{psmatrix} 1 & t_2 \\ & 1 \end{psmatrix},\begin{psmatrix} 1 & t_3 \\ & 1 \end{psmatrix}\right):t_1+t_2+t_3=0 \right\}. \end{split}
\end{align}
It follows from \cite[Lemma 2.3]{Getz:Liu:Triple} that
\begin{align}\label{H:stabilizer}
\begin{split}
H_{\gamma_0}(R)&=T_H(R)N_0(R).
\end{split}
\end{align}

\begin{lem} \label{lem:extend}
The action of $\SL_{\underline{2}}$ on $X_P$ induces an isomorphism  $N_0 \backslash \SL_{\underline{2}} \tilde{\to} O(\gamma_0)$ that extends to an isomorphism $\overline{N_0 \backslash \mathrm{SL}_{\underline{2}}}^{\mathrm{aff}} \to X_P.$
\end{lem}
\begin{proof}
We have an isomorphism $N_0 \backslash \SL_{\underline{2}} \tilde{\to} O(\gamma_0)$ by the \eqref{H:stabilizer}.  Moreover, by the stabilizer computations in \cite[Lemma 2.3]{Getz:Liu:Triple}
the codimension of $X_P-O(\gamma_0)$ in $X_P$ is at least $2.$ Thus we conclude by the algebraic Hartogs theorem \cite[Theorem 6.45]{Gortz_Wedhorn}. 
\end{proof}

In view of \eqref{commute} and Lemma \ref{lem:extend}, the action of $\mathrm{GSp}_6$ on $X_P$ is the extension of the action of $H$ on $N_0 \backslash \SL_{\underline{2}}$ induced by 
\begin{align} \label{H:act:0} \begin{split}
    \mathrm{SL}_{\underline{2}}(R) \times H(R) &\lto \mathrm{SL}_{\underline{2}}(R)\\
    (g,h) &\longmapsto \begin{psmatrix} \Delta(\nu(h))^{-1} & \\ & 1 \end{psmatrix} g h. \end{split}
\end{align}

There are three morphisms $\mathrm{Pl}_i:O(\gamma_0) \to \GG_a^2$ induced by the morphisms
\begin{align*}
\SL_{\underline{2}}(R) &\lto R^2\\
\left(\begin{psmatrix} a_i & b_i \\ c_i & d_i \end{psmatrix}\right) & \longmapsto \begin{psmatrix} c_i & d_i \end{psmatrix}.
\end{align*}
We sometimes write 
\begin{align} \label{cidi}
    \mathrm{Pl}_i(x)=(c_i(x),d_i(x)). \index{$\mathrm{Pl}_i$}
\end{align}

\subsection{The Pl\"ucker embedding} \label{sec:the:pairing}
The Pl\"ucker embedding $\mathrm{Pl}_P:X_P \to \wedge^3 \GG_a^6$ \index{$\mathrm{Pl}_P$}
has image in the irreducible subrepresentation $V_P < \wedge^3\GG_a^6$ of $\mathrm{Sp}_6$ spanned by the highest weight vector attached to the Siegel parabolic $P$ as in \cite[\S 3.2]{Getz:Hsu:Leslie}. The representation $V_P$ may be described as the kernel of the map 
\begin{align} \label{contraction}
\wedge^3 \GG_a^6 \lto \GG_a^6
\end{align}
given by contraction with the element of the dual of $\wedge^2 \GG_a^6(F)$ corresponding 
to the symplectic pairing defining $\mathrm{Sp}_6$ \cite[p. 258, \S 17.1]{Fulton:Harris}. 
Concretely, by \cite[p. 260, \S 17.2]{Fulton:Harris}, the map \eqref{contraction} is given on pure tensors by
\begin{align*}
v_1 \wedge v_2 \wedge v_3 &\longmapsto \sum_{i<j}\langle v_i,v_j\rangle(-1)^{i+j-1} v_{k}
\end{align*}
where $\langle\,,\,\rangle$ is the pairing defining $\mathrm{Sp}_{6}$ and $k$ is the unique index that is not $i$ or $j.$  For convenience, we write $e_{ij}=e_i \wedge e_j \in \wedge^2\ZZ^6$ and $e_{ijk}:=e_i \wedge e_j \wedge e_k \in \wedge^3 \ZZ^6.$ 
The following is a basis of $V_P(F):$
   \begin{align*}
        e_{123}, e_{126}, e_{135}, e_{156}, e_{234}, e_{246}, e_{345}, e_{456},e_1',e_2',e_3',e_4',e_5',e_6'
        \end{align*}
        where
        \begin{align*}
        e'_{1}:=e_{124}+e_{236}, e'_{2}:=e_{125}-e_{136}, e'_{3}:=e_{134}-e_{235}, e_4':=e_{145}+e_{356}, e'_5:=e_{245}-e_{346}, e'_6:=e_{146}-e_{256}.
\end{align*}
We denote by $V_3$ the span of the first $8$ vectors and by $V_P'$ the span of the last $6$ vectors.  
There is a symplectic pairing 
\begin{align} \label{VP} \begin{split}
\langle\, ,\rangle_P:V_P(R) \times V_P(R) &\lto \wedge^6 R \tilde{\lto}R\\
(v,v') &\longmapsto v \wedge v' \end{split}
\end{align}
where the isomorphism is given by stipulating that $e_1 \wedge e_2 \wedge e_3 \wedge e_4 \wedge e_5 \wedge e_6$ is sent to $1.$  We observe that for $h \in H(R)$ one has that 
\begin{align} \label{h:equiv}
    \langle v,v'.h^{-1}\rangle_P= \langle v,\nu(h)v'h^{-1}\rangle_P=\langle \nu(h)^{-2}vh,v'\rangle_P=\langle \nu(h)^{-1}v.h,v'\rangle_P.
\end{align}

The subspaces $V_3$ and $V_P'$ are orthogonal with respect to $\langle\,,\,\rangle_P.$  The matrix of $\langle\,,\,\rangle_P|_{V_3}$ with respect to the basis above is
\begin{align} \label{prod1}
    \begin{psmatrix} & & & & & & & 1 \\& & & & & & -1 &  \\ & & & & &-1 & & \\& &  & & 1 & & &
    \\& &  & -1 & & &&
    \\& & 1 & & & && \\& 1 &  & & & && \\ -1 & & & & & & &\end{psmatrix}.
\end{align}  

For $\left(\begin{psmatrix} a_i& b_i \\ c_i & d_i \end{psmatrix}\right) \in \SL_{\underline{2}}(R)$ we have 
$$
\gamma_0\begin{psmatrix} a_1 & & &b_1 & & \\ & a_2 & && b_2 & \\ & & a_3 & & &b_3\\ c_1 & & &d_1 & & \\ & c_2 & & &d_2 & \\ & & c_3 & & &d_3 \end{psmatrix}=\begin{psmatrix} *  \\
\begin{smallmatrix} a_1 & a_2 & a_3 & b_1 & b_2 & b_3\\ -c_1 & c_2 & 0 & -d_1 & d_2 & 0\\ -c_1 & 0 & c_3 & -d_1 & 0 & d_3 \end{smallmatrix}\end{psmatrix}.
$$
In terms of the coordinates for $V_P$ given above, this is

\begin{center}
\begin{tabular}{|l|l|}
\hline
$e_{123}$ & $a_1c_2c_3+c_1a_2c_3+c_1c_2a_3$    \\ \hline
$e_{126}$ & $a_1c_2d_3+c_1c_2b_3+c_1a_2d_3$    \\ \hline
$e_{135}$ & $-(a_1d_2c_3+c_1b_2c_3+c_1d_2a_3)$ \\ \hline
$e_{156}$ & $c_1d_2b_3+a_1d_2d_3+c_1b_2d_3$    \\ \hline
$e_{234}$ & $d_1a_2c_3+d_1c_2a_3+b_1c_2c_3$    \\ \hline
$e_{246}$ & $-(d_1a_2d_3+b_1c_2d_3+d_1c_2b_3)$ \\ \hline
$e_{345}$ & $b_1d_2c_3+d_1b_2c_3+d_1d_2a_3$    \\ \hline
$e_{456}$ & $b_1d_2d_3+d_1d_2b_3+d_1b_2d_3$    \\ \hline
$e_1'$ & $-c_2$                             \\ \hline
$e_2'$ & $-c_1$                             \\ \hline
$e_3'$ & $c_3$                              \\ \hline
$e_4'$ & $d_2$                              \\ \hline
$e_5'$ & $d_1$                             \\ \hline
$e_6'$ & $-d_3$                              \\ \hline
\end{tabular}
\end{center}
\smallskip
This computation implies that $V_P' \cong \GG_a^2 \oplus \GG_a^2 \oplus \GG_a^2$ as representations of $\SL_{\underline{2}}.$  
Since $V_P'$ is $H$-invariant, the same is true of $V_3.$
It is not hard to see that $V_3 \cong \GG_a^2 \otimes \GG_a^2 \otimes \GG_a^2.$

For $(a,t) \in (R^\times)^{\underline{1}} \times R^{\underline{1}},$ we have
\begin{align} \label{pl:comp}
\mathrm{Pl}_P \left(\gamma_0\begin{psmatrix}1/a & at\Delta([a^{-1}]) \\ & a \end{psmatrix} \right)&= \tfrac{[a]}{a^2_1}e_{156}-\tfrac{[a]}{a^2_2}e_{246}+\tfrac{[a]}{a^2_3}e_{345}+ \sum_{i=1}^3t_ie_{456}+a_2e_4'+a_1e_5'-a_3e_6'
\end{align}
where $[a]:=a_1a_2a_3.$

\index{$v_0$}
Let $v_0:=e_{156}-e_{246}+e_{345}.$  This is the image of $\gamma_0$ under the $H$-equivariant map
\begin{align} \label{X2V3}
    X_P \stackrel{\mathrm{Pl}_P}{\lto} V_P \lto V_3,
\end{align}
where the second map is the canonical projection. Let $O(v_0)$ be the $H$-orbit of $v_0$ in $V_3.$
Recall from \eqref{H:act:0} that  we have an action of $H$ on $O(\gamma_0).$  Let $\mathcal{E}$ be defined as in \eqref{eqn: E group}.

\begin{lem} \label{lem:XtoV3}
The orbit maps $H \to O(\gamma_0)$ and $H \to O(v_0)$ induce $H$-equivariant isomorphisms $N_0 \backslash \SL_{\underline{2}} \tilde{\to} O(\gamma_0)$ and $\mathcal{E} N_0\backslash \SL_{\underline{2}} \tilde{\to}O(v_0).$  One has a diagram of $H$-equivariant maps
\begin{equation}
    \begin{tikzcd}
 N_0 \backslash \SL_{\underline{2}} \arrow[r] \arrow[d]& X_P \arrow[d]\\
    \mathcal{E}N_0 \backslash \SL_{\underline{2}} \arrow[r]& V_3,
    \end{tikzcd}
\end{equation}
    where the left vertical arrow is the canonical map and the right vertical arrow is \eqref{X2V3}.
\end{lem}

\begin{proof}
We have already shown that there is an isomorphism $N_0 \backslash \SL_{\underline{2}} \tilde{\to} O(\gamma_0)$ that extends to $\overline{N_0 \backslash \SL_{\underline{2}}}^{\mathrm{aff}} \tilde{\to} X_P$ in Lemma \ref{lem:extend}.

There is a natural action of $S_6,$ the symmetric group on $6$ letters, on the standard basis of $F^6,$ and hence an action of $S_6$ on $\wedge^3 F^6.$  Under this action, $v_0$ is fixed by $(123)(456).$  Hence $(\SL_{\underline{2}})_{v_0}$ is stabilized by the automorphism given on points by $(g_1,g_2,g_3) \longmapsto (g_2,g_3,g_1).$

 The stabilizer $(\SL_{\underline{2}})_{v_0}$ contains $N_0.$
We claim that the stabilizer $(\SL_{\underline{2}})_{v_0} \leq B,$ where $B$ is the Borel subgroup of upper triangular matrices in $\SL_{\underline{2}}.$  
If not, then by the comments in the previous paragraph, the $R$-points of the stabilizer contains an element of one of the following forms:
\begin{align} 
  \label{case1} &\left( \begin{psmatrix} c_1t & * \\ c_1 & * \end{psmatrix},\begin{psmatrix} c_2t & * \\ c_2 & * \end{psmatrix},\begin{psmatrix} c_3t & * \\ c_3 & * \end{psmatrix}\right),
\\   
 \label{case2}&\left( \begin{psmatrix} c_1t & * \\ c_1 & * \end{psmatrix},\begin{psmatrix} c_2t & * \\ c_2 & * \end{psmatrix},\begin{psmatrix} a_3 & 0 \\ 0 & a_3^{-1} \end{psmatrix}\right),\\
 \label{case3}&\left( \begin{psmatrix} c_1t & * \\ c_1 & * \end{psmatrix},\begin{psmatrix} a_2 & 0 \\ 0 & a_2^{-1} \end{psmatrix},\begin{psmatrix} a_3 & 0 \\ 0 & a_3^{-1} \end{psmatrix}\right).
\end{align}
Here $c_i \in R^\times$ and we have used the fact that $N_0$ is contained in the stabilizer. 

We use the table above to analyze these cases. If we consider the action on $v_0$ by an element of the form \eqref{case1},  we conclude that $t=0$ by considering the coordinate $e_{123}$; considering the coordinate $e_{126}$, we see that this element cannot stabilize $v_0$.  
Similarly, we see that elements of the form \eqref{case2} and \eqref{case3}  cannot be in the stabilizer  by considering the coordinate $e_{123}$ and $e_{126}$ respectively.

Thus the stabilizer of $v_0$ in $\SL_{\underline{2}}$ is equal to the stabilizer of $v_0$ in $B.$  Using \eqref{pl:comp} we see that this stabilizer is $\mathcal{E}N_0.$
Thus we obtain an isomorphism $
   \mathcal{E} N_0 \backslash \SL_{\underline{2}} \tilde{\to} O(v_0).$   The commutativity of the diagram in the lemma is clear.
\end{proof}

The action of $H$ on $V_3$ factors through $\mathcal{E}.$  Letting $Z_{\GL_{\underline{2}}}$ act trivially, we obtain an action of $H^e$ on $V_3.$  The stabilizer $H_{v_0}$ of $v_0$ in $H$ is $T_H \mathcal{E} N_0$ by Lemma \ref{lem:XtoV3} and hence the stabilizer of $v_0$ in $H^e$ is 
\begin{align} \label{Htildestab}
H^e_{v_0}= Z_{\GL_{\underline{2}}} \wedge^{\mathcal{E}} \mathcal{E}T_H N_0\cong Z_{\GL_{\underline{2}}} \times T_HN_0.
\end{align}

While the map $H \to O(v_0)$ is surjective as a map of schemes, the map $H(F) \to O(v_0)(F)$ is not surjective in general. On the other hand, we have the following lemma:
\begin{lem} \label{lem:surj}
The map $H^e(F) \to O(v_0)(F)$ is surjective, as is the restriction $\mathrm{SL}_{\underline{2}}^e(F) \to O(v_0)(F).$
\end{lem}

\begin{proof}
 Since $H^1(F,Z_{\GL_{\underline{2}}} \times T_HN_0)=H^1(F,Z_{\GL_{\underline{2}}} \times N_0)=1$, the lemma follows from the basic exact sequence in Galois cohomology.
\end{proof}

\noindent In fact, 
Lemma \ref{lem:surj} is the motivation for introducing $H^e$ in the first place.

\subsection{The space $Y$}
\label{sec:Y:def}
Let $\beta \in F^\times$ and let $\xi=\begin{psmatrix}1 & \\ & \Delta(\beta)\end{psmatrix}.$\index{$\beta$} \index{$\xi$}
Let
\begin{align} \label{Y:tilde}
\widetilde{Y}(R):&=\left\{ ((u_1,u_2),x) \in \mathcal{M}(R) \times O(\gamma_0)(R):  u_1c(x)+\Delta(\beta) u_2d(x)=0\right\}.
\end{align}
The last equality here is shorthand for $u_{i1}c_i(x)+\beta u_{i2}d_i(x)=0$ for $1 \leq i \leq 3$, where we are using notation as in \eqref{cidi}

The group $\mathcal{E}$ acts freely on $O(\gamma_0),$ hence on $\mathcal{M} \times O(\gamma_0),$ and this action preserves  $\widetilde{Y}.$  The space $\widetilde{Y}$ is quasi-affine.  Thus the fppf quotient
$
\mathcal{E} \backslash \widetilde{Y}$
 is a scheme, and $\widetilde{Y} \to \mathcal{E} \backslash \widetilde{Y}$ is an $\mathcal{E}$-torsor \cite[Tag 07S7]{stacks-project}.
By Lemma \ref{lem:XtoV3} there is a canonical isomorphism $\mathcal{E}\backslash O(\ga_0)\tilde{\to} O(v_0)\subset V_3.$  This induces an immersion $\mathcal{E} \backslash \widetilde{Y}\to \mathcal{M} \times V_3.$  Identifying $\mathcal{E} \backslash \widetilde{Y}$ with its image, we set \index{$Y$}
\begin{align} \label{Y:def}    
Y:=\overline{\mathcal{E} \backslash \widetilde{Y}} \subset \mathcal{M} \times V_3.
\end{align}

For $h \in \GL_{2}(R)$ or $h \in \GL_{\underline{2}}(R)$ we set\index{$h^\iota$}
\begin{align}
h^{\iota}:=\xi h^{-t} \xi^{-1},
\end{align} 
where we abuse notation and set $\xi=\begin{psmatrix}1 & \\ & \beta\end{psmatrix}$ when $h\in \GL_2(R)$. This automorphism preserves $Z_{\GL_{\underline{2}}},H,\mathcal{E}$ and hence induces an automorphism of $H^e.$

The group $\GL_{\underline{r-2}} \times H^e$ acts on $\mathcal{M} \times V_3$ via
\begin{align} \label{R:def} \begin{split}
  \mathcal{R}:  \mathcal{M}(R) \times V_3(R) \times \GL_{\underline{r-2}}(R) \times H^e(R) &\lto \mathcal{M}(R) \times V_3(R)\\
    ((m,v),g,h) &\longmapsto (g^{t}mp_1(h)^{\iota},v.p_2(h)). \end{split}
\end{align}
The action preserves $\mathcal{M}^{\circ} \times V_3,$ as well as $Y\subset \mathcal{M}\times V_3$. Identifying $\mathcal{M}^{\circ}$ with $\mathcal{M}^{\circ} \times \{0\},$ we obtain an action of $\GL_{\underline{r-2}} \times H^e$ on $\mathcal{M}^{\circ}$ as well.

We now explain how to view $Y$ as a vector bundle over $\mathcal{M}^{\circ}.$
Let \index{$m_0$}
\begin{align}
    m_0:=(e_{\underline{r-2}},0)
\end{align}
and let
\begin{align} \label{V:def}
V(R):=\{ae_{156}+be_{246}+ce_{345}+te_{456} :a,b,c,t \in R\};
\end{align}
note that $v_0\in V(F)$.

\begin{lem} \label{lem:Y:vb}
The projection to the first factor  realizes $Y_{\mathcal{M}^{\circ}}=Y\times_{\mathcal{M}}\mathcal{M^\circ}$ as a $\GL_{\underline{r-2}} \times H^e$-equivariant vector bundle over $\mathcal{M}^\circ.$ 
One has 
\begin{align} \label{Y:fiber}
Y_{m_0}=\{m_0\} \times V.
 \end{align}
\end{lem}

\begin{proof}
We prove the second statement first. We have
\begin{align} \label{fiber}
\widetilde{Y}_{m_0}(R)=\left\{\left(m_0,\gamma_0\begin{psmatrix}1/a & a\Delta(t)[a^{-1}] \\ & a \end{psmatrix} \right) :(a,t) \in  (R^\times)^3 \times R\right\}.
\end{align}
The computation \eqref{pl:comp} implies that the image of  $\widetilde{Y}_{m_0}$ in $\{m_0\} \times V_3$ is a dense subset of $\{m_0\} \times V.$  The subscheme $\{m_0\} \times V \subset \{m_0\} \times V_3$ is closed; hence $\{m_0\} \times V$ is the closure of the image of $\widetilde{Y}_{m_0}.$
Since $Y_{m_0}$ is a closed subset of $\{m_0\} \times V_3$ containing the image of $\widetilde{Y}_{m_0},$ we deduce that $\{m_0\} \times V \subseteq Y_{m_0}.$ On the other hand, the dimension of $Y_{m_0}$
is equal to the dimension of $\widetilde{Y}_{m_0},$ which is $4$ by \eqref{fiber}.  Thus $Y_{m_0}=\{m_0\} \times V.$

Temporarily let
\begin{align}
Q(R):=\left\{\left(\begin{psmatrix} I_{\underline{r-3}} & \\ v & 1 \end{psmatrix},\begin{psmatrix} a^{-1} & \\ b & a \end{psmatrix}\right)\in \GL_{\underline{r-2}}(F)\times \SL_{\underline{2}}(F):(v,a,b) \in R^{\underline{r-3}} \times (R^\times)^{\underline{1}} \times R^{\underline{1}}\right)\}.
\end{align}
For suitable $g_1,\dots,g_k \in \GL_{\underline{r-2}}(F) \times \SL_{\underline{2}}(F),$ one has a Zariski open cover $
\bigcup_{i=1}^k m_0Qg_i$ of $\mathcal{M}^{\circ},$ and each $m_0Qg_i$ is a $Q$-torsor (over $F$).  This open cover provides us with local trivializations
$$
m_0Qg_i \times V \tilde{\lto} Y_{m_0Qg_i}.
$$
\end{proof}

\begin{rem}
One can build zeta integrals analogous to those in \eqref{glob:zeta} with $Y$ replaced by $\widetilde{Y}$  and $H^e$ replaced by $H.$  This was in fact our original construction, and it is in some sense simpler.  However, the space $\widetilde{Y}$ is not a vector bundle over $\mathcal{M}^\circ,$ and we do not know how to  construct Fourier transforms and prove Poisson summation formula for the fibers of  $\widetilde{Y}_{\mathcal{M}^\circ} \to \mathcal{M}^\circ.$  In contrast, it is known how to construct Fourier transforms and prove the Poisson summation formulae for the fibers of any vector bundle, including $Y_{\mathcal{M}^\circ} \to \mathcal{M}^\circ.$
\end{rem}
\quash{
Let $J:=\left(\begin{psmatrix} &1\\  -1& \end{psmatrix},\begin{psmatrix} &1\\ -1&  \end{psmatrix},\begin{psmatrix} &1\\  -1 & \end{psmatrix}\right) \in \SL_2^3(\ZZ)$ and set
\begin{align} \label{J:pair}
    \langle\,,\,\rangle:=\langle \cdot,\cdot.  J\xi^{-1} \rangle_P: V_3 \times V_3 \lto \GG_a.
\end{align}

\begin{lem}
    The pairing \eqref{J:pair} is symmetric.
\end{lem}
\begin{proof}
Using the basis for $V_3(F)$ used to compute \eqref{prod1}, this is straightforward to check. 
\end{proof}

Using \eqref{h:equiv}, for $h \in H(R)$ one has 
\begin{align}  \label{h:equiv2} \begin{split}
\langle v,v'.h^{-1}\rangle&=\langle v,v'.h^{-1}. J\xi^{-1} \rangle_P=\langle v'.J\xi^{-1}.\nu(h)^{-1}h^{-\iota}\rangle_P\\
&=\langle \nu(\nu(h)^{-1}h^{-\iota})v.(\nu(h)h^{\iota}),v'\rangle\\
&=\langle v.h^{\iota},v'\rangle.
\end{split}
\end{align}}

There is a morphism
\begin{align} 
  \mathrm{pr}: \mathcal{M}^{\circ} \lto \mathbb{P}^1
\end{align}
that assigns to $(m_1,m_2) \in \mathcal{M}^{\circ}(R)$ the  morphism 
\begin{align*}
R^2 &\lto Rm_1+Rm_2 <R^{r-2}\\
   (v_1,v_2) &\longmapsto v_1m_1+v_2m_2.
\end{align*}
The image is a projective rank $1$-submodule of $R^{r-2}$ by 
\cite[Lemma 16.17]{Gortz_Wedhorn}.
 We therefore obtain a map \index{$\mathrm{pr}_i$}
\begin{align}  \label{Plperp}
  \mathrm{pr}_i: \mathcal{M}_i^{\circ} \lto \mathbb{P}^1.
\end{align}
At the level of $F$-points, this map is described as follows: for $(m_1,m_2) \in \mathcal{M}_i^{\circ}(F),$ there is an $x \in F^{r_i-2}$ and $(c_1,c_2) \in F^2-\{(0,0)\}$ such that $(m_1,m_2)=(c_1x,c_2x).$  Then $\mathrm{pr}_i(m_1,m_2)$ is the line spanned by $(c_1,c_2).$  
 Set \index{$\mathrm{pr}$}
\begin{align} \label{Plbrak}
\mathrm{pr}:=\prod_{i=1}^3\mathrm{pr}_i:\mathcal{M}^{\circ} \lto \mathbb{P}^{\underline{1}} .
\end{align}
 For $m \in \mathcal{M}^\circ(R),$ the free $R$-submodule $Y_{m}(R)<V_3(R)$ of rank $4$ depends only on $\mathrm{pr}(m).$  Define
\begin{align} \label{Yell}
Y_{\mathrm{pr}(m)}:=Y_m.
\end{align}

\subsection{Actions} \label{ssec:actions}

We now prepare to compute the open orbit of $\GL_{\underline{r-2}} \times H^e$ in $Y$ and the generic stabilizer with respect to this action.

Let $B_2\leq \GL_2$ be the Borel subgroup of upper triangular matrices and let \index{$B$}
\begin{align} \label{B:def}
    B:=B_{\underline{2}} \cap \SL_{\underline{2}}, \quad T:= T_{\underline{2}} \cap \SL_{\underline{2}}
\end{align}
Thus $B^e=Z_{\GL_{\underline{2}}} \wedge^{\mathcal{E}}B$ is a Borel subgroup of $\SL_{\underline{2}}^e.$
Let $B_H<H$ be the Borel subgroup of upper triangular matrices.  

When working with zeta integrals later, it will be important to have an understanding of the $B^e$-orbit of $(m_0,v_0)$ at the level of points.    
We have a map 
\begin{align} \label{orbit}
B^e \lto B_{\underline{2}} \times V_3
\end{align}
given on points in an $F$-algebra $R$ by
\begin{align*}
    B^e(R) &\lto B_{\underline{2}}(R) \times  V_3(R)\\
    h &\longmapsto  (p_1(h),v_0.p_2(h)).
\end{align*}
The group $B^e$ acts on $B^e$ on the right, and on $B_{\underline{2}} \times V_3$ via 
\begin{align*}
(b,v).h=(bp_1(h),v.p_2(h)).
\end{align*}
The map \eqref{orbit} is $B^e$-equivariant.

On the other hand, we have a  closed subscheme $J \subset B_{\underline{2}} \times V$ given by
\begin{align*} 
J(R):&=\left\{\left( \begin{psmatrix} x & y \\ & w 
    \end{psmatrix},c_1e_{156}-c_2e_{246}+c_3e_{345}+te_{456}\right) :  x_iw_i^{-1}c_{i+1}c_{i+2}=1, \sum_{i=1}^3w_i^{-1}y_i[c]=t\right\}.
    \end{align*}
Here  $(x,y,w,c,t) \in (R^\times)^{\underline{1}} \times R^{\underline{1}} \times (R^\times)^{\underline{1}} \times (R^\times)^{\underline{1}} \times R$ and we are regarding indices modulo $3$ in the obvious sense.  One has
\begin{align}\label{J}    
   &J(R)=\left\{ \left(\begin{psmatrix}1 & y \\ & 1 \end{psmatrix}\left(\begin{psmatrix} x_1 &  \\ & x_1c_2c_3
    \end{psmatrix},\begin{psmatrix} x_2 &  \\ & x_2c_3c_1
    \end{psmatrix},\begin{psmatrix} x_3 &  \\ & x_3c_1c_2 
\end{psmatrix}\right),\left(c_1,c_2,c_3,[c]\sum_{i=1}^3y_i\right)\right)\right\}
\end{align}
where $x,c \in (R^\times)^{\underline{1}}, y \in R^{\underline{1}}$ and we have written elements of $V\subset V_3$ as tuples $(a,b,c,d)$ with respect to the basis $e_{156},-e_{246},e_{345},e_{456}.$

\begin{lem} \label{lem:Bisom}
The map \eqref{orbit} induces an isomorphism $B^e \to J.$
\end{lem}
\begin{proof}
By generalities on algebraic group actions \cite[Chapter 7]{Milne:AGbook}, it suffices to check that $J(\overline{F})$ is the $B^e(\overline{F})$-orbit of $(I_{\underline{2}},v_0)$ and the stabilizer of this point is trivial.  
The pullback of \eqref{orbit} to $Z_{\GL_{\underline{2}}} \times B$ is given on points in an $F$-algebra $R$ by 
\begin{align}  \label{induced} \begin{split}
Z_{\GL_{\underline{2}}}(R) \times B(R) &\lto B_{\underline{2}}(R) \times V_3(R)\\
\left(z, \begin{psmatrix} 1/a & at\Delta([a^{-1}])\\ & a \end{psmatrix}\right) &\longmapsto \left(z \begin{psmatrix} 1/a & at\Delta([a^{-1}])\\ & a \end{psmatrix},\tfrac{[a]}{a^2_1}e_{156}-\tfrac{[a]}{a^2_2}e_{246}+\tfrac{[a]}{a^2_3}e_{345}+ \sum_{i=1}^3t_ie_{456}\right). \end{split}
\end{align}
Here we have used \eqref{pl:comp}.  The fiber over $(I_{\underline{2}},v_0)$ is $\mathcal{E}$ (embedded diagonally) and hence the fiber of \eqref{orbit} over $(I_{\underline{2}},v_0)$ is trivial.  Thus it suffices to show that the map \eqref{induced} has image in $J$ and is surjective at the level of $\overline{F}$-points. 
This can be checked directly.
\end{proof}

For $x,c \in (F^\times)^{\underline{1}},$ let 
\begin{align} \label{ell:def}
    \ell\begin{psmatrix} x&\\ & x\Delta([c])/c\end{psmatrix} \in \SL_{\underline{2}}^e(F)
\end{align}
be an element satisfying $p_1(\ell\begin{psmatrix} x&\\ & x\Delta([c])/c\end{psmatrix})=\begin{psmatrix} x&\\ & x\Delta([c])/c\end{psmatrix}.$  The $\ell$ stands for ``lift.''

Let $(\GL_{\underline{r-2}} \times H^e)_{y_0}$ be the stabilizer of $y_0:=(m_0,v_0),$ and let \index{$y_0$}
\begin{align} \label{Ycirc}
    Y^{\circ}:=O(y_0) \subset Y
\end{align}
be the orbit of $y_0.$

\begin{lem} \label{lem:orb:stab1}
The map
$$
\mathrm{Id} \times p_1:(\GL_{\underline{r-2}} \times H^e)_{y_0} \lto \GL_{\underline{r-2}} \times \GL_{\underline{2}}
$$
is an isomorphism onto its image, and the set of $R$-points of the image is 
\begin{align*}
   \left\{\left(\begin{psmatrix}g & z \\ & \Delta(\lambda)x \end{psmatrix},\begin{psmatrix} \Delta(\lambda)x& \Delta(\lambda)xy \\ & x  \end{psmatrix}\right):
   (g,z,x,\lambda,y) \in \GL_{\underline{r-3}}(R) \times R^{\underline{r-3}} \times (R^\times)^{\underline{1}} \times R^\times \times  R^{\underline{1}}, \sum_{i=1}^3 y_i=0\right\}.
\end{align*} 
The unique open orbit of $\GL_{\underline{r-2}} \times H^e$ on $Y$ under $\mathcal{R}$ is $Y^\circ.$
\end{lem}

\begin{proof}
By \eqref{Htildestab} we have
\begin{align} \label{contained}
    (\GL_{\underline{r-2}} \times H^e)_{y_0} \leq \GL_{\underline{r-2}} \times (Z \wedge^{\mathcal{E}}\mathcal{E}T_HN_0).
\end{align}
Since the restriction of $p_1$ to $Z \wedge^{\mathcal{E}}\mathcal{E}T_HN_0$ is an isomorphism onto its image, we deduce that the same is true of the morphism in the lemma.  The characterization of the image follows from \eqref{contained}.

We deduce that the dimension of $Y^\circ=O(y_0)$ is equal to the dimension of $Y.$  
We claim that $Y$ is irreducible. 
First, $\mathcal{M}^{\circ}$ is irreducible since any orbit map $\GL_{\underline{r-2}} \times \SL_{\underline{2}} \to \mathcal{M}^{\circ}$ induced by a choice of base point is surjective.  
The vector bundle $Y_{\mathcal{M}^\circ}$ is then irreducible by \cite[Tag 004Z]{stacks-project} because $Y_{\mathcal{M}^{\circ}} \to \mathcal{M}^{\circ}$ is smooth, hence open. 
 Finally, $Y$ is the closure of $Y_{\mathcal{M}^{\circ}},$ hence irreducible \cite[Tag 004U]{stacks-project}.  Since $Y^\circ$ is an orbit under a group action the underlying topological space of $Y^\circ$ is a locally closed subset of $Y.$  On the other hand $Y^\circ \subset Y$ is of codimension $0$ and $Y$ is irreducible; hence $Y^\circ$ is open and dense  in $Y.$   The other orbits are in the closure of $Y^\circ$ and hence of strictly smaller dimension. 
\end{proof}

The morphisms $\mathrm{Pl}_i:\SL_{\underline{2}} \to \GG_a^2$ do not descend to $\mathcal{E}N_0 \backslash \SL_{\underline{2}}.$  However,
$$
\mathrm{Pl}_i^{\otimes 2}:=\mathrm{Pl}_i \otimes \mathrm{Pl}_i:\SL_{\underline{2}} \lto \GG_a^2 \otimes \GG_a^2
$$
does descend to $\mathcal{E}N_0 \backslash \SL_{\underline{2}}.$  Since $\mathcal{E}N_0 \backslash \SL_{\underline{2}} \tilde{\lto} O(v_0)$ by Lemma \ref{lem:XtoV3} we obtain maps \index{$\mathrm{Pl}_i^{\otimes 2}$}
\begin{align} \label{Pli2}
\mathrm{Pl}_i^{\otimes 2}:O(v_0) \tilde{\lto}\mathcal{E}N_0 \backslash \SL_{\underline{2}} \lto \GG_a^2 \otimes \GG_a^2 .
\end{align}\quash{
For each $i,$ we have an action 
\begin{align}
\mathcal{R}_i:\GG_a^2 \otimes \GG_a^2 \times H \lto \GG_a^2 \otimes \GG_a^2
\end{align}
given on pure tensors by $v \otimes w\mathcal{R}_i((h_1,h_2,h_3))= vh_i \otimes wh_i.$  This action descends to $H/\mathcal{E},$ and hence pulls back to $H^e.$  Let us continue to denote the induced action of $H^e$ by $\mathcal{R}_i.$
Using  \eqref{H:act:0}, for $(v,h) \in O(v_0)(R) \times H^e(R)$ we have 
\begin{align} \label{Pliequi}
\mathrm{Pl}_i^{\otimes 2}(v.p_2(h))=\mathrm{Pl}_i^{\otimes 2}(v)\mathcal{R}_i(h).
\end{align}}
Using \eqref{induced} we see that 
\begin{align} \label{Pli2:comp}
\mathrm{Pl}_i^{\otimes 2}\left(c_1e_{156}-c_2e_{246}+c_3e_{345}+ te_{456} \right)= \frac{[c]}{c_i}e_i \otimes e_i.
\end{align}
\quash{

The map
\begin{align}
\mathrm{br}:=\mathrm{Pl}_1 \otimes \mathrm{Pl}_2 \otimes \mathrm{Pl}_3:\SL_{\underline{2}} \lto \GG_a^2 \otimes \GG_a^2 \otimes \GG_a^2
\end{align}
also descends to $\mathcal{E}N_0 \backslash \SL_{\underline{2}}$ and hence defines a map 
\begin{align} \label{br}
\mathrm{br}:O(v_0) \lto \GG_a^2 \otimes \GG_a^2 \otimes \GG_a^2.
\end{align}
The br is short for ``bottom row.''
Using \eqref{induced} we see that 
\begin{align}
\mathrm{br}\left(c_1e_{156}-c_2e_{246}+c_3e_{345}+ te_{456} \right)=c_1c_2c_3e_2 \otimes e_2 \otimes e_2.
\end{align}}


\subsection{Some affine $\Psi$-bundles} \label{sec:Whittaker:ind}
In this subsection we will only consider constructions related to a single factor $\GL_{r_i}$ of $\GL_{\underline{r}}.$  We therefore omit the $i$ from notation.

For $F$-algebras $R$  let \begin{align} \label{Prdef}
P_{r-2,2}(R):=\left\{\begin{psmatrix} g & x \\ & h \end{psmatrix},(g,h,x) \in \GL_{r-2}(R) \times \GL_2(R) \times M_{r-2,2}(R)\right\}.
\end{align}
Thus $P_{r-2,2}$ is a maximal parabolic subgroup of $\GL_{r}.$  Let $N_{r-2,2} \leq P_{r-2,2}$ be the unipotent radical.  We identify the Levi subgroup of block diagonal matrices in $P_{r-2,2}$ with $\GL_{r-2} \times \GL_2.$  This provides us with an identification $P_{r-2,2}=N_{r-2,2} \rtimes (\GL_{r-2} \times \GL_2).$

Consider the quotient 
\begin{align} \label{Xr}
    X_{r}^\circ:=N_{r-2,2} \backslash \GL_{r}.
\end{align}
We let $X_r:=\overline{X_r^\circ}^{\mathrm{aff}}$ be the affine closure of $X_r^\circ.$  The action 
\begin{align} \label{GL:act}\begin{split}
    \GL_{r}(R) \times \GL_r(R) \times  \GL_{r-2}(R) \times \GL_2(R) \lto \GL_{r}(R)\\
    (g_0,g,g',h) \longmapsto \begin{psmatrix}g' & \\ & h \end{psmatrix}^{-1}g_0g \end{split}
\end{align}
descends to $X_{r}^\circ.$  With respect to this action $X_{r}^\circ$ is a spherical variety.  
Let
\begin{align}\label{Mipairing}
   \begin{split}\langle \,,\, \rangle:M_{r-2,2}(R) \times M_{r-2,2}(R) &\lto R\\
    \left((u_1,u_2),(u_1',u_2')\right) &\longmapsto 
    u_1^tu_1'+\beta u_2^tu_2'.
    \end{split}
\end{align}
For each $m \in M_{r-2,2}(R),$ we use the obvious isomorphism $N_{r-2,2} \tilde{\to} M_{r-2,2}$ to regard $\langle m,\cdot\rangle$ as a homomorphism $\langle m,\cdot \rangle:(N_{r-2,2})_{R} \to \GG_{aR}.$ 
The group $\GL_{r-2} \times \GL_2$ acts on $N_{r-2,2}$ by conjugation.  
For each $m \in \mathcal{M}(F),$ let $(\GL_{r-2} \times \GL_2)_m$ be the stabilizer of the character $\langle m,\cdot \rangle:N_{r-2,2} \to \GG_a.$  Thus
\begin{align*}
    (\GL_{r-2} \times \GL_2)_m(R)=\{(g',h) \in \GL_{r-2}(R) \times \GL_2(R): g'^{t}mh^{\iota}=m\}.
\end{align*}\quash{
This group is the stabilizer of $m$ with respect to the action
\begin{align} \label{M:action} \begin{split}
\mathcal{M}(R) \times \GL_{r-2}(R) \times \GL_2(R) \lto \mathcal{M}(R)\\
(m,g,h) \longmapsto g^tmh^{\iota}. \end{split}
\end{align}}

We have an action of $\GL_r \times P_{r-2,2}$ on $\GL_{r} \times \mathcal{M}$ given on points in an $F$-algebra $R$ by
\begin{align} \label{act} \begin{split}
 \GL_{r}(R) \times  \mathcal{M}(R) \times \GL_{r}(R) \times P_{r-2,2}(R) &\lto  \GL_{r}(R) \times  \mathcal{M}(R)\\
    \left(g_0,m,g,\begin{psmatrix} I_{r-2} & z \\ & I_2 \end{psmatrix}\begin{psmatrix}g' & \\ & h \end{psmatrix}\right) &\longmapsto \left( \begin{psmatrix} g' & \\ & h \end{psmatrix}^{-1} \begin{psmatrix} I_{r-2} & z \\ & I_2 \end{psmatrix}^{-1}g_0g,g^tmh^{\iota}\right). \end{split}
\end{align}
The canonical map $\GL_r \times \mathcal{M} \to X^{\circ}_r \times \mathcal{M}$ is equivariant, where $\GL_{r} \times P_{r-2,2}$ acts through its Levi quotient $\GL_{r} \times \GL_{r-2} \times \GL_{2}$ on the codomain.  

There is a morphism
\begin{align}
\Psi_{r}:\GL_r \times \mathcal{M} \times N_{r-2,2} \lto \GG_a
\end{align}
given on points by $\left(g,m,\begin{psmatrix} I_r & z\\ & I_2 \end{psmatrix} \right) \mapsto \langle m,z\rangle.$
The pair $(\GL_r \times \mathcal{M} \to X^{\circ}_r \times \mathcal{M},\Psi_r)$ is an affine $\Psi$-bundle in the sense of \S \ref{ssec:gen whit}.

To study functions on $X_r^\circ,$ it is helpful to construct an affine embedding of it.
For 
$$
\begin{psmatrix} v_1\\\vdots \\ v_{r-2}\end{psmatrix} \in M_{r-2,r}(R) \textrm{ and }v_{r-1},v_{r} \in R^r
$$
we abbreviate 
$$
\begin{psmatrix} v_1\\\vdots \\ v_{r-2}\end{psmatrix} \wedge v_{r-1} \wedge v_r:=\begin{psmatrix} v_1 \wedge v_{r-1} \wedge v_r\\
\vdots \\v_{r-2} \wedge v_{r-1} \wedge v_r\end{psmatrix} \in M_{r-2,r}(R) \wedge R^r \wedge R^r \leq M_{r-2,1}(\wedge^3 R^r).
$$

We have a map
\begin{align*} \begin{split}
\mathrm{Pl}_{X_r}:\GL_{r}(R) &\lto M_{r-2,r}(R) \wedge R^{r} \wedge R^{r}  \times   M_{2,r}(R) \\
\begin{psmatrix} v_1\\
\vdots\\
v_r \end{psmatrix} &\longmapsto \left(\begin{psmatrix} v_1\\\vdots \\ v_{r-2}\end{psmatrix} \wedge v_{r-1} \wedge v_r,\begin{psmatrix}v_{r-1}\\  v_{r} \end{psmatrix} \right)\end{split}
\end{align*}
that descends to a morphism
\begin{align*} 
   \mathrm{Pl}_{X_r} :X_r^\circ \lto M_{r-2,r} \wedge \GG_a^{r} \wedge \GG_a^{r}  \times   M_{2,r}.
\end{align*}
The map is equivariant under $\GL_r \times \GL_{r-2} \times \GL_2.$  Here the action on the codomain is a representation of $\GL_r \times \GL_{r-2} \times \GL_2$ on $M_{r-2,r} \wedge \GG_a^r \wedge\GG_a^r \times M_{2,r}.$  It is given explicitly on pure tensors by  
\begin{align}
\left(\begin{psmatrix} v_1\\\vdots \\ v_{r-2}\end{psmatrix} \wedge v_{r-1} \wedge v_r,\begin{psmatrix}v'_{r-1}\\  v'_{r}\end{psmatrix}\right)(g,g',h):=\left(\left(\frac{g'^{-1}}{\det h}\begin{psmatrix} v_1 g\\\vdots \\ v_{r-2}g\end{psmatrix}\right) \wedge v_{r-1}g \wedge v_rg,h^{-1}\begin{psmatrix}v'_{r-1}g\\  v'_{r}g \end{psmatrix}\right).
\end{align}

\begin{lem} \label{lem:Pl:Iso}
    The map $\det \times \mathrm{Pl}_{X_r}:X_r^\circ \to \GG_m \times M_{r-2,r} \wedge \GG_a^r \wedge \GG_a^r \times M_{2,r}$ induces an isomorphism between $X_r$ and the closure of $\det \times \mathrm{Pl}_{X_r^\circ}(X_r^\circ)$ in $\GG_m \times M_{r-2,1} \wedge \GG_a^r \wedge \GG_a^r \times M_{2,r}.$
\end{lem}

\begin{proof}
It is enough to show that the restriction
\begin{align} \label{is:immersion}
\mathrm{Pl}_{X_r}|_{N_{r-2,2} \backslash \SL_r}:N_{r-2,2} \backslash \SL_r \lto M_{r-2,r} \wedge \GG_a^r \wedge \GG_a^r \times M_{2,r}
\end{align}
induces an isomorphism of the affine closure $\overline{N_{r-2,2} \backslash \SL_r}$ with the closure of $\mathrm{Pl}_{X_r}(N_{r-2,2} \backslash \SL_r)$ in $M_{r-2,r} \wedge \GG_a^r \wedge \GG_a^r \times M_{2,r}.$
For our use in a moment, we point out that \eqref{is:immersion} is the inclusion of an orbit of $\SL_r,$ hence is an immersion by \cite[Proposition 7.17]{Milne:AGbook}.

Let $L$ be the Levi quotient of $P_{r-2,2} \cap \SL_r.$  We have an equivariant closed immersion of $L \times \SL_r$-representations
\begin{align*}
&M_{r-2,r} \wedge \GG_a^r \wedge \GG_a^r \oplus M_{2,r}\lto\mathrm{Hom}((\wedge^3 V_{\mathrm{st}})^{N_{r-2,2}},\wedge^3V_{\mathrm{st}}) \oplus \mathrm{Hom}(V_{\mathrm{st}}^{N_{r-2,2}},V_{\mathrm{st}})
\end{align*}
where $V_{\mathrm{st}}$ is the standard representation, viewed as a space with $\mathrm{SL}_{r}$-action on the right.  Since \eqref{is:immersion} is an $L \times \SL_r$-equivariant immersion, we deduce the lemma from 
the discussion after Definition 1 in  \cite[\S 3.2]{Arzhantsev:Timashev}. 
\end{proof}

We continue to denote by $\det \times \mathrm{Pl}_{X_r}$ the extension of $\det \times \mathrm{Pl}_{X_r} $ to $X_r.$

\subsection{Products} \label{sec:fiber}

We now return to studying $\GL_{\underline{r}}=\GL_{r_1} \times \GL_{r_2} \times \GL_{r_3}.$
Let 
\begin{align}
   X^\circ_{\underline{r}}:=\prod_{i=1}^3X^\circ_{r_i}, \quad X_{\underline{r}}:=\prod_{i=1}^3X_{r_i}.
   \end{align}
Let
 $$
 X:=X_{\underline{r}} \times Y.
 $$
The actions \eqref{act} induce an action of $\GL_{\underline{r}} \times P_{\underline{r-2},\underline{2}}$ on $\GL_{\underline{r}} \times \mathcal{M}$ that descends to $X^{\circ}_{\underline{r}} \times \mathcal{M}$ and extends to $X_{\underline{r}} \times \mathcal{M}.$   Moreover, the canonical map $\GL_{\underline{r}} \times \mathcal{M} \to X^{\circ}_{\underline{r}} \times \mathcal{M}$ is equivariant.
We have an induced action
\begin{align} \begin{split} \label{ind:action}
   ( \GL_{\underline{r}} \times Y)(R) \times \left(\GL_{\underline{r}} \times N_{\underline{r-2},\underline{2}} \rtimes (\GL_{\underline{r-2}} \times H^e)\right)(R) &\lto (\GL_{\underline{r}} \times Y)(R) \\
    \left(g_0,(m,v),g,\begin{psmatrix} I_{\underline{r-2}} & z\\ & I_{\underline{2}} \end{psmatrix} \rtimes (g',h)\right) &\longmapsto \left( \begin{psmatrix} g'^{-1} & g'^{-1}z\\ & p_1(h)^{-1}\end{psmatrix} g_0g,(g'^tmp_1(h)^\iota,v.p_2(h))\right). \end{split}
\end{align}
Here the semi-direct product structure is induced by the morphism 
$$
\mathrm{Id} \times p_1:\GL_{\underline{r-2}} \times H^e \lto \GL_{\underline{r-2}} \times \GL_{\underline{2}}
$$
where the codomain is viewed as our fixed Levi subgroup of $P_{\underline{r-2},\underline{2}}.$  This action descends to $X^{\circ}_{\underline{r}} \times Y$ and extends to $X,$ and the canonical map $\GL_{\underline{r}} \times Y \to X^{\circ}_{\underline{r}} \times Y$ is equivariant.  

For $F$-algebras $R$ we have a morphism
\begin{align} \label{Psi} \begin{split}
\Psi:\GL_{\underline{r}}(R) \times Y(R) \times N_{\underline{r-2},\underline{2}}(R) &\lto R\\
\left(\left(g_0,(m,v)\right),\begin{psmatrix} I_{\underline{r-2}} & z\\ & I_{\underline{2}}\end{psmatrix}\right) &\longmapsto \langle m,z\rangle.\end{split}
\end{align}
It is $G$-invariant.  
Here the relevant action is given by 
\begin{align*}
  &\left(\left(g_0,(m,v)\right),\begin{psmatrix} I_{\underline{r-2}} & z\\ & I_{\underline{2}}\end{psmatrix},g,g',h\right) \\ &\longmapsto \left( \begin{psmatrix} g & \\ & p_1(h)\end{psmatrix}^{-1}g_0g,(g'^{t}mp_1(h)^\iota,v.p_2(h)),\begin{psmatrix} I_{\underline{r-2}} & g'^{-1}zp_1(h)\\ & I_{\underline{r-2}} \end{psmatrix}\right).
\end{align*}
The pair $(\GL_{\underline{r}} \times Y^\circ\to X^\circ,\Psi)$ is thus an affine $\Psi$-bundle in the sense of \S\ref{ssec:gen whit} with $\mathbb{G}_a^n=N_{\underline{r-2},\underline{2}}$ and $G=\GL_{\underline{r}} \times \GL_{\underline{r-2}} \times H^e$.

Using Lemma \ref{lem:orb:stab1} we obtain the following
\begin{lem} \label{lem:orb} \index{$x_0$}
    The scheme $X^\circ$ is the $G$-orbit of $x_0:=(I_{\underline{r}},y_0).$ \qed
\end{lem}
 Let us compute the stabilizer $G_{x_0}$ of $x_0$. Let $Q \leq \GL_{\underline{r}}$ be the subgroup whose points in an $F$-algebra $R$ are given by 
\begin{align} \label{Q:def} \begin{split}
    Q(R):&=\left\{\begin{psmatrix} g & z_1& z_2  & z_3 \\ & \Delta(\lambda) &x_1 &x_2\\ & & \Delta(\lambda) & y\\ & & & 1\end{psmatrix}:
    \sum_{i=1}^3y_i=0 \right\}.\end{split}
\end{align}
Here 
$(g,z_1,z_2,z_3,x_1,x_2,y,\lambda) \in \GL_{\underline{r-3}}(R) \times (R^{\underline{r-3}})^3 \times (R^{\underline{1}})^2 \times R^{\underline{1}} \times R^\times.$

Lemma \ref{lem:orb:stab1} also implies the following:
\begin{lem} \label{lem:orb:stab2} Projection onto the first factor induces an isomorphism 
    \begin{align}
p_1:G_{x_{0}}\lto Z_{\GL_{\underline{r}}}Q.
\end{align} 
We further have isomorphisms
$$
\begin{tikzcd}
 &  \arrow[ld,swap,"\sim"]Z_{G_{x_0}} \arrow[d,"\sim"]\arrow[rd,"\sim"]&\\
    Z_{\GL_{\underline{r}}} & Z_{\GL_{\underline{r-2}}} & Z_{\GL_{\underline{2}}}
\end{tikzcd}
$$
given by projection to the three factors. \qed
\end{lem}

 Lemma \ref{lem:orb:stab2} and Hilbert's theorem 90 for the additive group and general linear group imply that $X^\circ(F)$ is the $G(F)$-orbit of $(I_{\underline{r}},y_0).$  

\subsection{Models over the integers} \label{ssec:models}

Assume $F$ is a number field with ring of integers $\OO_F.$  Let $S \supset \infty$ be a finite set of places of $F.$  Let $W$ be a quasi-affine scheme of finite type over $F.$  Our convention is that a \textbf{model of $W$ over $\OO_F^S$} is a quasi-affine scheme flat and of finite type over $\OO_F^S$ whose generic fiber is $W.$  
Many of the schemes and group schemes defined above admit natural models over $\OO_F$ that we will denote with the same letters by abuse of notation, e.g. $\GL_r,N_{r,r-2},P_{r-2,2},H,Z_{\GL_{\underline{2}}},\mathcal{E}.$  For some of the constructions, this is not so obvious and hence we explicate in this subsection.

The following lemma must be well-known. We include a proof for lack of a reference:

\begin{lem} \label{lem:representability}
If $G$ is a smooth affine group scheme over a Dedekind domain $\OO$ and $G' \leq G$ is a smooth closed affine subgroup, then the fppf quotient sheaf $G' \backslash G$ is representable by a smooth scheme over $\OO.$ It is a group scheme if $G'$ is normal in $G.$
\end{lem}
\begin{proof}
By \cite[Th\'eor\`eme 4.C]{Anantharaman} the quotient $G' \backslash G$ is representable.  
The fppf quotient sheaf is the quotient of $G \times G$ by the equivalence relation $Z$ given on points in an $\OO$-algebra by $Z(R):=\{(g'g,g):(g,g') \in G(R) \times G'(R)\}.$  
We point out that the morphism $Z \to G$ given by projection to the first factor is smooth and surjective. 
Moreover, the property of being smooth and surjective is preserved under base change \cite[tag 02WE]{stacks-project} and is fpqc local \cite[tag 02YJ]{stacks-project}.   
   This implies that $G \to G' \backslash G$ is also smooth and surjective \cite[\S 4]{Raynaud:Passage}.  Since $G$ is smooth over $\OO$ by assumption, we deduce $G' \backslash G$ is smooth \cite[tag 02K5]{stacks-project}.
 To prove the last assertion, one checks that
 if $G$ is a group scheme that is fppf over an arbitrary base and $G' \subset G$ is a normal subgroup scheme, also fppf over the base, then the fppf sheaf $G' \backslash G$ is a sheaf of groups.
\end{proof}

Thus $H^e:=Z_{\GL_{\underline{2}}}\wedge^{\mathcal{E}}H$ and $\SL_{\underline{2}}^e:=Z_{\GL_{\underline{2}}} \wedge^{\mathcal{E}} \SL_{\underline{2}}$ are smooth affine group schemes of finite type over over $\OO_F.$  We extend $\mathcal{M}$ to a scheme over $\OO_F$ by taking the schematic closure in $M_{\underline{r-2},\underline{2}}.$ We then define 
$$
\mathcal{M}^\circ:=\mathcal{M} \cap \prod_{i=1}^3\left(M_{r_i-2,2}-\{0\}\right).
$$
The morphism $\mathrm{pr}$ extends to a morphism $\mathrm{pr}:\mathcal{M}^\circ \to \mathbb{P}^{\underline{1}}$ over $\OO_F.$ 
  We extend $Y$ to a scheme over $\OO_F$ by taking the schematic closure in $\mathcal{M} \times V_3.$  Then the action of $\GL_{\underline{r-2}} \times H^e$ on $Y$ over $F$ extends over $\OO_F.$

The affine scheme $N_{\underline{r-2},\underline{2}} \backslash \GL_{\underline{r}}$ over $\OO_F$ is a model of $X_{\underline{r}}^\circ$ over $\OO_F$ (which we continue to denote by $X_{\underline{r}}^\circ$). We extend the affine closure $X_{\underline{r}}$ to a scheme over $\OO_F$ by taking the schematic closure under the image of the Pl\"{u}cker embedding.   By Lemma \ref{lem:orb:stab1} the stabilizer $(\GL_{\underline{r-2}} \times H^e)_{y_0}$ admits an obvious smooth model over $\OO_F.$  Hence
$(\GL_{\underline{r-2}} \times H^e)_{y_0} \backslash (\GL_{\underline{r-2}} \times H^e)$ is representable over $\OO_F$ by Lemma \ref{lem:representability}.  We may identify $Y^{ \circ}$ with the generic fiber $((\GL_{\underline{r-2}} \times H^e)_{y_0^{}} \backslash (\GL_{\underline{r-2}} \times H^e))_F$  using the action map and thereby regard $(\GL_{\underline{r-2}} \times H^e)_{y_0} \backslash (\GL_{\underline{r-2}} \times H^e)$
as a model of $Y^{\circ }$ over $\OO_F.$  
Likewise $\mathcal{E}N_0\backslash \SL_{\underline{2}}$ provides a model for $O(v_0)$ over $\OO_F.$   Moreover the morphisms $\mathrm{Pl}_i^{\otimes 2}$ of \eqref{Pli2} extend over $\OO_F.$

We point out that for each finite place $v$ of $F,$ the map $\GL_{\underline{r}}(\OO_{F_v})\to X_{\underline{r}}^\circ(\OO_{F_v})$ is surjective
because $N_{\underline{r-2},\underline{2}}$ is smooth over $\OO_{F_v}$ with connected fibers \cite[Lemma B.2.4]{Getz:Hahn}.  Similarly, for finite $v,$ the orbit map $\GL_{\underline{r-2}}(\OO_{F_v}) \times H^e(\OO_{F_v}) \to Y^{\circ}(\OO_{F_v})$ induced by the choice of basepoint $y_0$ is surjective because $(\GL_{\underline{r-2}} \times H^e)_{y_0^{}}$ is smooth over $\OO_{F}$ with connected fibers.

In the sequel we will usually not be explicit about whether we are viewing the schemes discussed in this section over $F$ or over $\OO_F.$  It should be clear from the context.

\section{Schwartz spaces and transforms}
\label{sec:SS}
After some preliminaries, in \S\ref{ssec:desi:loc:ss} we formulate explicit desiderata for certain local Schwartz spaces.  In \S \ref{ssec:adelic:S} we state the Poisson summation conjecture in the setting corresponding to the zeta integrals discussed in the introduction.

\subsection{Eigenmeasures} \label{ssec:vol:forms}

Until later in this subsection we fix a place $v$ of $F$ and omit it from notation, writing $F:=F_v.$
For every $r\in \ZZ_{>0}$, we choose a maximal compact subgroup $K_r$ of $\GL_r(F)$, which we take to be $\GL_r(\mathcal{O}_F)$ if $F$ is non-Archimedean. Let $B_r$ be the Borel subgroup of upper triangular matrices in $\GL_r$, let $T_r \leq \GL_r$ be the maximal torus of diagonal matrices, and let $U_r$ be the unipotent radical of $B_r.$ We call a parabolic subgroup of $\GL_r$ standard if it contains $B_r.$

We normalize the Haar measure on $\GL_r(F)$ so that for $(n,t,k) \in U_{r}(F) \times T_r(F) \times K_r$ one has
$$
d(ntk)=\delta_{B_r}(t)^{-1}\prod_{1 \leq i <j \leq r}dn_{ij} \prod_{i=1}^rd^\times t_{ii} dk,
$$
where $dk$ is the Haar measure on $K_r$ with $dk(K_r)=1$, and $(\cdot)_{ij}$ denotes the $(i,j)$th entry of a matrix.  

We let $d_{X_{r_i}},$ $d_{X_{\underline{r}}}$ be the unique invariant positive Radon measures on $X^\circ_{r_i}(F)$ and $X_{\underline{r}}^\circ(F)$ such that the usual decompositions induced by the Iwasawa decomposition hold, namely
\begin{align*}
d_{X_{r_i}}(pk)=d_\ell \dot{p}dk, \quad (p,k) \in P_{r_i-2,2}(F) \times K_{r_i}\\
d_{X_{\underline{r}}}(pk)=d_\ell \dot{p} dk\quad (p,k) \in P_{\underline{r-2},\underline{2}}(F) \times K_{\underline{r}}.
\end{align*}
Here $d_\ell \dot{p}$ is a quotient of a left Haar measure $d_{\ell}p$ on $P_{r_i-2,2}(F)$ (resp.~$P_{\underline{r-2},\underline{2}}(F)$) by a Haar measure $dn$ on $N_{r_i-2,2}(F)$ (resp.~$N_{\underline{r-2},\underline{2}}(F)$).
Then 
\begin{align*} 
d_{X_{r_i}}\left(\begin{psmatrix}g' & \\ & h \end{psmatrix}^{-1}xg \right)=\frac{|\det g'|^2}{|\det h|^{r_i-2}}d_{X_{r_i}}(x)
\end{align*}
for $(g,g',h) \in \GL_{r_i}(F)\times \GL_{r_i-2}(F) \times \GL_2(F).$

Choose a maximal compact subgroup $K_{H^e} \leq H^e(F)$ such that $K_{H^e}=H^e(\OO_{F})$ if $F$ is non-Archimedean.  We choose left Haar measures $d_\ell b$ on $B_{\underline{2}}^e(F)$ and Haar measures $dh$ on $H^e(F)$  such that the decomposition
$
dh=d_\ell b dk
$
holds.  Here $dk$ is a Haar measure on $K_{H^e}.$

Using Lemma \ref{lem:orb:stab1} we see that the character $\mathrm{Ad}_{(\GL_{\underline{r-2}} \times H^e)_{y_0}}$ admits an extension 
\begin{align} \label{chi} \index{$\chi$}
\chi(g',h):=\frac{\nu(h)^2\prod_{i=1}^3\det g_i'}{\prod_{i=1}^3(\nu(h)\det p_1(h)_i)^{(r_i-2)/2}}
\end{align}
to $\GL_{\underline{r-2}} \times H^e.$  Here $p_1(h)=(p_1(h)_1,p_1(h)_2,p_1(h)_3).$
Thus upon choosing a right Haar measure on $(\GL_{\underline{r-2}} \times H^e)_{y_0}(F),$ we obtain a $|\chi|$-invariant measure on $(\GL_{\underline{r-2}} \times H^e)_{y_0}(F) \backslash \GL_{\underline{r-2}}(F) \times H^e(F).$  The action $\mathcal{R}$ of \eqref{R:def} yields an isomorphism of analytic manifolds
$$
(\GL_{\underline{r-2}} \times H^e)_{y_0}(F) \backslash (\GL_{\underline{r-2}}(F) \times H^e(F)) \lto Y^{\circ}(F).
$$ 
Using this isomorphism, we obtain a $|\chi|$-invariant measure $d_Y$ on $Y^{\circ}(F).$ That is, for $(g',h) \in \GL_{\underline{r-2}}(F) \times H^e(F)$ we have \index{$d_Y$}
\begin{align} \label{omegay:inv} \begin{split}
d_Y(g'^tmp_1(h)^{\iota},v.p_2(h))&=|\chi(g',h)|d_Y(m,v).\end{split}
\end{align}

Let $K_0$ be a maximal compact subgroup of $\SL_{\underline{2}}^e(F);$ we assume $K_0=\SL_{\underline{2}}^e(\OO_{F})$ when $v$ is non-Archimedean.  
Let $B:=B_{\underline{2}} \cap \SL_{\underline{2}}$ as in \eqref{B:def}.
Using Lemma \ref{lem:orb:stab1}, we see that there is a surjection
\begin{align*}
\begin{split}
  K_{\underline{r-2}} \times N_0(F) \backslash B^e(F) \times K_{0} &\lto Y^{\circ}(F)\\
    (k',b,k) &\longmapsto (k'^tm_0p_1(bk)^{\iota},v_0.p_2(bk)) \end{split}
\end{align*}
with compact fibers. One has
\begin{align} \label{dyiwa}
d_{Y}(k'^tm_0p_1(bk)^{\iota},v_0.p_2(bk))=\frac{dk'd_\ell \dot{b} dk}{|\det p_1(b)|^{\underline{(r-2)/2}}}.
\end{align}
Here $d_\ell \dot{b}$ is the quotient of a left Haar measure on $B^e(F)$ by a Haar measure on $N_0(F).$

\subsection{The Harish-Chandra space} \label{ssec:lin:bund}  
There is a unique $G(F)$-equivariant Hermitian line bundle $\mathcal{H}$ over $X^\circ(F)$ whose smooth sections are smooth functions $f$ on $\GL_{\underline{r}}(F) \times Y^{ \circ}(F)$ such that 
\begin{align*}
&  f\left(\begin{psmatrix} I_{\underline{r-2}} & z \\ & I_{\underline{2}} \end{psmatrix}g,m,v\right)=\psi(\langle m,z \rangle)f(g,m,v)
\end{align*} 
for 
$$
(g,z,(m,v)) \in \GL_{\underline{r}}(F) \times  M_{\underline{r-2},\underline{2}}(F) \times Y^\circ(F).
$$
The group $\GL_{\underline{r}}(F) \times \GL_{\underline{r-2}}(F) \times H^e(F)$ acts on $C^\infty(X^\circ(F),\mathcal{H})$ via 
\begin{align*}
\mathcal{R}(g,g',h) f(g_0,m,v)=f\left(\begin{psmatrix} g' & \\ & p_1(h)\end{psmatrix}^{-1}g_0g,g'^tmp_1(h)^{\iota},v_0.p_2(h)\right),
\end{align*}
and this action preserves the subspaces of compactly supported sections and Schwartz sections. The unitarily normalized action is given by 
\begin{align} \label{Ru} \index{$\mathcal{R}_u$}
\mathcal{R}_u(g,g',h):=\frac{|\det g'|}{|\det p_1(h)|^{\underline{(r-2)/2}}}|\chi(g',h)|^{1/2}\mathcal{R}(g,g',h).
\end{align}

We now prepare to define the Harish-Chandra Schwartz space of $Y^{\circ}(F)$ and $(X^{\circ}(F),\mathcal{H}).$
First we require some explicit refinements of well-known bounds for $\Xi$-functions.  In the non-Archimedean case, we always define $\Xi_{\GL_r}$ using the maximal compact subgroup $K_r:=\GL_{r}(\OO_F),$ and if $H \leq \GL_r$ is a unimodular algebraic subgroup, we normalize the Haar measure $dh$ on $H(F)$ so that $dh(H(F) \cap K_r)=1.$
\begin{lem} \label{lem:leq1}
One has $\Xi_{\GL_r}(g) \leq 1.$
\end{lem}
\begin{proof}
    Let $I(0)$ be the normalized induction of the trivial representation of $B_{r}(F)$ to $\GL_{r}(F)$ and let $\varphi_0$ be the spherical vector normalized to have $L^2$-norm $1.$  Then by the Cauchy-Schwarz inequality $
|\Xi_{\GL_{r}}(g)| =|\langle I(0)(g)\varphi_0,\varphi_0\rangle|\leq \norm{I(0)(g)\varphi_0}\norm{\varphi_0}=1.
$
\end{proof}

We recalled log norms in \S \ref{ssec:quasi}. 

\begin{lem} \label{lem:unif}Assume $F$ is non-Archimedean.
For any $r,$ there is a $d_0 \geq 0$ and $C>0$ independent of $F$ such that for all $d \geq d_0,$
$$
 C\delta_{B_r}^{1/2}(t)\le \Xi_{\GL_r}(t)\le \delta_{B_r}^{1/2}(t)\sigma(t)^d.
$$
\end{lem}

\begin{proof}
It follows from \cite[Proposition 4.6.1]{Macdonald} that there is a polynomial $Q\in \CC[X^\ast(T_r)]$ such that 
$\Xi_{\GL_r}(\varpi^{\lambda})\le \delta_{B_r}^{1/2}(\varpi^{\lambda})Q(\lambda)$ for $\lambda$ in the positive Weyl chamber with respect to $B_r.$ The degree of $Q$ is independent of $F$ and the coefficients of $Q$ can be bounded by constants independent of $F.$
Therefore, there exists $d_0$ large independent of $F$ such that $\Xi_{\GL_r}(t)\le \delta_{B_r}^{1/2}(t)\sigma(t)^{d_0}.$ The lower bound can be proved similarly using the formula for $Q$. More precisely, the terms $|P_{W_0}(q^{-1})|$ and $\pi(p)$ in loc. cit are uniformly bounded above and below, and each factor in $\Pi\circ F(x+p)$ is uniformly bounded below.
\end{proof}

Let 
\begin{align}
P_{r-1,1}(R):=\{\begin{psmatrix} g & y \\ & a \end{psmatrix}: (g,y,a) \in \GL_{r-1}(R) \times R^{r-1} \times R^\times\}.
\end{align}
It is a maximal parabolic subgroup of $\GL_r.$  Let $N_{P_{r-1,1}}$ be its unipotent radical and let $M$ be its Levi subgroup.

\begin{lem}\label{lem:unif2} Fix a number field $F.$
There is a $d' \geq 0$ and $C>0$ such that for $d \geq d'$ and any place $v$ of $F$ 
    $$
   \delta_{P_{r-1,1}}^{1/2}(m)\int_{N_{P_{r-1,1}}(F_v)} \Xi_{\GL_r}(mu)\sigma(mu)^{-d}du\le C\Xi_M(m).
    $$
    for $m \in M(F).$
\end{lem}
\begin{proof} It suffices to treat non-Archimedean $v$. This is a refinement of \cite[Proposition II.4.5]{Waldspurger:plancherel}.
We claim that the constants $C_1,C_2,C_3,D$ in the proof of \cite[Proposition II.4.5]{Waldspurger:plancherel}
can be uniformly bounded for all finite $v.$  
These constants ultimately depend on the constants in \cite[Lemmes II.1.1, II.3.1, II.3.2, II.4.4]{Waldspurger:plancherel}.  The constants in \cite[Lemme II.1.1]{Waldspurger:plancherel} can be made uniform by Lemma \ref{lem:unif}.  The proof of \cite[Lemme II.3.1]{Waldspurger:plancherel} implies immediately that its constants may be taken to be uniform.  The constant in \cite[Lemme II.3.2]{Waldspurger:plancherel} is $1$ in our case because $\GL_r$ is split.
The constants in \cite[Lemme II.4.4]{Waldspurger:plancherel} rely only on constants we have already observed are uniform.  

Given these comments, by the proof of \cite[Proposition II.4.5]{Waldspurger:plancherel}, it suffices to prove that the integral
\begin{align}
    \int_{F_{v}^{r-1}} \delta_{B_r}^{1/2}\left(m_0(\begin{psmatrix} I_{r-1} & \\ u &1\end{psmatrix}\right)\sigma(u)^{-d}du
\end{align}
converges for $d$ greater than a constant independent of $v$ and is bounded by a constant independent of $v.$    Here $m_0$ is defined as in \cite[p. 240]{Waldspurger:plancherel} (it is denoted $m_r$ in  \cite[(6.1.1)]{GGH}). 
The integral above is
\begin{align*}
    &1+\sum_{k=1}^\infty \int_{y \in F_v^{r-1}: \norm{y}=q_v^k}\delta_{B_{r}}^{1/2}\left(m_0\begin{psmatrix} I_{r-1} & \\  y &  1 \end{psmatrix}\right) (1+k\log q_v)^{-d}dy\\
    &=1+\sum_{k=1}^\infty \left(q_v^{k(r-1)}-q_v^{(k-1)(r-1)}\right)\delta_{B_{r}}^{1/2}\left(m_0\begin{psmatrix} 1& & \\ &  I_{r-2} &   \\  \varpi^{-k}& &1 \end{psmatrix}\right) (1+k\log q_v)^{-d}\\
    &=1+\sum_{k=1}^\infty \left(q_v^{k(r-1)}-q_v^{(k-1)(r-1)}\right)\delta_{B_{r}}^{1/2}\left(m_0\begin{psmatrix} -\varpi^k& & 1 \\ &  I_{r-2} &   \\  & &\varpi^{-k} \end{psmatrix}\right) (1+k\log q_v)^{-d}\\
    &=1+\sum_{k=1}^\infty (1-q_v^{-1})
     (1+k\log q_v)^{-d}.
\end{align*}
The desired uniformity is clear.
\end{proof}

We define
\begin{align} \label{XiY}
\begin{split}
    &\Xi_{Y^\circ}(m_0p_1(h),v_0.p_2(h))\\&=\min\left\{\frac{|\det p_1(a)|^{\underline{(r-2)/4}}\delta_{B^e}(a)^{1/4}}{\max(1, \delta_{B_e}(a)^{-1/2}|y|)^{1/2}}: h \in N_0(F)\begin{psmatrix} 1 & \Delta(y)\\ & 1 \end{psmatrix}aK_{H^e} \textrm{ for some }(y,a) \in F \times T^e(F) \right\}.\end{split}
\end{align}

\begin{lem} \label{lem:Y:L2}
There is a $d>0$ such that $\Xi_{Y^\circ}\sigma_{Y^{\circ}}^{-d} \in L^2(Y^{\circ}(F)).$
\end{lem}
\begin{proof}
By \eqref{iso:inv} the quantity $\norm{\Xi_{Y^\circ}\sigma_{Y^\circ}^{-d}}_2$ is dominated by
\begin{align*}
    \int_{F \times T^e(F)} \sigma\left(\begin{psmatrix} 1 & \Delta(y) \\ & 1 \end{psmatrix}a\right)^{-d}\left(\frac{|\det p_1(a)|^{\underline{(r-2)/4}}\delta_{B^e}(a)^{1/4}}{\max(1, \delta_{B^e}(a)^{-1/2}|y|)^{1/2}}\right)^2 \frac{dyda}{\delta_{B^e}(a)|\det p_1(a)|^{\underline{(r-2)/2}}}.
\end{align*}
Changing variables in $y$ and using \eqref{prod:log} and \eqref{iso:inv}, this is dominated by 
\begin{align*}
    \int_{F \times T^e(F)} \frac{\sigma\left(y\right)^{-dd_0}\sigma(a)^{-dd_0}}{\max(1, |y|)} dyda
\end{align*}
for some absolute constant $d_0>0.$  This integral
 is finite for $d$ large enough.
\end{proof}
For
$$
(z,g_1,y_1,t_{r-2},t_r,y_2,(k,k'),b,k_0) \in \GL_{\underline{1}}(F) \times \GL_{\underline{r-3}}(F) \times F^{\underline{r-3}} \times \GL_{\underline{1}}(F)^2 \times F^{\underline{1}} \times K_{\underline{r},\underline{r-2}}  \times B^e(F) \times K_{H^e}
$$
let 
\begin{align} \label{Xi}
\Xi_{X^\circ}:X^{\circ}(F)&\lto \RR_{>0}\index{$\Xi_{X^\circ}$}
\end{align}
be the function sending 
\begin{align} \label{coord}
x=\left(\begin{psmatrix} k'^{-1} t_{r-2}\begin{psmatrix}g_1 & y_1&\\ &1 \end{psmatrix}& \\ & p_1\left(k_0\right)^{-1}\begin{psmatrix} 1 & y_2t_r^{-1}\\ & t_r^{-1} \end{psmatrix}\end{psmatrix} zk,k'^tm_0p_1(bk_0)^{\iota}, v_0.p_2(bk_0)\right)
\end{align}
to
\begin{align}
    |t_{r-2}|^{\underline{r-2}}|t_r|^{\underline{(2(r-2)+1)/4}}|\det g_1|^{\underline{5/4}}\Xi^{1/2}_{\GL_{\underline{r-3}}}(g_1)\Xi_{\GL_{\underline{r-2}}}^{1/2}\begin{psmatrix}
        g_1 &y_1 \\
          & 1
    \end{psmatrix}\Xi_{\GL_2}^{1/2}\begin{psmatrix} 1 & y_2t_r^{-1} \\ & t_r^{-1}\end{psmatrix}\Xi_{Y^\circ}(y_0b).
    \end{align}

Let
\begin{align} \label{sigprime}
\sigma_{X^{\circ}}'\left(\begin{psmatrix} k'^{-1} t_{r-2}\begin{psmatrix}g_1 & y_1&\\ &1 \end{psmatrix}& \\ & p_1\left(k_0\right)^{-1}\begin{psmatrix} 1 & y_2t_r^{-1}\\ & t_r^{-1} \end{psmatrix}\end{psmatrix}zk,y\right):=\sigma(t_r)\sigma\begin{psmatrix} g_1 & y_1 \\ & 1 \end{psmatrix} \sigma(z)\sigma\begin{psmatrix} 1 & y_2t_r^{-1} \\ & t_r^{-1} \end{psmatrix}\sigma_{Y^\circ}\left(y\right).
\end{align}
In view of \eqref{prod:log} and \eqref{iso:inv}, one has
\begin{align} \label{d0:bound}
    \sigma_{X^\circ}'(x)^{d_0^{-1}} \leq \sigma_{X^\circ}(x)\leq \sigma_{X^\circ}'(x)^{d_0}
\end{align}
for some absolute constant $d_0>0.$  
\begin{lem} \label{lem:L2}
There is a $d>0$ such that 
$
\Xi_{X^\circ}\sigma_{X^\circ}^{-d} \in L^2(X^{\circ}(F),\mathcal{H}).$
\end{lem}

\begin{proof}
In view of Lemma \ref{lem:Y:L2},  for $d$ sufficiently large  $\norm{\Xi_{X^\circ}\sigma_{X^\circ}^{-d}}^2$ is dominated by 
\begin{align*}
&\int \frac{|t_{r-2}|^{\underline{2(r-2)}}|t_{r}|^{\underline{(2(r-2)+1)/2}}|\det g_1|^{5/2}\Xi_{\GL_{\underline{r-3}}}(g_1)\Xi_{\GL_{\underline{r-2}}}\begin{psmatrix}
        g_1 &y_1 \\
          & 1
    \end{psmatrix}\Xi_{\GL_2}\begin{psmatrix} 1 & y_2t_r^{-1} \\ & t_r^{-1}\end{psmatrix}}{\left(\sigma(t_{r-2})\sigma\begin{psmatrix} g_1 & y_1 \\ & 1 \end{psmatrix}\sigma(z)\sigma\begin{psmatrix} 1 & y_2t_r^{-1} \\ & t_r^{-1} \end{psmatrix}\right)^{d/d_0}} \frac{d^\times zd^\times t_{r-2}d^\times t_r dg_1 dy_1dy_2}{|t_{r-2}|^{2(r-2)}|t_r|^{r-1}|\det g_1|^{3}},
\end{align*}
where the integral is over 
$$
(z,g_1,y_1,t_{r-2},t_r,y_2) \in \GL_{\underline{1}}(F) \times \GL_{\underline{r-3}}(F) \times F^{\underline{r-3}} \times \GL_{\underline{1}}(F)^2 \times F^{\underline{1}}.
$$
Executing the integrals over $t_{r-2},z$ and simplifying, this is 
\begin{align*}
&\int_{\GL_{\underline{r-1}}(F) \times F^{\underline{r-1}}} \frac{\Xi_{\GL_{\underline{r-3}}}(g_1)\Xi_{\GL_{\underline{r-2}}}\begin{psmatrix}
        g_1 &y_1 \\
          & 1
    \end{psmatrix}}{\sigma\begin{psmatrix} g_1 & y_1 \\ & 1 \end{psmatrix}^{d/d_0}}\frac{ dg_1 dy_1}{|\det g_1|^{\underline{1/2}}} \int_{(F^\times)^{\underline{1}} \times F^{\underline{1}}}\frac{\Xi_{\GL_{2}}\begin{psmatrix} 1 & y_2t_r^{-1}\\ & t_r^{-1}\end{psmatrix}}{|t_r|^{\underline{1/2}}\sigma\begin{psmatrix} 1 & y_2t_r^{-1} \\ & t_r^{-1} \end{psmatrix}^{d/d_0}}dy_2d^\times t_r.
\end{align*}
This is bounded for $d>0$ large enough by \cite[Proposition 1.5.1 (iv),(v)]{BP:Ast}.
\end{proof}
For each $d \in \RR$ and $f \in C^\infty(X^\circ(F),\mathcal{H}),$ let 
\begin{align*}
\nu_d(f):=\mathrm{sup}_{x \in X^\circ(F)}\frac{|f(x)|\sigma^d(x)}{\Xi_{X^\circ}(x)}.
\end{align*}
Using these seminorms, one defines the Harish-Chandra space $\mathcal{C}(X^\circ(F),\mathcal{H})$ in the usual manner.  Specifically, when $F$ is non-Archimedean, 
\begin{align*}
\mathcal{C}(X^{\circ}(F),\mathcal{H}):=\bigcup_K\{f \in C^\infty(X^{\circ}(F),\mathcal{H})^K: \nu_{d}(f) \leq \infty \textrm{ for all }d>0\}.\index{$\mathcal{C}(X^\circ(F),\mathcal{H})$}
\end{align*}
Here the union is over all compact open subgroups of $G(F).$  The seminorms $\nu_d$ give $\mathcal{C}(X^{\circ}(F),\mathcal{H})^K$ the structure of Frech\'et space and $\mathcal{C}(X^\circ(F),\mathcal{H})$ is naturally an LF-space.  When $F$ is Archimedean, 
\begin{align*}
\mathcal{C}(X^{\circ}(F),\mathcal{H}):=\{f \in C^\infty(X^{\circ}(F),\mathcal{H}): \nu_{d}(D.f)< \infty \textrm{ for all }D \in \mathcal{U}(\mathfrak{g}) \textrm{ and }d>0\}.
\end{align*}
This is naturally a Frech\'et space.

For some purposes (including Lemma \ref{lem:prod:ineq} below), it  is helpful to relate $\Xi_{X^\circ}$ to norms with respect to suitable affine embeddings.   
In the Archimedean case, let $(\cdot,\cdot)$ be a $K_H$-invariant inner product on ~$F^2 \otimes F^2$ and let $\norm{\cdot}:=(\cdot,\cdot)^{[F:\RR]}.$  Here $K_H$ acts via projection to the $i$th factor followed by the tensor product.
In the non-Archimedean case, let $\norm{v}$ be the maximum of the norms of the coefficients of $v$ with respect to the basis $e_{j_1} \otimes e_{j_2}$ for $j_i \in \{1,2\}.$ 
For $v \in O(v_0)(F)$ set
\begin{align} \label{Plt}
\mathrm{Pl}^{\underline{s}}(v):=\prod_{i=1}^3\norm{\mathrm{Pl}_i^{\otimes 2}(v)}^{s_i},
\end{align}
where $\mathrm{Pl}_i^{\otimes 2}$ is defined as in \eqref{Pli2}.

In  the Archimedean case, let $(\cdot,\cdot)$ be a $K_{r_i}$-invariant inner product on $\wedge^2F^{r_i}$
and let $\norm{\cdot}:=(\cdot,\cdot)^{[F:\RR]}.$  In the non-Archimedean case, for each $i$ the standard basis $e_1,\dots,e_{r_i}$ of $F^{r_i}$ induces a basis of $ \wedge^2F^{r_i}$ in the obvious manner, and we let $\norm{\,\cdot\,}$ be the box norm with respect to this basis. 

In the Archimedean case, let $(\cdot,\cdot)$ be a $K_{\underline{r-2}} \times K_H$-invariant inner product on $M_{r_i-2,2}(F)$ and set $\norm{\cdot}:=(\cdot,\cdot)^{[F:\RR]}.$ 
In the non-Archimedean case, let $\norm{m}$ be the maximum of the norms of the entries of $m \in M_{r_i-2,2}(F).$  
For $m \in M_{\underline{r-2},\underline{2}}(F)$ and $\underline{s} \in \RR^{\underline{1}},$ we abbreviate
\begin{align} \label{norm:m}
    \norm{m}^{\underline{s}}:=\prod_{i=1}^3\norm{m_i}^{s_i}.
\end{align}

\begin{lem} \label{lem:Xi:norm}
After rescaling the inner products by an element of $\RR_{>0}$ in the Archimedean case if necessary, if $d>0$ is greater than an absolute constant one has
    \begin{align*}
    \frac{\Xi_{X^\circ}(g,(m,v))}{\sigma_{X^\circ}(g,(m,v))^d}& \leq \norm{m}^{\underline{(2-r)/2}} \mathrm{Pl}^{\underline{(r-3)/4}}(v)\prod_{i=1}^3\frac{|\det g_i|}{\norm{e_{r_i-1}g_i \wedge e_{r_i}g_i}^{r_i/2}}.
    \end{align*}
\end{lem}

\begin{proof}
Using Lemma \ref{lem:Bisom} and \eqref{Pli2:comp}, one checks that the right hand side of the inequality, evaluated at \eqref{coord}, is
\begin{align*}
    |\det p_1(b)|^{\underline{(r-2)/4}}\delta_{B^e}^{1/4}(b) |t_{r-2}|^{\underline{r-2}}|t_r|^{\underline{(r-2)/2}}|\det g_1|^{\underline{1}}.
\end{align*}
Thus it suffices to prove that 
\begin{align} \label{ineq}
   \frac{|\det g_1|^{\underline{1/4}}\Xi^{1/2}_{\GL_{\underline{r-3}}}(g_1)\Xi_{\GL_{\underline{r-2}}}^{1/2}\begin{psmatrix}
        g_1 &y_1 \\
          & 1
    \end{psmatrix}|t_r|^{\underline{1/4}}\Xi_{\GL_{\underline{2}}}^{1/2}\begin{psmatrix} 1 & y_2t_r^{-1}\\ & t_r^{-1}\end{psmatrix}\Xi_{Y^\circ}(y_0b)}{\sigma\begin{psmatrix} g_1 & y_1 \\ & 1 \end{psmatrix}^d\sigma\begin{psmatrix} 1 & y_2t_{r}^{-1} \\ & t_r^{-1} \end{psmatrix}^d\sigma_{Y^\circ}(y_0b)^{d}|\det p_1(b)|^{\underline{(r-2)/4}}\delta_{B^e}^{1/4}(b)} \leq 1
\end{align}
for $d$ sufficiently large.
By the definition \eqref{XiY} of $\Xi_{Y^\circ}(y_0b)$ one has
$$
\frac{\Xi_{Y^\circ}(y_0b)}{|\det p_1(b)|^{\underline{(r-2)/4}}\delta_{B^e}(b)^{1/4}} \leq 1.
$$
On the other hand
for $d>0$ sufficiently large \cite[Lemma 2.3 and 2.4]{GGH} imply
\begin{align} \label{to:make:unif} \begin{split}
\frac{|\det g_1|^{\underline{1/4}}\Xi_{\GL_{\underline{r-2}}}^{1/2}\begin{psmatrix}g_1 & y_1\\ & 1 \end{psmatrix}}{\sigma\begin{psmatrix} g_1 & y_1 \\ & 1\end{psmatrix}^d} &\ll |\det g_1|^{\underline{1/4}}\Xi^{1/2}_{\GL_{\underline{r-2}}}\begin{psmatrix}
         g & \\
         & 1
     \end{psmatrix} \ll \Xi_{\GL_{\underline{r-2}}}^{1/2}(g_1)\sigma(g_1)^d,\\
     \frac{|t_r|^{\underline{1/4}}\Xi_{\GL_{\underline{2}}}\begin{psmatrix} 1 & y_2t_r^{-1}\\ & t_r^{-1}\end{psmatrix}}{\sigma\begin{psmatrix}1 & y_2t_r^{-1} \\ & t_r^{-1}\end{psmatrix}^d} &\ll |t_r|^{\underline{1/4}}\Xi_{\GL_{\underline{2}}}^{1/2}\begin{psmatrix} t_r & \\ & 1\end{psmatrix} \ll \sigma(t_r)^d.\end{split}
\end{align}
Combined with Lemma \ref{lem:leq1},
 these estimates suffice in the Archimedean setting because we are allowed to renormalize the inner products.  

To complete the proof, it suffices to show that the implicit constants in \eqref{to:make:unif} can be taken to be $1$ outside a finite set of places of $F$ and that $d$ can be chosen uniformly.  In view of lemmas \ref{lem:unif} and \ref{lem:unif2}, inspection of the proofs of \cite[Lemma 2.3 and 2.4]{GGH} implies that we can choose $d$ and the implicit constants to be independent of $F.$  This yields \eqref{ineq} with an unspecified constant $C$ in place of the $1$ on the right hand side.  Since the inequality \eqref{ineq} is trivially valid if $\sigma(g_1)=\sigma(t_r)=\sigma(y_1)=\sigma(y_2)=1,$ by increasing $d$ if necessary we can assume $C=1.$
\end{proof}

\subsection{$\Xi$-rapidly decreasing functions}
First assume that $F$ is non-Archimedean.   We say $f \in C^\infty(X^\circ(F),\mathcal{H})$ is \textbf{$\Xi$-rapidly decreasing} if $f$ is fixed by a compact open subgroup $K <G(F),$ and for all $d>0$ there is a $\Phi_{f,d} \in \mathcal{S}_{\mathrm{ES}}(X(F))$ such that for all $x\in X^\circ(F)$
 \begin{align*}
    \frac{|f|\sigma_{X^\circ}^d}{\Xi_{X^\circ}}(x) \leq \Phi_{f,d}(x). 
    \end{align*}

Now assume $F$ is Archimedean.  
We say that $f \in C^\infty(X^\circ(F),\mathcal{H})$ is \textbf{$\Xi$-rapidly decreasing} if for all $d>0$ and $D \in U(\mathfrak{g})$ there are functions $\Phi_{D,f,d} \in \mathcal{S}_{\mathrm{ES}}(X(F))$ such that for all $x\in X^\circ(F)$
    $$
   \frac{|D.f|\sigma_{X^\circ}^d}{\Xi_{X^\circ}}(x) \leq \Phi_{D,f,d}(x).
    $$
\begin{rem}
The notion of $\Xi$-rapidly decreasing defined here is more restrictive than the analogous definition in \cite{GGH}.  
\end{rem}

\subsection{Adelic sections} \label{ssec:adelic:sec}
In this subsection, we revert to the global setting; thus $F$ is a number field.  We define a Hermitian line bundle $\mathcal{H}$ over $X^\circ(\A_F)$ and the action of $G(\A_F)$ on $C^\infty(X^\circ(\A_F),\mathcal{H})$ via the obvious analogues of the local definitions.

Let
\begin{align} \label{Xi:glob}
\Xi_{X^\circ}:=\prod_v\Xi_{X_{F_v}^\circ} \quad\textrm{and}\quad \index{$\Xi_{X^\circ}$}
\sigma_{X^\circ}:=\prod_v\sigma_{X_{F_v}^\circ}.
\end{align}
Moreover, for $(m,w) \in \mathcal{M}^\circ(\A_F) \times O(v_5)(\A_F)$ and $\underline{s} \in \RR^{\underline{1}},$ let
\begin{align}
\norm{m}^{\underline{s}}:=\prod_{v}\norm{m}_v^{\underline{s}} \quad\textrm{ and } \quad\mathrm{Pl}^{\underline{s}}(w):=\prod_{v}\mathrm{Pl}^{\underline{s}}_v(w)
\end{align}
where the local factors are defined as in \eqref{Plt} and \eqref{norm:m}.

Let $v$ be a place of $F.$
Consider the box norm $\norm{x}_v:=\max\{|x_1|_v,\dots,|x_n|_v\}$ on $F^n_v.$
  It follows from the product formula that
$
\prod_v \norm{x}_v \geq 1
$
for $x \in F^{n}.$
Using this observation, Lemma \ref{lem:Xi:norm}, and the product formula, one obtains the following lemma:

\begin{lem} \label{lem:prod:ineq}
If $x=(g,(m,v)) \in  X^\circ(F),$ then $\frac{\Xi_{X^\circ}(x)}{\sigma_{X^\circ}(x)^d} \ll \mathrm{Pl}^{\underline{(r-3)/4}}(v)$ for $d$ sufficiently large.  \qed
\end{lem}

\subsection{Desiderata for the local Schwartz spaces}
\label{ssec:desi:loc:ss} We now 
formulate the conjectural properties of local Schwartz spaces we require.  Our aim is to isolate just enough desiderata to state and prove the Poisson summation conjecture in our setting (Conjecture \ref{PS:conj}  below).  The reader can consult \cite[\S 4.4, \S 8.1]{WWLi:Zeta}, for example, for additional conjectures on Schwartz spaces in the setting of spherical varieties. 

Let $v$ be a place of $F$ which we omit from notation, writing $F:=F_v.$  
We posit the existence of a Schwartz space \index{$\mathcal{S}(X(F) \times Y(F),\mathcal{H})$} 
$$
\mathcal{S}(
X(F),\mathcal{H}) < \mathcal{C}(X^\circ (F),\mathcal{H})
$$
containing $\mathcal{S}(X^\circ(F),\mathcal{H})$ (see Section \ref{ssec:fspaces}) that is stable under the action of $G(F)$  and is
equipped with a Fourier transform \index{$\mathcal{F}^{}$}
$$
\mathcal{F}:\mathcal{S}(X(F),\mathcal{H}) \lto \mathcal{S}(X(F),\mathcal{H})
$$
satisfying the following desiderata:
\begin{enumerate}    
    \myitem[(1)] \label{FT1} The Fourier transform $\mathcal{F}^{}$ extends to an isometry
    $$
    \mathcal{F}:L^2(X^\circ(F),\mathcal{H}) \tilde{\lto} L^2(X^\circ(F),\mathcal{H}). 
$$
    \myitem[(2)] \label{FT2}
For $(g,g',h) \in G(F)$  one has that 
$$
        \mathcal{F}^{} \circ \mathcal{R}_u(g,g',h)=\mathcal{R}_u(g^{-t},g'^{-t},h^{\iota}) \circ  \mathcal{F}.
        $$
\end{enumerate}  
We point out that functions in $\mathcal{S}(X(F),\mathcal{H})$ need only be defined on $X^\circ(F).$

We say we are in the \textbf{unramified setting} if $F$ is non-Archimedean, unramified over the completion of $\QQ< F$ with respect to the norm on $F$, $\psi$ is unramified and $6=2\cdot 3 \in \OO_F^\times.$  In the unramified setting, we assume the existence of nonzero functions $b \in \mathcal{S}(X(F),\mathcal{H})$ with support in $X(\OO_F)$ that satsify the following:
 \begin{enumerate}   
    \myitem[(3)] \label{FT3} In the unramified setting $\mathcal{F}^{}(b^{})=b^{}.$
\end{enumerate}
The desiderata \ref{FT1}, \ref{FT2}, and \ref{FT3} are standard expectations for local Schwartz spaces. We will give another important spectral desideratum \ref{zeta:basic} for $b$ in  \S \ref{sec:add:des} below.

To work with functions in the Schwartz space, the following additional assumptions are useful:
\begin{enumerate}
    \myitem[(4)]  \label{rapidlydecreasing} Every function  in $\mathcal{S}(X(F),\mathcal{H})$ is $\Xi$-rapidly decreasing. Moreover, in the unramified setting, for any $d>0$ and $x\in X^\circ(F)$
    $$
    \frac{|b|\sigma_{X^\circ}^d}{\Xi_{X^\circ}}(x) \leq \one_{X(\OO_F)}(x)
    $$
    if the residual characteristic is sufficiently large in a sense depending on $d.$
\myitem[(5nA)] \label{nA:unif:smooth} If $F$ is non-Archimedean, every element of $\mathcal{S}(X(F),\mathcal{H})$ is fixed by a compact open subgroup of $G(F).$
    \myitem[(5A)]  \label{Frechet} If $F$ is Archimedean, $\mathcal{S}(X(F),\mathcal{H})$ is a Fr\'echet space and the inclusion $\mathcal{S}(X^{\circ}(F),\mathcal{H}) \to \mathcal{S}(X(F),\mathcal{H})$ is continuous.
\end{enumerate}

\subsection{The adelic Schwartz spaces}
\label{ssec:adelic:S}
We now revert to the global setting, so that $F$ denotes a number field.  We assume the existence of local Schwartz spaces $\mathcal{S}(X(F_v),\mathcal{H})$ for all places $v$ that satisfy 
the desiderata of the previous section. For the Archimedean places, we define
\begin{align}
\mathcal{S}(X(F_\infty),\mathcal{H}):=\widehat{\otimes}_{v|\infty}\mathcal{S}(X(F_v),\mathcal{H})
\end{align}
where the hat denotes the completed projective tensor product with respect to the Fr\'echet structures of \ref{Frechet}. Let $S$ be a finite set of places of $F$ including the infinite places such that for $v \not \in S$ we are in the unramified setting.  We  define
\begin{align}
\mathcal{S}(X(\A_F),\mathcal{H}):=\mathcal{S}(X(F_\infty),\mathcal{H}) \otimes \sideset{}{'}\bigotimes_{v \nmid \infty} \mathcal{S}(X(F_v),\mathcal{H})
\end{align}
where the restricted tensor product is taken with respect to $(b_v)_{v\notin S}.$
For any finite set of places $S'$ of $F$, we define $\mathcal{S}(X(F_{S'}),\mathcal{H})$ and $\mathcal{S}(X(\A_F^{S'}),\mathcal{H})$ analogously. 

\begin{lem} \label{lem:abs:conv00}
For $f \in \mathcal{S}(X(\A_F),\mathcal{H}),$ 
\begin{align*}
\sum_{x \in X^\circ(F)}|f(x)|<\infty.
\end{align*}
\end{lem}

\begin{proof}Using \ref{rapidlydecreasing}, we see that for any $d>0,$ there is a $\Phi_f \in \mathcal{S}_{\mathrm{ES}}(X(\A_F))$ such that  
\begin{align*}
    \sum_{x \in X^\circ(F)}|f(x)| \leq \sum_{x \in X^\circ(F)}\Phi_f(x)\frac{\Xi_{X^\circ}(x)}{\sigma_{X^\circ}(x)^d}.
\end{align*}
By Lemma \ref{lem:prod:ineq} for $d>0$ large enough, this is dominated by 
$$
\sum_{x=(g,(m,v)) \in X^\circ(F)}\Phi_f(x)\mathrm{Pl}^{\underline{(r-3)/4}}(v) \ll 1.
$$
\end{proof}
Since \ref{FT3} holds, the local Fourier transforms induce an adelic Fourier transform $\mathcal{F}:\mathcal{S}(X(\A_F),\mathcal{H}) \tilde{\longrightarrow} \mathcal{S}(X(\A_F),\mathcal{H}).$  For the proof of the functional equation in Theorem \ref{thm:func:eqn} below, we require the following assumption:
\begin{enumerate}
 \myitem[(6)]  \label{FTsquared} One has that $\mathcal{F} \circ \mathcal{F}=\mathcal{R}(g)$ for some $g \in G(F).$
\end{enumerate}

In a moment we will make use of the following support assumption on elements of $\mathcal{S}(X(\A_F),\mathcal{H}):$
\begin{enumerate}
    \myitem[(Sup)] \label{sup} There are two places $v_1$ and $v_2$ of $F$ such that $f=f_{v_1}f_{v_2}f^{v_1v_2} \in \mathcal{S}(X(\A_F),\mathcal{H})$ with $f_{v_1} \in \mathcal{S}(X^\circ(F_{v_1}),\mathcal{H})$ and $\mathcal{F}(f_{v_2}) \in \mathcal{S}(X^\circ(F_{v_2}) ,\mathcal{H}).$
\end{enumerate}
We also make use of assumptions \ref{nA:ratio}, \ref{Arch:ratio}, and \ref{zeta:basic} from \S \ref{sec:add:des} below.  These are the hypothesis that Schwartz functions (and in the case of \ref{zeta:basic}, the basic function) behave as expected with respect to the local zeta integrals from \S\ref{sec:local:int}. 

The following is the most optimistic form of the Poisson summation conjecture in the setting relevant to this paper:

\begin{conj} \label{PS:conj} 
One can choose Schwartz spaces $\mathcal{S}(X(F_v),\mathcal{H})$ for all places $v$ of $F$ such that \ref{FT1}, \ref{FT2}, \ref{FT3}, \ref{rapidlydecreasing}, \ref{nA:unif:smooth},  \ref{Frechet}, \ref{FTsquared}, \ref{nA:ratio}, \ref{Arch:ratio}, and \ref{zeta:basic} hold, and, for all $f \in \mathcal{S}(X(\A_F),\mathcal{H})$ satisfying \ref{sup}, the identity
\begin{align*}
    \sum_{x \in X^\circ(F)}f(x)=\sum_{x \in X^\circ(F)}\mathcal{F}(f)(x)
\end{align*}
holds.
\end{conj}

\section{The local integrals}
\label{sec:local:int}

Let $F$ be a local field and $\psi:F \to \CC^\times$ be a nontrivial additive character. In this section, we introduce and study the local zeta integrals attached to $X(F).$   We establish absolute convergence of the local zeta integrals for $\Xi$-rapidly decreasing functions; see Proposition \ref{prop:ram:conv}.

\subsection{Preliminaries}
Let $\pi$ be a generic unitary irreducible admissible representation of $\GL_r(F)$ with central character $\omega.$  Recall that a \textbf{gauge} $\mathcal{G} \in C^\infty(\GL_r(F))$ is a function of the form 
\begin{align} \label{gauge}
\mathcal{G}(nak):=\left|\prod_{j=1}^{r-1} \alpha_j(a) \right|^{-\ell}\phi(\alpha_1(a),\cdots, \alpha_{r-1}(a)) \index{$\mathcal{G}$}
\end{align}
for $(n,a,k)\in U_r(F)\times T_r(F)\times K_r$. Here $\alpha_i$ is the root $a\mapsto a_ia_{i+1}^{-1},$ $\ell \in \RR_{\geq 0},$ and $\phi\in \mathcal{S}(F^{r-1}).$  We call $\ell$ the \textbf{exponent} of the gauge.  
We can also write the expression above as
\begin{align*}
\mathcal{G}\left(n\begin{psmatrix}\prod_{j=1}^ra_j & & &\\ & \ddots & \\ & & a_{r-1}a_r\\& && a_r\end{psmatrix}k\right):=\left|\prod_{j=1}^{r-1}a_j \right|^{-\ell}\phi(a_1,\cdots, a_{r-1}).
\end{align*}

Let $\mathcal{W}(\pi,\psi)$ be the Whittaker model of $\pi$ with respect to $\psi.$
By \cite[\S 2]{JPSSGL3I} and \cite[Proposition 3.1]{JacquetPerfectRS} every $W\in \mathcal{W}(\pi,\psi)$ is dominated by a gauge $\mathcal{G} \in C^\infty(\GL_r(F))$.  We let $\ell(\pi) \in \RR_{\geq 0}$ \index{$\ell(\pi)$}be the infimum of the set of $\ell$ such that an estimate of the form \eqref{gauge} holds for all $W \in \mathcal{W}(\pi,\psi)$, where $\phi$ is allowed to vary as $W$ and $\ell$ vary. We refer to $\ell(\pi)$ as the \textbf{exponent} of $\pi.$   We say that \textbf{the gauge is unramified} if $F$ is non-Archimedean and $\phi=\one_{\OO_F^{r-1}}.$

\begin{lem} \label{lem:gauge:1:factor} 
Let $(\mathcal{G},\mathcal{G}') \in C^\infty(\GL_r(F)) \times C^\infty(\GL_{r-2}(F))$ be gauges with exponents $(\ell,\ell')$ and let 
$$
(\Phi_1,\Phi_2) \in \mathcal{S}(F^{r+2}) \times  \mathcal{S}(F^\times).
$$
Let $\rho_0,\rho',\rho'',s,s',z \in \RR$ and let $\chi_0:T_{r-2}(F)^2 \times (F^\times)^2 \to \RR_{>0}$ be a quasi-character.
Consider 
    \begin{align*}
\int &\mathcal{G}\begin{psmatrix} t't & &\\ &t_{{r-1}} & \\ & & \lambda^{-1}c^{-1}t_{{r}}\end{psmatrix}\mathcal{G}'(t')\Phi_1(tt_{r-1}t_{r},t_{r-1},t_{r},t_{r-2}',c)\Phi_2\left(t_{r-1}t_r \det t\right)
 \\& 
\times |\lambda|^{s/3+\rho_0}|c|^{s'+\rho'}|\det t'|^{z+\rho''}\chi_0(t',t,t_{{r-1}},t_{{r}})dtdt'd^\times t_{{r-1}} d^\times t_{{r}} d^\times c d^\times \lambda  
\end{align*}
where the integral is over $(t',t,t_{r-1},t_{r},c,\lambda) \in T_{r-2}(F)^2 \times (F^\times)^4.$  
Here $t t_{r-1}t_r$ is shorthand for 
$$
(t_1t_{r-1}t_r,\,t_{2}t_{r-1}t_r,\,\dots,\,t_{r-2}t_{r-1}t_r).
$$
The integral converges if 
$s \gg 1,$ $s'- s/3 \gg 1,$ and $z - s(3(r-1))^{-1}\gg 1.$
Here all implied constants are allowed to depend on $\rho_0,\rho',\rho'',\ell,\ell',r,\chi_0.$
If the gauges $\mathcal{G}$ and $\mathcal{G}'$ are unramified and
$
(\Phi_1,\Phi_2)=(\one_{\OO_{F}^{r+2}},\one_{\OO_F^\times}),
$
then the integral is bounded by $\zeta(3/2)^{2r}$ in the range given above.
\end{lem}

\begin{proof}
Changing variables $\lambda \mapsto c^{-1}t_rt_{r-1}^{-1}\lambda$ and executing the integral over $\lambda,$ we see that the integral is bounded by
$\zeta(s/3+\rho_0-\ell)$ times
\begin{align*}
\int &\mathcal{G}''\begin{psmatrix} t't & \\ &t_{{r-1}} \end{psmatrix}\mathcal{G}'(t')\Phi_1(tt_{r-1}t_{r},t_{r-1},t_{r},t_{{r-2}}',c)\Phi_2\left(t_{r-1}t_r \det t\right)
 \\& 
\times |c|^{{s'}+\rho'-s/3-\rho_0}|\det t'|^{z+\rho''}|t_r|^{s/3+\rho_0}|t_{r-1}|^{-s/3-\rho_0}
\chi_0(t',t,t_{r-1},t_r)dtdt'd^\times t_{{r-1}} d^\times t_{{r}} d^\times c,
\end{align*}
where $\mathcal{G}''$ is a gauge on $C^\infty(\GL_{r-1}(F))$ of exponent $\ell$ which we may take to be unramified if $\mathcal{G}$ is unramified.  Here the integral is over $(t',t,t_{r-1},t_{r},c) \in T_{r-2}(F)^2 \times (F^\times)^3.$  Now we execute the integral over $c$ to see that this is bounded by  $\zeta(s'+\rho'-s/3-\rho_0)$ times
\begin{align*}
\int &\mathcal{G}''\begin{psmatrix} t't & \\ &t_{{r-1}} \end{psmatrix}\mathcal{G}'(t')\Phi_1'(tt_{r-1}t_{r},t_{r-1},t_{r},t_{{r-2}}')\Phi_2\left(t_{r-1}t_r \det t\right)
 \\& 
\times |\det t'|^{z+\rho''}|t_r|^{s/3+\rho_0}|t_{r-1}|^{-s/3-\rho_0}
\chi_0(t',t,t_{r-1},t_r)dt dt' d^\times t_{{r-1}} d^\times t_{{r}},
\end{align*}
where $\Phi_1' \in \mathcal{S}(F^{r+1})$ and the integral is over $T_{r-2}(F)^2 \times (F^\times)^2.$  If $\Phi_1=\one_{\OO_F^{r+2}},$ then we may take $\Phi_1'=\one_{\OO_F^{r+1}}.$ 

We change variables $t \mapsto t_{r-1}t$ to see that the above is 
\begin{align*}
\int &\mathcal{G}''\begin{psmatrix} t't & \\ &1\end{psmatrix}\mathcal{G}'(t')\Phi_1'(tt_{r-1}^2t_{r},t_{r-1},t_{r},t_{{r-2}}')\Phi_2\left(t_{r-1}^{r-1}t_r \det t\right)
 \\& 
\times |\det t'|^{z+\rho''}|t_r|^{s/3+\rho_0}|t_{r-1}|^{-s/3-\rho_0}
\chi_0(t',t_{r-1}t,t_{r-1},t_r)dtdt'd^\times t_{{r-1}} d^\times t_{{r}}   .
\end{align*}
Choosing appropriate coordinates on $T_{r-2}(F),$ we may write
\begin{align*}
\begin{psmatrix}
 t't & \\ & 1
\end{psmatrix}=\begin{psmatrix} \prod_{j=1}^{r-2}t'_jt_j & & & \\ & \ddots & & \\
& & t_{r-2}'t_{r-2} & \\ & & & 1\end{psmatrix}.
\end{align*}
Using the definition of a gauge, there are $(\phi,\phi') \in \mathcal{S}(F^{r-2}) \times \mathcal{S}(F^{r-3})$ such that the integral above is
\begin{align*}
\int &\phi(t_1't_1,\dots,t_{r-2}'t_{r-2})\phi'(t_1',\dots,t_{r-3}')\Phi_1'(tt_{r-1}^2t_r,t_{r-1},t_{r},t_{{r-2}}') \prod_{i=1}^{r-2}|t_it_i'|^{-\ell}\prod_{i=1}^{r-3}|t_i|^{-\ell'}
\\
& \times\Phi_2\left(t_{r-1}^{r-1}t_r t_1t_2^2\cdots t_{r-2}^{r-2}\right)
\chi_0'(t',t,t_{r-1},t_r) |t'_1t_2^{\prime 2}\cdots t_{r-2}^{\prime r-2}|^{z}
|t_r|^{s/3}|t_{r-1}|^{-s/3}
dtd^\times t_{{r-1}} d^\times t_{{r}} dt', 
\end{align*}
where $\chi_0': T_{r-2}(F)^2 \times (F^\times)^2 \to \RR_{>0}$ is a quasi-character.
Here  $(\phi,\phi')=(\one_{\OO_{F}^{r-2}},\one_{\OO_F^{r-3}})$ when $\mathcal{G}$ and $\mathcal{G}'$ are unramified.  We now change variables $t \mapsto tt'^{-1}$ to arrive at
\begin{align*}
\int &\phi(t_1,\dots,t_{r-2})\phi'(t_1',\dots,t_{r-3}')\Phi_1'\left(\frac{t_{1}\dots t_{r-2}t_{r-1}^2t_{r}}{t_1'\dots t_{r-2}'},\dots, \frac{t_{r-2}t_{r-1}^2t_r}{t_{r-2}'},t_{r-1},t_{r},t_{{r-2}}'\right)\prod_{i=1}^{r-2}|t_i|^{-\ell}\prod_{i=1}^{r-3}|t_i'|^{-\ell'}\\
& \times \Phi_2\left(\frac{t_{r-1}^{r-1}t_r t_1\cdots t_{r-2}^{r-2}}{t_1'\cdots t^{\prime r-2}_{r-2}}\right)
\chi_0'(t',t/t',t_{r-1},t_r)|t'_1t_2^{\prime 2}\cdots t_{r-2}^{\prime r-2}|^{z-2}
|t_r|^{s/3}|t_{r-1}|^{-s/3}
dtdt'd^\times t_{{r-1}} d^\times t_{{r}}.
\end{align*}
Now $\Phi_2(x) \ll_\sigma |x|^\sigma$ for any $\sigma,$ and the implied constant is $1$ if $\Phi_2=\one_{\OO_F^\times}.$
Thus the integral above is bounded by 
\begin{align*}
\int &\phi_1(t_1,\dots,t_{r-2},t_{r-1},t_r)\phi'_1(t_1',\dots,t_{r-2}') \prod_{i=1}^{r-2}|t_i|^{-\ell}\prod_{i=1}^{r-3}|t_i'|^{-\ell'}\\
& \times\left|\frac{t_{r-1}^{r-1}t_r t_1\cdots t_{r-2}^{r-2}}{t_1'\cdots t^{\prime r-2}_{r-2}}\right|^{\sigma}
\chi_0'(t',t/t',t_{r-1},t_r)|t'_1t_2^{\prime 2}\cdots t_{r-2}^{\prime r-2}|^{z-2}
|t_r|^{s/3}|t_{r-1}|^{-s/3} dtdt'd^\times t_{{r-1}} d^\times t_{{r}}  
\end{align*}
for some $\phi_1 \in \mathcal{S}(F^r)$ and $\phi_1' \in \mathcal{S}(F^{r-2}).$  When $(\phi,\phi',\Phi_1',\Phi_2)=(\one_{\OO_F^{r-2}},\one_{\OO_F^{r-3}},\one_{\OO_F^{r+1}},\one_{\OO_F^\times}),$  we may take $(\phi_1,\phi_1')=(\one_{\OO_{F}^r},\one_{\OO_F^{r-2}}).$  In general, the $\phi_1,\phi_1'$ depend on $\phi,\phi',\Phi_1',\Phi_2$ and $\sigma.$

Assume $\sigma>0.$  For the following discussion, all constants are allowed to depend on $r,\ell,\ell',\chi_0'.$ The integral over $t_r$ converges if $s \gg 1$ and the integral over $t_{r-1}$ converges if $(r-1)\sigma -s/3 \gg 1.$ Assuming $\sigma \gg 1,$ the integrals over $t_{1},\dots,t_{r-2}$ converge.  Finally if $z -\sigma \gg 1,$
the integral over $t_1',\dots,t_{r-2}'$ converges.  

Now if $s \gg 1$ and  $z - s(3(r-1))^{-1}\gg 1,$ then we can choose $\sigma$ so that the bounds just described are valid.  This proves the convergence assertion in the lemma.  The refined bound in the unramified case also follows.
\end{proof}

When $\pi$ is replaced by a generic unitary irreducible admissible representation of $\GL_{\underline{r}}(F),$ we use the obvious analogues of the definitions above.  In this setting a gauge $\mathcal{G} \in C^\infty(\GL_{\underline{r}}(F))$ is a product of gauges in $C^\infty(\GL_{r_i}(F)).$
 
 \subsection{The local zeta integrals}\label{sec:local zeta}
Let $\pi=\otimes_{i=1}^3 \pi_i$ (resp. $\pi'=\otimes_{i=1}^3 \pi_i'$) be a generic unitary irreducible admissible representation of $\GL_{\underline{r}}(F)$ (resp. $\GL_{\underline{r-2}}(F)$). We let $\omega,\omega_i,\omega',\omega_{i}'$ be the central characters of 
$\pi$, $\pi_i$, $\pi'$ and $\pi_i'$, respectively. We assume $\omega=\omega^{\prime-1}.$ Write
\begin{align*}
\ell(\pi):=(\ell(\pi_1),\ell(\pi_2),\ell(\pi_3)) \quad \textrm{and} \quad \ell(\pi'):=(\ell(\pi_1'),\ell(\pi_2'),\ell(\pi_3')).
\end{align*} 
The character $\nu:H \to\GG_m$ descends to $H/\mathcal{E}$ and hence extends to a character $\nu:H^e \to \GG_m.$
We define a quasi-character of $\GL_{\underline{r-2}}(F)\times H^e(F)$ by\index{$\eta_{\underline{s},s}(g',h)$}
\begin{align}\label{eta}
\eta_{\underline{s},s}(g',h):=|\nu(h)|^{A(\underline{s},s)}\prod_{i=1}^3\frac{|\det g'_i|^{C_i(\underline{s},s)}}{|\det p_1(h_i)|^{C_i(\underline{s},s)(r_i-2)/2}}
\end{align}
where  \index{$C_i(\underline{s},s)$} \index{$A(\underline{s},s)$}
\begin{align} \label{CA} \begin{split}
C_i(\underline{s},s):&=\frac{1}{(r_i-2)/2}\left(s_i+\frac{s-\sum_{j=1}^3 s_j}{2}\right)
=\frac{1}{r_i-2}(s+s_i-s_{i+1}-s_{i+2}),\\
A(\underline{s},s):&=s-\sum_{i=1}^3C_i(\underline{s},s)(r_i-2)/2=\frac{1}{2}(-s+s_1+s_2+s_3). \end{split}
\end{align}
In the definition of $C_i(\underline{s},s),$ the indices should be understood modulo $3.$
We point out that if $aI_{\underline{2}} \in Z_{\GL_{\underline{2}}}(F)$ so that  $p_1(aI_{\underline{2}})=aI_{\underline{2}},$ then
\begin{align} \label{eta:invariance}
\eta_{\underline{s},s}(ag',ah) =\eta_{\underline{s},s}(g',h).
\end{align}

Let $\mathcal{W}(\pi,\psi):=\otimes_{i=1}^3\mathcal{W}(\pi_i,\psi)$ and $ \mathcal{W}(\pi',\overline{\psi}):= \otimes_{i=1}^3\mathcal{W}(\pi'_i,\overline{\psi}).$
For $(W,W') \in \mathcal{W}(\pi,\psi) \times \mathcal{W}(\pi',\overline{\psi})$ and $f \in C^\infty(X^\circ(F),\mathcal{H}),$
we define
\index{$Z(W,W',f,\underline{s},s)$}
\begin{align} \label{loc:zeta02}
\begin{split}
&Z(W,W',f,\underline{s},s):=\iiint W(g)W'(g') \mathcal{R}_u(g,g',h)f(x_0)\eta_{\underline{s},s}(g',h)dg dg'dh \end{split}
\end{align}
where the integrals, from left to right, are over $Z(F) N_0(F)\backslash H^e(F),$
 $N_{\underline{r-2}}(F) \backslash \GL_{\underline{r-2}}(F)$, and  $N_{\underline{r-2},\underline{2}}(F) \backslash \GL_{\underline{r}}(F).$ 
 
Whenever this integral converges absolutely we have
\begin{align} \label{COV} \begin{split}
&Z(W,W',f,\underline{s},s) \\:&=\iiint W(g)W'(g') f_0\left(\begin{psmatrix} g' & \\ & p_1(h) \end{psmatrix}^{-1}g,g'^tm_0p_1(h)^{\iota},v_0.p_2(h)\right)\frac{|\nu(h)||\det g'|^{\underline{3/2}}\eta_{\underline{s},s}(g',h)dg dg'dh}{|\det p_1(h)|^{\underline{3(r-2)/4}}|\Delta(\nu(h))|^{\underline{(r-2)/4}}}\\
&=\iiint W\left(\begin{psmatrix}g' & \\ & p_1(h) \end{psmatrix}g\right)W'(g') f_0\left(g,g'^tm_0p_1(h)^{\iota},v_0.p_2(h)\right)\frac{|\nu(h)||\det p_1(h)|^{\underline{(r-2)/4}}\eta_{\underline{s},s}(g',h)dg dg'dh}{|\Delta(\nu(h))|^{\underline{(r-2)/4}}|\det g'|^{\underline{1/2}}},\end{split}
\end{align}
where the integrals are as before.

\begin{prop} \label{prop:ram:conv}
Assume that $f \in C^\infty(X^\circ(F),\mathcal{H})$ is $\Xi$-rapidly decreasing.  Then the  integral defining $Z(W,W',f,\underline{s},s)$ converges absolutely provided
\begin{align} \label{bounds}
\mathrm{Re}(s) \gg 1, \quad 
\mathrm{Re}(s_i)- \frac{\mathrm{Re}(s)}{3} \gg 1, \quad 
\mathrm{Re}(C_i(\underline{s},s)) -\frac{\mathrm{Re}(s)}{3(r_i-1)}\gg 1.
\end{align}
Here all implied constants are allowed to depend on $\ell(\pi),\ell(\pi'),\underline{r}.$
\end{prop}

\begin{proof}
Without loss of generality, we assume that $(\underline{s},s) \in \RR^{\underline{1}} \times \RR$ and that $f$ is nonnegative.  Let
\begin{align*}
f':=\int_{K_{\underline{r}} \times K_{\underline{r-2}} \times K_{\underline{2}}}\mathcal{R}(k,k',k_0)f dk dk'dk_0.
\end{align*}
Replacing $W,W'$ by gauges $\mathcal{G},\mathcal{G}'$ and using the Iwasawa decomposition together with \eqref{COV}, we see that the integral $Z(W,W',f,\underline{s},s)$ is bounded by 
\begin{align*}
\iiint \mathcal{G}\left(\begin{psmatrix} t' & \\ &p_1(b)\end{psmatrix}g\right)\mathcal{G}'(t')f'(
g,t'm_0p_1(b)^{\iota},v_0.p_2(b))\frac{|\nu(b)||\det p_1(b)|^{\underline{(r-2)/4}}\eta_{\underline{s},s}(t',b)}{|\Delta(\nu(b))|^{\underline{(r-2)/4}}\delta_{B_{\underline{r-2}}}(t')|\det t'|^{\underline{1/2}}}dgdt' d_\ell b
\end{align*}
where the integrals, from left to right, are over 
$Z(F) N_0(F) \backslash B_{H}^e(F),$ $T_{\underline{r-2}}(F),$
 $N_{\underline{r-2},\underline{2}}(F) \backslash \GL_{\underline{r}}(F).$  
 
Recall the definition of $\ell(h)$ for $h \in \GL_{\underline{2}}(F)$ from \eqref{ell:def}.
Using Lemma \ref{lem:Bisom} and notation as in \eqref{J}, we rewrite the integral above as $|3|^{-1}$ times
\begin{align*}
\iint &\mathcal{G}\left(\begin{psmatrix} t' & &\\ &\Delta(\lambda) & \\ & & \Delta([c])c^{-1}\end{psmatrix}g\right)\mathcal{G}'(t')f'\left(
g,t'm_0\Delta(\lambda)^{-1},\left(c_1,c_2,c_3,[c]y\right)\right)\\& \times \frac{\eta_{\underline{s},s}\left(t',\begin{psmatrix} \Delta(\lambda) & \\ & 1 \end{psmatrix}\ell\begin{psmatrix} 1 & \\ & \Delta([c])c^{-1} \end{psmatrix}\right) \chi(c,\lambda)}{|\det t'|^{\underline{1/2}}\delta_{B_{\underline{r-2}}}(t')}dgd^\times \lambda d^\times c dydt'
\end{align*}
for some quasi-character $\chi:\GL_{\underline{1}}(F) \times F^\times \to \RR_{>0},$
where the outer integral is over $F^\times \times \GL_{\underline{1}}(F) \times F \times T_{\underline{r-2}}(F)$ and the inner integral is over $N_{\underline{r-2},\underline{2}}(F) \backslash \GL_{\underline{r}}(F).$ 
Applying the Iwasawa decomposition on $\GL_{\underline{r}}(F),$ this is 
\begin{align*}
\int &\mathcal{G}\begin{psmatrix} t't & &\\ &\Delta(\lambda)t_{\underline{r-1}} & \\ & & \Delta([c])c^{-1}t_{\underline{r}}\end{psmatrix}\mathcal{G}'(t')f'\left(
\begin{psmatrix} tn' & \\ & \begin{psmatrix} t_{\underline{r-1}} & t_{\underline{r-1}} x\\ & t_{\underline{r}} \end{psmatrix}\end{psmatrix},t'm_0\Delta(\lambda)^{-1},\left(c_1,c_2,c_3,[c]y\right)\right)\\& \times \frac{\eta_{\underline{s},s}\left(t',\begin{psmatrix} \Delta(\lambda) & \\ & 1 \end{psmatrix}\ell\begin{psmatrix} 1& \\ & \Delta([c])c^{-1} \end{psmatrix}\right)\chi(c,\lambda)}{|\det t'|^{\underline{1/2}}\delta_{B_{\underline{r-2}}}(t')\delta_{P_{\underline{r-2},\underline{2}}}\begin{psmatrix}
    t & &\\
    & t_{\underline{r-1}} &\\
    & & t_{\underline{r}}
\end{psmatrix}}d^\times \lambda d^\times c d^\times t_{\underline{r-1}} d^\times t_{\underline{r}} dy dt dt' dn'dx
\end{align*} 
where the integral is now over $F^\times \times \GL_{\underline{1}}(F)^3 \times F \times T_{\underline{r-2}}(F)^2 \times N_{\underline{r-2}}(F) \times F^{\underline{1}}.$ 
Changing variables $t' \mapsto \Delta(\lambda)t'$ and then $\lambda \mapsto [c]\lambda,$ using \eqref{eta:invariance} and the fact that $\mathcal{G}$ and $\mathcal{G}'$ are invariant under the center, the above is
\begin{align*}
\int &\mathcal{G}\begin{psmatrix} t't & &\\ &t_{\underline{r-1}} & \\ & & \Delta(\lambda^{-1})c^{-1}t_{\underline{r}}\end{psmatrix}\mathcal{G}'(t')f'\left(
\begin{psmatrix} tn' & \\ & \begin{psmatrix} t_{\underline{r-1}} & t_{\underline{r-1}} x\\ & t_{\underline{r}} \end{psmatrix}\end{psmatrix},t'm_0,\left(c_1,c_2,c_3,[c]y\right)\right)\\& \times \frac{\eta_{\underline{s},s}\left(t',\begin{psmatrix} \Delta(\lambda[c]) & \\ & 1 \end{psmatrix}\ell\begin{psmatrix} \Delta(\lambda[c])^{-1} & \\ & \Delta(\lambda)^{-1}c^{-1} \end{psmatrix}\right)|t_{\underline{r-1}}t_{\underline{r}}|^{\underline{r-2}}\chi(c,[c]\lambda) }{|\det \Delta(\lambda[c])t'|^{\underline{1/2}}|\det t|^{\underline{2}}\delta_{B_{\underline{r-2}}}(t')}d^\times \lambda d^\times c d^\times t_{\underline{r-1}}d^\times t_{\underline{r}} dydtdt'dn'dx.
\end{align*} 
Using  \eqref{Pli2:comp}, Lemma \ref{lem:Xi:norm} and the fact that $f$ is $\Xi$-rapidly decreasing, this is bounded by 
\begin{align*}
\int &\mathcal{G}\begin{psmatrix} t't & &\\ &t_{\underline{r-1}} & \\ & & \Delta(\lambda^{-1})c^{-1}t_{\underline{r}}\end{psmatrix}\mathcal{G}'(t')\Phi\left(
\begin{psmatrix} tn' & \\ & \begin{psmatrix} t_{\underline{r-1}} & t_{\underline{r-1}} x\\ & t_{\underline{r}} \end{psmatrix}\end{psmatrix},t'm_0,\left(c_1,c_2,c_3,[c]y\right)\right)\\& \times \frac{\eta_{\underline{s},s}\left(t',\begin{psmatrix} \Delta(\lambda[c]) & \\ & 1 \end{psmatrix}\ell\begin{psmatrix} \Delta(\lambda[c])^{-1} & \\ & \Delta(\lambda)^{-1}c^{-1} \end{psmatrix}\right)|t_{\underline{r-1}}t_{\underline{r}}|^{\underline{(r-2)/2}}\chi'(c,\lambda)}{|\det t'|^{\underline{1/2}}|\det t|^{\underline{1}}\delta_{B_{\underline{r-2}}}(t')|t_{\underline{r-2}}'|^{\underline{(r-2)/2}}}
d^\times \lambda d^\times c d^\times t_{\underline{r-1}} d^\times t_{\underline{r}}
 dydtdt'dn'dx
\end{align*}
for some $\Phi \in \mathcal{S}_{\mathrm{ES}}(X(F))$ and some quasi-character $\chi':\GL_{\underline{1}}(F) \times F^\times \to \RR_{>0}.$  
We can use the embedding from Lemma \ref{lem:Pl:Iso} to relate $\Phi$ to a Schwartz function on an appropriate affine space.  Using this observation and 
executing the integrals over $n$, $x$ and $y$ the above is 
\begin{align} \label{before:small:comp:cons}\begin{split}
\int &\mathcal{G}\begin{psmatrix} t't & &\\ &t_{\underline{r-1}} & \\ & & \Delta(\lambda^{-1})c^{-1}t_{\underline{r}}\end{psmatrix}\mathcal{G}'(t')\Phi_1(tt_{\underline{r-1}}t_{\underline{r}},t_{\underline{r-1}},t_{\underline{r}},t_{\underline{r-2}}',c)
\Phi_2((\det t)t_{\underline{r-1}}t_{\underline{r}})
\\& \times
\eta_{\underline{s},s}\left(t',\begin{psmatrix} \Delta(\lambda[c]) & \\ & 1 \end{psmatrix}\ell\begin{psmatrix} \Delta(\lambda[c])^{-1} & \\ & \Delta(\lambda)^{-1}c^{-1} \end{psmatrix}\right)\frac{\chi_0(t,t',t_{\underline{r-1}},t_{\underline{r}})\chi'(c,\lambda)}{|c|^{\underline{1}}}dt  dt'd^\times t_{\underline{r-1}}d^\times t_{\underline{r}} d^\times cd^\times \lambda  .\end{split}
\end{align}
Here $\chi_0:T_{\underline{r-2}}(F)^2 \times \GL_{\underline{1}}(F)^2 \to \RR_{>0}$ is a quasi-character,
$\Phi_1 \in \mathcal{S}(F^{\underline{r+2}}),$ $\Phi_2\in \mathcal{S}(\GL_{\underline{1}}(F)),$  and the integral is over $(\lambda,c,t_{\underline{r-1}},t_{\underline{r}},t,t') \in F^\times \times \GL_{\underline{1}}(F)^3 \times T_{\underline{r-2}}(F)^2.$

By \eqref{small:comp} below we have
\begin{align} \label{small:comp:cons} \begin{split} 
\eta_{\underline{s},s}\left(t',\begin{psmatrix} \Delta(\lambda[c]) & \\ & 1 \end{psmatrix}\ell\begin{psmatrix} \Delta(\lambda[c])^{-1} & \\ & \Delta(\lambda)^{-1}c^{-1} \end{psmatrix}\right)
&=|\lambda|^{s}|c|^{\underline{s}}\prod_{i=1}^3|\det t_i’|^{C_i(\underline{s},s)}. \end{split}
\end{align}
We now use \eqref{small:comp:cons} to write \eqref{before:small:comp:cons} as 
\begin{align*}
\int &\mathcal{G}\begin{psmatrix} t't & &\\ &t_{\underline{r-1}} & \\ & & \Delta(\lambda^{-1})c^{-1}t_{\underline{r}}\end{psmatrix}\mathcal{G}'(t')\Phi_1(tt_{\underline{r-1}}t_{\underline{r}},t_{\underline{r-1}},t_{\underline{r}},t_{\underline{r-2}}',c)
\Phi_2((\det t)t_{\underline{r-1}}t_{\underline{r}})
 \\& 
\times |\lambda|^{s+\rho_0}|c|^{\underline{s+\rho'}}\left(\prod_{i=1}^3|\det t_i'|^{C_i(\underline{s},s)}\right)
\chi_0(t,t',t_{\underline{r-1}},t_{\underline{r}})dtdt'd^\times t_{\underline{r-1}} d^\times t_{\underline{r}}   d^\times cd^\times \lambda
\end{align*}
for suitable $(\rho_0,\underline{\rho}') \in \RR \times \RR^{\underline{1}}.$  
Using the fact that gauges are smooth, together with Lemma \ref{lem:Sobolev} if $F$ is Archimedean and trivial considerations if $F$ is non-Archimedean, we can dominate the integral over $ \lambda \in F^\times$ by an integral over $(F^\times)^3.$  This allows us to dominate the expression above by the product over $1 \leq i \leq 3$ of the expressions in Lemma \ref{lem:gauge:1:factor} and deduce the lemma.
\end{proof}

Using the unramified bound in Lemma \ref{lem:gauge:1:factor} and the proof of Proposition \ref{prop:ram:conv}, one deduces the following refinement in the unramified case:

\begin{lem} \label{lem:unr:conv} Suppose $F$ is non-Archimedean. Assume that $\psi, \pi,\pi'$ are unramified, and that $W,W'$ are unramified Whittaker functions, normalized so that $W(I_{\underline{r}})=W'(I_{\underline{r-2}})=1.$   Assume that $f \in C^\infty(X^\circ(F)),\mathcal{H})$ 
is fixed by $\GL_{\underline{r}}(\OO_F) \times \GL_{\underline{r-2}}(\OO_F) \times H^e(\OO_F)$ 
and 
$
(|f|\sigma_{X^\circ}^d\Xi_{X^\circ}^{-1})(x) \leq \one_{X(\OO_F)}(x).
$
Provided the bounds \eqref{bounds} hold, if the residual characteristic is sufficiently large in a sense depending on $d,$ then
one has
$
|Z(W,W',f,\underline{s},s)| \leq \zeta(3/2)^{\sum_{i=1}^32r_i}.$ \qed
\end{lem}

When studying local zeta functions, it is important to know that one can choose the data to be non-vanishing.  Since we expect $\mathcal{S}(X^\circ(F),\mathcal{H})$ to be contained in the Schwartz space, we focus on zeta functions attached to functions in this space.

For $\Phi \in C^\infty(\GL_{\underline{r}}(F) \times Y^{\circ}(F)),$ define
\begin{align} \label{IPhi}
    I_{\Phi}(g,m,v):=\int_{M_{\underline{r-2},\underline{2}}(F)}\overline{\psi}\left(\langle m,z\rangle \right)\Phi\left( \begin{psmatrix} I_{\underline{r-2}} & z\\ & I_2 \end{psmatrix} g, (m,v)\right)dz\,
\end{align}
whenever this integral is absolutely convergent.  If $\Phi \in C_c^\infty(\GL_{\underline{r}}(F) \times Y^{\circ}(F)),$ then $I_{\Phi} \in C_c^\infty(X^\circ(F),\mathcal{H}).$

By \eqref{COV}, if $\Phi:=\Phi_1 \otimes f_1 \in C^\infty(\GL_{\underline{r}}(F)) \otimes C^\infty(Y^{\circ}(F))$ one has that 
\begin{align} \label{int:exp}
\begin{split}
Z(W&,W',I_\Phi,\underline{s},s)\\=&\int_{Z_{\GL_{\underline{2}}}(F)N_0(F) \backslash H^e(F)}\int_{U_{\underline{r-2}}(F) \backslash \GL_{\underline{r-2}}(F)}\int_{ N_{\underline{r-2},\underline{2}}(F) \backslash \GL_{\underline{r}}(F)}W\left(\begin{psmatrix} g' & \\ & p_1(h) \end{psmatrix}g\right)W'(g')\\& 
 \times\int_{M_{\underline{r-2},\underline{2}}(F)}\overline{\psi}(\langle g'^tm_0p_1(h)^{\iota} ,z \rangle)\Phi_1\left(\begin{psmatrix}I_{\underline{r-2}} & z\\ & I_{\underline{2}} \end{psmatrix}g\right)dz\\
 & \times f_1\left(g'^tm_0 p_1(h)^{\iota},v_0.p_2(h)\right) \frac{|\nu(h)||\det p_1(h)|^{\underline{(r-2)/4}}\eta_{\underline{s},s}(g',h)dg dh dg'}{|\Delta(\nu(h))|^{\underline{(r-2)/4}}|\det g'|^{\underline{1/2}}}
   \\
    =&\int_{Z_{\GL_{\underline{2}}}(F)N_0(F) \backslash H^e(F)}\int_{U_{\underline{r-2}}(F) \backslash \GL_{\underline{r-2}}(F)}\pi(\Phi_1)W\begin{psmatrix} g' & \\ & p_1(h) \end{psmatrix}W'(g')\\&\times 
 f_1\left(g'^tm_0 p_1(h)^{\iota},v_0.p_2(h)\right) \frac{|\det p_1(h)|^{\underline{(r-2)/4}}|\nu(h)|}{|\Delta(\nu(h))|^{\underline{(r-2)/4}}|\det g'|^{\underline{1/2}}}
    \eta_{\underline{s},s}(g',h)dg'dh
    \end{split}
    \end{align}
whenever the integrals are absolutely convergent.
We use this computation in Lemma \ref{lem:non:van}, and again in the proof of Theorem \ref{thm:loc:comp} below.

For the next lemma, it is convenient to remark that for any element $\sum_{i=1}^m W_i \otimes W_i' \otimes f_i$ of the algebraic tensor product $\mathcal{W}(\pi,\psi) \otimes \mathcal{W}(\pi',\overline{\psi}) \otimes \mathcal{S}(X^\circ(F),\mathcal{H}),$ the zeta function
\begin{align}
Z\left(\sum_{i=1}^m W_i \otimes W_i' \otimes f_i,\underline{s},s\right):=\sum_{i=1}^mZ(W_i,W_i',f_i,\underline{s},s)
\end{align}
is defined. 

\begin{lem}\label{lem:non:van}
There exists $\sum_{i=1}^m W_i \otimes W'_i \otimes f_i \in \mathcal{W}(\pi,\psi) \otimes \mathcal{W}(\pi',\bar{\psi}) \otimes \mathcal{S}(X^\circ(F),\mathcal{H})$ such that $Z(\sum_{i=1}^m W_i \otimes W_i' \otimes f_i,\underline{s},s)$ is absolutely convergent and  holomorphic as a function of $(\underline{s},s) \in \CC^4.$ We can also assume that $Z(\sum_{i=1}^mW_i \otimes W_i' \otimes f_i,\underline{s},s) \neq 0$ for any particular $(\underline{s},s) \in \CC^4.$
When $F$ is non-Archimedean, we may even choose $W,W',f$ so that $Z(W,W',f,\underline{s},s)=1$ for all $(\underline{s},s).$
\end{lem}

\begin{proof} 
Using the computation
 \eqref{int:exp}, we see that it suffices to show that there is $\sum_{i=1}^m W_i \otimes W_i' \otimes f_i \in \mathcal{W}(\pi,\psi) \otimes  \mathcal{W}(\pi',\overline{\psi}) \otimes  C_c^\infty(Y^{\circ}(F))$ such that the function
\begin{align} \label{first:int:rn} \begin{split}
&\int_{Z_{\GL_{\underline{2}}}(F)N_0(F) \backslash H^e(F)}\int_{U_{\underline{r-2}}(F) \backslash \GL_{\underline{r-2}}(F)}\sum_{i=1}^mW_{i}\begin{psmatrix} g' & \\ & p_1(h)\end{psmatrix}W'_i(g')\\& 
 \times f_i\left(g'^tm_0 p_1(h)^{\iota},v_0.p_2(h)\right) 
\eta_{\underline{s},s}\left(g',h\right)\frac{|\det p_1(h)|^{\underline{(r-2)/4}}|\nu(h)|}{|\Delta(\nu(h))|^{\underline{(r-2)/4}}|\det g'|^{\underline{1/2}}}d\dot{g}'d\dot{h} \end{split}
 \end{align}
of $(\underline{s},s)$ satisfies the assertions of the lemma.

Temporarily denote by $N$ the unipotent radical of the usual mirabolic subgroup $\mathcal{P}_{\underline{r-2}}$ (see \eqref{mira} below).  Let $N^{\mathrm{op}}$ be the opposite of $N,$ thus for $F$-algebras $R$ one has 
\begin{align}
    N^{\mathrm{op}}(R):=\left\{\begin{psmatrix} I_{\underline{r-3}} & \\ x & 1 \end{psmatrix}:x \in R^{\underline{r-3}} \right\}.
\end{align}
The Bruhat decomposition implies that the obvious product map $\mathcal{P}_{\underline{r-2}} \times Z_{\GL_{\underline{r-2}}} \times N^{\mathrm{op}} \to \GL_{\underline{r-2}}$ is an open immersion.  Using this decomposition, we rewrite the integral above as
\begin{align} \label{decomp} \begin{split}
&\int
\sum_{i=1}^m f_i(n^tm_0p_1(h)^{\iota},v_0.p_2(h))\eta_{\underline{s},s}\left(I_{\underline{r-2}},h\right)|\det p_1(h)|^{\underline{(r-2)/4}}
\\& \times \left(\int W_{i}\begin{psmatrix} \Delta(\lambda)pn & \\ & \begin{psmatrix} \Delta(\lambda) & \\ & 1 \end{psmatrix}p_1(h)\end{psmatrix}W'_i(\Delta(\lambda)pn)
 \frac{ \eta_{\underline{s},s}\left(\Delta(\lambda)p,\begin{psmatrix} \Delta(\lambda) & \\ & 1 \end{psmatrix}\right)d_\ell p d^\times \lambda }{|\Delta(\lambda)|^{\underline{(r-2)/2}}|\det p|^{\underline{1/2}}|\lambda|}
\right)
  dndh, \end{split}
 \end{align}
where the outer integral is over  $N^{\mathrm{op}}(F) \times N_0(F) \backslash \SL_{\underline{2}}^e(F),  $ and the inner integral is over 
$U_{\underline{r-2}}(F) \backslash \mathcal{P}_{\underline{r-2}}(F) 
 \times F^\times. $  Here the integral over $Z_{\GL_{\underline{r-2}}}(F)$ occurring in the Bruhat decomposition has been removed due to the $Z_{\GL_{\underline{2}}}(F)$ occurring in the outer integral in \eqref{first:int:rn}.

Let $(W,W') \in \mathcal{W}(\pi,\psi) \times \mathcal{W}(\pi',\overline{\psi})$ be arbitrary.  Using the Dixmier-Malliavin lemma in the Archimedean case or elementary considerations in the non-Archimedean case, we can choose $\widetilde{f}_i \in C_c^\infty(N^{\mathrm{op}}(F) \times \SL_{\underline{2}}^e(F))$ and $(W_i,W_i') \in \mathcal{W}(\pi,\psi) \times \mathcal{W}(\pi',\overline{\psi})$ such that 
\begin{align*}
&\sum_{i=1}^m\int_{N^{\textrm{op}}(F) \times \SL_{\underline{2}}^e(F)} \widetilde{f}_i(n,h)W_i\left(x \begin{psmatrix}n & \\ & p_1(h) \end{psmatrix}\right)W_i'(x'n) \eta_{\underline{s},s}\left(I_{\underline{r-2}},h\right)|\det p_1(h)|^{\underline{(r-2)/4}}dn dh\\&=W(x)W'(x')
\end{align*}
for all $(x,x') \in \GL_{\underline{r}}(F) \times \GL_{\underline{r-2}}(F).$
On the other hand, using Lemma \ref{lem:orb:stab1}, one checks that the map
\begin{align*}
N \times N_0 \backslash \SL_{\underline{2}}^e \lto Y^\circ
\end{align*}
given on points by $(n,h) \mapsto (nm_0p_1(h)^{\iota},v_0.p_2(h))$ is an open immersion.  Hence we can define $f_i \in C_c^\infty(Y^{\circ}(F))$ by stipulating that 
$$
f_i(n^tm_0p_1(h)^{\iota},v_0.p_2(h))=\int_{N_0(F)}\widetilde{f}_i(n,n'h)dn
$$
for $(n,h) \in N^{\mathrm{op}}(F) \times \SL_{\underline{2}}^e(F).$ 
Then \eqref{decomp} is equal to 
\begin{align} \label{red:to:W}
    \int_{U_{\underline{r-2}}(F) \backslash \mathcal{P}_{\underline{r-2}}(F) \times F^\times} W\begin{psmatrix} \Delta(\lambda)p & \\ & \begin{psmatrix} \Delta(\lambda) & \\ & 1 \end{psmatrix}\end{psmatrix}W'(\Delta(\lambda)p)
 \frac{ \eta_{\underline{s},s}\left(\Delta(\lambda)p,\begin{psmatrix} \Delta(\lambda) & \\ & 1 \end{psmatrix}\right)d_\ell p d^\times \lambda }{|\Delta(\lambda)|^{\underline{(r-2)/2}}|\det p|^{\underline{1/2}}|\lambda|}.
\end{align}
We point out that in the non-Archimedean case we may take $m=1.$  

We claim that we can choose $(W,W')\in \mathcal{W}(\pi,\psi) \times \mathcal{W}(\pi',\overline{\psi})$ such that \eqref{red:to:W} is holomorphic as a function of $\underline{s},s$ and nonzero for any particular choice of $\underline{s},s.$   We claim moreover that in the non-Archimedean case we can choose $(W,W')$ so that \eqref{red:to:W} is identically $1.$  Given our previous work the claim implies the lemma.

 Given a smooth function $\phi:\GL_{\underline{r-1}}(F) \to \CC$ compactly supported modulo  $U_{\underline{r-1}}(F)$ on the left such that $\phi(ng)=\psi(n)\phi(g)$ for $n \in U_{\underline{r-1}}(F),$ there is a $W_{\phi} \in \mathcal{W}(\pi,\psi)$ such that $W_{\phi}\begin{psmatrix} g & \\ & 1 \end{psmatrix}=\phi(g)$ \cite[(2.2)]{JPSS:Conv} \cite[Theorem 1]{Kemarsky} \cite[Proposition 5]{Jacquet:quasi-split}.   The claim follows.
\end{proof}

\section{The unramified computation}\label{sec:unramified}

We assume throughout this section that we are in the unramified setting.

\subsection{Single variable computation}\label{subsec:unrsingle}
 Fix $r\ge 3.$ Assume $\pi,\pi'$ are generic unitary unramified representations of $\GL_{r}(F)$ and $\GL_{r-2}(F)$ respectively.
Let $(W,W') \in \mathcal{W}(\pi,\psi)^{\GL_{r}(\OO_F)} \times \mathcal{W}(\pi',\bar{\psi})^{\GL_{r-2}(\OO_F)}$ be the unique pair satisfying $W(I_{r})=W'(I_{r-2})=1.$  

Denote by $\alpha=(\alpha_{1},\dots,\alpha_{r})$ the Langlands class (or Satake parameter)  of $\pi$. Tautologically, there exists a $\kappa> 0$ such that
\begin{align} \label{kappa}
    q^{-\kappa} < |\alpha_i| < q^{\kappa}
\end{align}
for each $i.$
We define $\alpha'$ and $\kappa'$ for $\pi'$ similarly.  
Since $\pi$ and $\pi'$ are generic and unitary, we may take $\kappa=\kappa'=\tfrac{1}{2}$ \cite[Corollary 2.5]{JacquetShalikaEPI}.  If $\pi$ is tempered, then we may take any $\kappa.$

\begin{lem}\label{lem:rankin}
One has
\begin{align*}
    \int_{T_{r-2}(F) }  W\begin{psmatrix} t & &\\ & 1& \\ & & \varpi^{-\ell}\end{psmatrix} W'(t) \frac{|\det t|^{s-1}}{\delta_{B_{r-2}}(t)} dt=\frac{L(s,\pi\times \pi')}{q^{\ell(r-1)/2}}\sum_{n=0}^{r-2} (-1)^{n}\mathbb{S}_{(0,\dots,0, n-\ell)}(\alpha)\mathrm{Tr}(\wedge^n \alpha')q^{-ns}.
\end{align*}
\end{lem}
\noindent Both sides of this identity converge for $\mathrm{Re}(s)$ sufficiently large in a sense depending on $\pi$ and $\pi'.$  The identity holds for all $s$ if understood as an identity of meromorphic functions in $s.$
\begin{proof}

Recall the Schur polynomial for $\GL_r$: Let $v_{r}:=(r-1,r-2,\ldots, 0) \in \ZZ^r$. For $\lambda=(\lambda_i)\in \ZZ_{\ge 0}^r$ with $\lambda_1 \geq \cdots \geq \lambda_r$, define the polynomial
\begin{align*}
    a_{\lambda+v_r}(x):=a_{\lambda+v_r}(x_1,\ldots,x_r):=\det \begin{psmatrix} x_1^{\lambda_1+r-1} & x_2^{\lambda_1+r-1} &\cdots &  x_r^{\lambda_1+r-1}\\
    x_1^{\lambda_2+r-2} & x_2^{\lambda_2+r-2} &\cdots &  x_r^{\lambda_2+r-2}\\
    \vdots & \vdots & \ddots & \vdots \\
      x_1^{\lambda_r} & x_2^{\lambda_r}&\cdots &  x_r^{\lambda_r}\\
    \end{psmatrix} \in \ZZ[x_1,\dots,x_r].
\end{align*} 
The Schur polynomial for $\lambda$ is
\begin{align*} 
    \mathbb{S}_{\lambda_1, \dots, \lambda_r}(x_1,\dots,x_r):=\frac{a_{\lambda+v_r}(x_1,\ldots,x_r)}{a_{v_r}(x_1,\ldots,x_r)}.
\end{align*}
\cite[(A.4)]{Fulton:Harris}.
 The definition directly extends to all $\lambda_1 \geq \cdots \geq \lambda_r$, in which case $\mathbb{S}_{\lambda_1, \dots, \lambda_r}(x_1,\dots,x_r)$ is only a polynomial in $x_1^{\pm 1},\ldots, x_r^{\pm 1}$. By convention $\mathbb{S}_{\lambda_1, \dots, \lambda_r}=0$ unless $\lambda_1 \geq \cdots \geq \lambda_r$. We have
\begin{align} \label{CS}
\delta^{-\frac{1}{2}}_{B_{r}}\begin{psmatrix}\varpi^{\lambda_1} & & \\ & \ddots & \\ & & \varpi^{\lambda_r} \end{psmatrix}W\begin{psmatrix}\varpi^{\lambda_1} & & \\ & \ddots & \\ & & \varpi^{\lambda_{r}} \end{psmatrix}=\mathbb{S}_{\lambda_1,\dots,\lambda_{r}}(\alpha)
\end{align}
by Shintani's formula (a special case of the Casselman-Shalika formula) \cite{ShintaniWhittaker} (see also \cite[\S 3.1.3]{CogdellPCMI}).

For $t=\mathrm{diag}(t_1,\ldots, t_{r-2})\in T_{r-2}(F)$, let $\lambda_j=\mathrm{ ord }\, t_j$. By considering the support of $W$, we deduce that for $t$ in the support of the integral in the lemma, one has that
\begin{align} \label{lambda:cond}
& \lambda_{1} \geq \cdots \geq \lambda_{ r-2}\geq 0.
\end{align}
Thus, by \eqref{CS} 
\begin{align*}
  \int_{T_{r-2}(F) }  W\begin{psmatrix} t & &\\ & 1& \\ & & \varpi^{-\ell}\end{psmatrix} W'(t) \frac{|\det t|^{s-1}}{\delta_{B_{r-2}}(t)} dt= q^{-\ell(r-1)/2}\sum_{\lambda}\mathbb{S}_{(\lambda,0,-\ell)}(\alpha)\mathbb{S}_{\lambda}(\alpha')q^{-s|\lambda|}
\end{align*}
where $|\lambda|=\sum_{j=1}^{r-2} \lambda_j$ and the sum in $\lambda$ is over all $\lambda$ satisfying \eqref{lambda:cond}.

By the cofactor expansion of the determinant $a_{(\lambda,0,-\ell)+v_r}$ along the last row, we have
\begin{align*}
    a_{(\lambda,0,-\ell)+v_r}(x)=\left(\prod_{j=1}^rx_j\right)\sum_{j=1}^r (-1)^{r+j} x_j^{-\ell-1}a_{(\lambda,0)+v_{r-1}}(x_1,\ldots, x_{j-1}, \widehat{x}_{j},x_{j+1},\ldots x_r).
\end{align*}
Here the hat denotes an omitted entry.
Note that $a_{v_r}(x)=\prod_{1\le m<n\le r} (x_m-x_n)$ and hence
\begin{align} \label{avr}
    a_{v_r}(x)/a_{v_{r-1}}(x_1,\ldots, x_{j-1}, \widehat{x}_{j},x_{j+1},\ldots x_r)=(-1)^{j+1}\prod_{\substack{1 \leq n \leq r\\ j \neq n}}(x_j-x_n).
\end{align}
Therefore
\begin{align*}
    \mathbb{S}_{(\lambda,0,-\ell)}(x)=(-1)^{r+1}\left(\prod_{j=1}^rx_j\right)\sum_{j=1}^r \frac{x_j^{-\ell-1}}{\prod_{\substack{1 \leq n \leq r\\ j \neq n}}(x_j-x_n)} \mathbb{S}_{(\lambda,0)}(x_1,\ldots, x_{j-1}, \widehat{x}_{j},x_{j+1},\ldots x_r).
\end{align*}
Thus we have the following identity in $\CC[x_1^{\pm},\dots,x_{r}^{\pm},x_1'^{\pm},\dots,x'^{\pm}_{r-1}][[q^{-s}]]$:
\begin{align*} &q^{-\ell(r-1)/2}\sum_{\lambda}\mathbb{S}_{(\lambda,0,-\ell)}(x)\mathbb{S}_{\lambda}(x')q^{-s|\lambda|}\\
     &=\frac{(-1)^{r+1}}{q^{\ell(r-1)/2}}\left(\prod_{i=1}^rx_i \right)\sum_{j=1}^r \frac{x_j^{-\ell-1}}{\prod_{\substack{1 \leq n \leq r\\ j \neq n}}(x_j-x_n)} \sum_{\lambda} \mathbb{S}_{(\lambda,0)}(x_1,\ldots, x_{j-1}, \widehat{x}_{j},x_{j+1},\ldots x_r)\mathbb{S}_{\lambda}(x')q^{-s|\lambda|}\\
     &=\frac{(-1)^{r+1}}{q^{\ell(r-1)/2}}\left(\prod_{i=1}^rx_i \right)\sum_{j=1}^r \frac{x_j^{-\ell-1}}{\prod_{\substack{1 \leq n \leq r\\ j \neq n}}(x_j-x_n)} \prod_{\substack{1\le n \le r, 1\le m\le r-2\\ n \neq j}} (1-x_n x'_mq^{-s})^{-1}.
\end{align*}
Here we have used the Cauchy identity. Thus
\begin{align} \label{id} \begin{split}
    &q^{-\ell(r-1)/2}\sum_{\lambda}\mathbb{S}_{(\lambda,0,-\ell)}(x)\mathbb{S}_{\lambda}(x')q^{-s|\lambda|}\\
&=q^{-\ell(r-1)/2}\left(\prod_{i=1}^{r}\prod_{k=1}^{r-2}\left(1-x_ix_k'q^{-s} \right)^{-1}\right)(-1)^{r+1}\left(\prod_{i=1}^rx_i \right)
\sum_{j=1}^r \frac{x_j^{-\ell-1} \prod_{1\le m\le r-2} (1-x_j x'_mq^{-s})}{\prod_{\substack{1 \leq n \leq r\\ j \neq n}}(x_j-x_n)}. \end{split}
\end{align}

Now assume $c \in \ZZ$ and $c\ge 2-r.$  By \eqref{avr} we have

\begin{align*}
    \sum_{j=1}^r x_{j}^{-c}\prod_{\substack{1 \leq n \leq r \\n\neq j}} (x_j-x_n)^{-1}
    &=\sum_{j=1}^r(-1)^{j+1}x_j^{-c}\frac{a_{v_{r-1}}(x_1,\ldots, x_{j-1}, \widehat{x}_{j},x_{j+1},\ldots x_r)}{a_{v_r}(x)}\\
    &=\frac{\det \begin{psmatrix} x_1^{-c} & x_2^{-c} &\cdots &  x_r^{-c}\\
    x_1^{r-2} & x_2^{r-2} &\cdots &  x_r^{r-2}\\
    \vdots & \vdots & \ddots & \vdots \\
      1 & 1&\cdots &  1\\
    \end{psmatrix}}{a_{v_r}(x)}=\left(\prod_{i=1}^rx_i\right)^{-c}
    \frac{\det \begin{psmatrix} 1 & 1 &\cdots &  1\\
    x_1^{c+r-2} & x_2^{c+r-2} &\cdots &  x_r^{c+r-2}\\
    \vdots & \vdots & \ddots & \vdots \\
    x_1^{c} & x_2^{c}&\cdots &  x_r^{c}\\
    \end{psmatrix}}{a_{v_r}(x)}\\
    &=(-1)^{r-1}\left(\prod_{i=1}^rx_i\right)^{-c}\mathbb{S}_{(c-1,c-1,\ldots,c-1,0)}(x)
    =(-1)^{r-1}\left(\prod_{i=1}^rx_i\right)^{-1}\mathbb{S}_{(0,0,\ldots,0,1-c)}(x).
\end{align*}
Here in the penultimate inequality we use our assumption that $c \geq 2-r.$  If $0 \geq c \geq 2-r,$ then the equality is just $0=0.$  

In any case, we deduce that
\begin{align} \label{for:Schur} \begin{split}
     &(-1)^{r+1}\left(\prod_{i=1}^r x_i \right)\sum_{j=1}^r \frac{x_j^{-\ell-1}\prod_{1\le m\le r-2}(1-x_jx'_mq^{-s})}{\prod_{\substack{ 1 \leq n \leq r\\n\neq j}} (x_j-x_n)}\\
     &=(-1)^{r+1}\left(\prod_{i=1}^r x_i \right)\sum_{j=1}^r 
     \frac{x_j^{-\ell-1}\sum_{m=0}^{r-2}(-1)^m x_j^m \mathrm{Tr}\left( \wedge^m x' \right)q^{-ms}}{\prod_{\substack{1 \leq n \leq r\\n\neq j}} (x_j-x_n)}\\
     &=\sum_{m=0}^{r-2}(-1)^m\mathbb{S}_{(0, \ldots, 0,m-\ell)}(x)\mathrm{Tr}(\wedge^m x')q^{-ms}. \end{split}
\end{align}
Here we have identified $x'$ with a diagonal matrix in $M_{r-2,r-2}(\ZZ[x'])$. Combining this with \eqref{id} and evaluating at $(x,x')=(\alpha,\alpha'),$ we deduce the lemma.
\end{proof}

\subsection{Unramified local zeta integrals} \label{ssec:unr}
Let $\pi=\pi_1 \otimes \pi_2 \otimes \pi_3$ and $\pi'=\pi_1' \otimes \pi_2' \otimes \pi_3'$ be unramified unitary generic representations of $\GL_{\underline{r}}(F)$ and $\GL_{\underline{r-2}}(F),$ respectively.  Assume central characters $\omega,\omega'$ of $\pi,\pi'$ satisfy $\omega\omega'=1$. Let $W=W_1 \otimes W_2 \otimes W_3$ and $W'=W_1' \otimes W_2' \otimes W_3'$ where $(W_i,W_i') \in \mathcal{W}(\pi_i,\psi)^{\GL_{r_i}(\OO_F)} \times \mathcal{W}(\pi'_i,\psi)^{\GL_{r_i-2}(\OO_F)}$ satisfies $W_i(I_{r_i})=W_i'(I_{r_i-2})=1.$ 

Recall the definition of $\mathrm{Pl}^{\underline{t}}$ from \eqref{Plt}.  Let  
\begin{align} \label{naive} \begin{split}
    &b^{\mathrm{nai}}(g,m,v):=
    \int\overline{\psi}\left(\langle m,z\rangle \right)\one_{\GL_{\underline{r}}(\OO_F)}\left( \begin{psmatrix} I_{\underline{r-2}} & z\\ & I_2 \end{psmatrix} g\right) \one_{Y(\OO_F)}(m,v)dz\mathrm{Pl}^{\underline{(r-2)/4}}(v) \end{split}
\end{align}
where the integral is over $M_{\underline{r-2},\underline{2}}(F).$
This is the analogue of the characteristic function of the integral points in our situation.  

\begin{thm} \label{thm:loc:comp}
For $(\mathrm{Re}(\underline{s}),\mathrm{Re}(s)) \in \RR^{\underline{1}} \times \RR$ satisfying \eqref{bounds}, one has that
\begin{align*}
Z&(W,W',b^{\mathrm{nai}},\underline{s},s)\\=&\sum_{\substack{\underline{k} \in \ZZ_{\geq 0}}}\sum_{\ell=0}^\infty q^{-\ell (s+1/2)}\prod_{i=1}^3 L\left(C_i(\underline{s},s)+\frac{1}{2}
, \pi_i \times \pi_i'\right) 
   \sum_{n_i=0}^{r_i-2}\frac{(-1)^{n_i} \mathbb{S}_{0,\dots,0,n_i-\ell-k_i}(\alpha_i)\mathrm{Tr}(\wedge^{n_i} \alpha'_i)}{q^{n_i(C_i(\underline{s},s)+1/2)+k_i(s_i+1/2)}}.
    \end{align*}
 \end{thm}
\begin{proof}
Using \eqref{int:exp} we see that $Z(W,W', b^{\mathrm{nai}},\underline{s},s)$ is equal to
\begin{align*}
\int_{Z_{\GL_{\underline{2}}}(F)N_0(F) \backslash H^e(F)} &\int_{U_{\underline{r-2}}(F) \backslash \GL_{\underline{r-2}}(F)}
W \begin{psmatrix} g' & \\ &p_1(h) \end{psmatrix} W'(g') \one_{Y(\OO_F)}\left(g^tm_0p_1(h)^{\iota},v_0.p_2(h) \right)\\& \times \mathrm{Pl}^{\underline{(r-2)/4}}(v_0.p_2(h))
\eta_{\underline{s},s}\left(g',h\right)\frac{|\det p_1(h)|^{\underline{(r-2)/4}}|\nu(h)|}{|\Delta(\nu(h))|^{\underline{(r-2)/4}}|\det g'|^{\underline{1/2}}}d\dot{g}'d\dot{h} .
\end{align*} 
By the Iwasawa decomposition of $\GL_{\underline{r-2}}(F) \times H^e(F)$, this equals
\begin{align*}
&\int_{Z_{\GL_{\underline{2}}}(F)N_0(F) \backslash B^e(F)}\int_{F^\times}\int_{T_{\underline{r-2}}(F)}
W\begin{psmatrix}  t & \\ &\begin{psmatrix} \Delta(\lambda)& \\ & 1 \end{psmatrix}p_1(b) \end{psmatrix} W'(t) \one_{Y(\OO_F)}\left(tm_0\Delta(\lambda)^{-1}p_1(b)^{\iota},v_0p_2(b) \right)\\& \times \mathrm{Pl}^{\underline{(r-2)/4}}(v_0p_2(b))
\eta_{\underline{s},s}\left(t,\begin{psmatrix} \Delta(\lambda) & \\ & 1 \end{psmatrix}b\right)\frac{ |\det p_1(b)|^{\underline{(r-2)/4}}dt d^\times \lambda d_\ell b}{\delta_{B_{\underline{r-2}}}(t)|\det t|^{\underline{1/2}}|\lambda|}.
\end{align*}
Here we are viewing $\begin{psmatrix} \Delta(\lambda) & \\ & 1 \end{psmatrix} \in H(F)$ as an element of $H^e(F)$ via the inclusion $H(F) \to H^e(F).$  

Using Lemma \ref{lem:Bisom}, \eqref{Pli2:comp} and notation as in \eqref{ell:def}, we can rewrite the above as
\begin{align*}
\int_{T_{\underline{r-2}}(F) \times F^\times \times F \times (F^\times)^3} &\psi(\lambda y)
W\begin{psmatrix} t &  & \\ & \Delta(\lambda) & \\ & & \Delta([c])c^{-1}  \end{psmatrix} W'(t) \one_{\OO_F^{\underline{1}} \times \OO_F^4}\left(t_{\underline{r-2}}\Delta(\lambda)^{-1},c_1,c_2,c_3,[c] y\right)
\\& \times \eta_{\underline{s},s}\left(t,\begin{psmatrix} \Delta(\lambda) & \\ & 1 \end{psmatrix}\ell\begin{psmatrix} 1 & \\ & \Delta([c])c^{-1} \end{psmatrix}\right)
\frac{ |\Delta([c])c^{-1}|^{\underline{(r-2)/2}}|c|^{\underline{2}}
d^\times \lambda dt  dy 
d^\times c }{\delta_{B_{\underline{r-2}}}(t)|\det t|^{\underline{1/2}}|\lambda|}.
\end{align*}
Here we are writing elements of $V(F)$ as tuples $(a,b,c,d)$ with respect to the basis 
\[
\{e_{156},-e_{246},e_{345},e_{456}\}.
\]

Changing variables $t \mapsto \Delta(\lambda)t,$ using \eqref{eta:invariance}, and recalling that $\omega_i\omega_i'=1,$ the above is
\begin{align*}
\int_{T_{\underline{r-2}}(F) \times F^\times \times F \times (F^\times)^3}&\psi(\lambda y)
W\begin{psmatrix} t &  & \\ & 1 & \\ & &\Delta(\lambda[c]^{-1})^{-1} c^{-1}  \end{psmatrix} W'(t) \one_{\OO_F^5}\left(t_{\underline{r-2}},c_1,c_2,c_3,[c] y\right)\\& \times \eta_{\underline{s},s}\left(t,\begin{psmatrix} \Delta(\lambda) & \\ & 1 \end{psmatrix}\ell\begin{psmatrix} \Delta(\lambda)^{-1} & \\ & \Delta(\lambda^{-1}[c])c^{-1} \end{psmatrix}\right)
\frac{|\Delta([c])c^{-1}|^{\underline{(r-2)/2}}d^\times \lambda dt  dy |c|^{\underline{2}}d^\times c }{\delta_{B_{\underline{r-2}}}(t)|\Delta(\lambda)|^{\underline{(r-2)/2}}|\det t|^{\underline{1/2}} |\lambda|}.
\end{align*}
Executing the integral over $y$ and changing variables $\lambda \mapsto \lambda [c],$ we arrive at
\begin{align} \label{before:small:comp} \begin{split}
&\int_{T_{\underline{r-2}}(F) \times F^\times \times (F^\times)^3}
W\begin{psmatrix} t &  & \\ & 1 & \\ & &\Delta(\lambda)^{-1}c^{-1}  \end{psmatrix} W'(t) \one_{\OO_F^4}\left(c_1,c_2,c_3,\lambda\right)\\& \times \eta_{\underline{s},s}\left(t,\begin{psmatrix} \Delta(\lambda[c]) & \\ & 1 \end{psmatrix}\ell\begin{psmatrix} \Delta(\lambda[c])^{-1} & \\ & \Delta(\lambda^{-1})c^{-1} \end{psmatrix}\right)
\frac{|\Delta([c])c^{-1}|^{\underline{(r-2)/2}} d^\times \lambda dt   d^\times c }{\delta_{B_{\underline{r-2}}}(t)| \Delta(\lambda[c])|^{\underline{(r-2)/2}}|\det t|^{\underline{1/2}}|\lambda|}.  \end{split}
\end{align}
Here we have also used the support properties of Whittaker functions to remove the redundant condition that $t_{\underline{r-2}} \in \OO_F^{\underline{1}}.$

Using the definition \eqref{eta} of $\eta_{\underline{s},s},$ a small computation implies
\begin{align} \label{small:comp} \begin{split}
&\eta_{\underline{s},s}\left(t,\begin{psmatrix} \Delta(\lambda[c]) & \\ & 1 \end{psmatrix}\ell\begin{psmatrix} \Delta(\lambda[c])^{-1} & \\ & \Delta(\lambda)^{-1}c^{-1} \end{psmatrix}\right)
\\&=|\lambda[c]|^{A(\underline{s},s)}
\left(\prod_{i=1}^3\frac{|\det t|^{C_{i}(\underline{s},s)}}{| \lambda c_i|^{-C_i(\underline{s},s)(r_i-2)/2}}\right) \\ 
&=|\lambda|^{A(\underline{s},s)+\sum_{i=1}^3C_i(\underline{s},s)(r_i-2)/2}
 \prod_{i=1}^3|\det t_i|^{C_i(\underline{s},s)}|c_i|^{A(\underline{s},s)+C_i(\underline{s},s)(r_i-2)/2}
\\&=|\lambda|^{s}|c|^{\underline{s}}\prod_{i=1}^3|\det t_i|^{C_i(\underline{s},s)}. \end{split}
\end{align}
Thus \eqref{before:small:comp} is
\begin{align}
    \label{before:Schur} \begin{split}
\int_{T_{\underline{r-2}}(F) \times F^\times \times (F^\times)^3}&
W\begin{psmatrix} t &  & \\ & 1 & \\ & &\Delta(\lambda)^{-1}c^{-1}  \end{psmatrix} W'(t) \one_{\OO_F^4}\left(c_1,c_2,c_3,\lambda\right)\\& \times|\Delta (\lambda)c|^{\underline{(1-r)/2}} |\lambda|^{s+1/2} |c|^{\underline{s+1/2}}\prod_{i=1}^3|\det t_i|^{C_i(\underline{s},s)-\tfrac{1}{2}}
\frac{ d^\times \lambda dt   d^\times c }{\delta_{B_{\underline{r-2}}}(t)}. \end{split}
\end{align}
Letting $(\ell,k_i)=(v(\lambda),v(c_i))$ and using Lemma \ref{lem:rankin}, we see that the integral \eqref{before:Schur} is
\begin{align*}
\sum_{\substack{\underline{k} \in \ZZ_{\geq 0}}}\sum_{\ell=0}^\infty q^{-\ell (s+1/2)}\prod_{i=1}^3 L\left(C_i(\underline{s},s)+\frac{1}{2}
, \pi_i \times \pi_i'\right) 
   \sum_{n_i=0}^{r_i-2}\frac{(-1)^{n_i} \mathbb{S}_{0,\dots,0,n_i-\ell-k_i}(\alpha_i)\mathrm{Tr}(\wedge^{n_i} \alpha'_i)}{q^{n_i(C_i(\underline{s},s)+1/2)+k_i(s_i+1/2)}}.
\end{align*}
\end{proof}

\begin{cor} \label{cor:loc:comp}
Assume $\pi$ and $\pi'$ are tempered.  
Let $\varepsilon>0$ and assume that 
$\mathrm{Re}(s) >\varepsilon,$ $\mathrm{Re}(s_i) > \varepsilon,$ and $C_i(\underline{s},s)>\varepsilon$ (for each $i$).  Then
$$
\frac{Z(W,W',b^{\mathrm{nai}},\underline{s},s)}{
L(s+\tfrac{1}{2},\pi^\vee, \otimes^3)\prod_{i=1}^3L(C_i(\underline{s},s)+\tfrac{1}{2},\pi_i \times \pi_i')L(s_i+\tfrac{1}{2},\pi_i^\vee) }=1+O_\varepsilon(q^{-1-\varepsilon/2}).
$$
\end{cor}

\begin{proof}
Since $\pi$ and $\pi'$ are tempered, one has $|\alpha_i|=|\alpha_i'|=1$ for all $i.$  Thus  
\begin{align*}
&1+q^{-(s+1/2)}\prod_{i=1}^3\mathbb{S}_{0,\dots,0,-1}(\alpha_i)+\sum_{i=1}^3\mathbb{S}_{0,\cdots,0,-1}(\alpha_i)q^{-(s_i+1/2)}+O_{\varepsilon}(q^{-1-\varepsilon})\\&=L(s+\tfrac{1}{2},\pi^\vee,\otimes^3)\prod_{i=1}^3L(s_i+\tfrac{1}{2},\pi_i^\vee)
\end{align*}
if $\mathrm{Re}(s)>\varepsilon$ and $\mathrm{Re}(s_i)>\varepsilon.$

Thus by Theorem \ref{thm:loc:comp} it suffices to observe that
\begin{align*}
\sum_{\substack{\underline{k} \in \ZZ_{\geq 0}^3}}&\sum_{\ell=0}^\infty q^{-\ell (s+1/2)}   \prod_{i=1}^3  \sum_{n_i=0}^{r_i-2}(-1)^{n_i} \mathbb{S}_{0,\dots,0,n_i-\ell-k_i}(\alpha_i)\mathrm{Tr}(\wedge^{n_i}\alpha'_i)q^{-n_i(C_i(\underline{s},s)+1/2)-k_i(s_i+1/2)}\\& =1+q^{-(s+1/2)}\prod_{i=1}^3\mathbb{S}_{0,\dots,0,-1}(\alpha_i)+\sum_{i=1}^3\mathbb{S}_{0,\cdots,0,-1}(\alpha_i)q^{-(s_i+1/2)}+O_{\varepsilon}(q^{-1-\varepsilon})
\end{align*}
under the assumptions of the corollary.
\end{proof}

\begin{lem} \label{lem:Qn}
    Let $I \subset \RR$ be an open interval, and let
    $$
    Q_n:=\{(\alpha_1,\dots,\alpha_n) \in I^n: \alpha_i=\alpha_j \textrm{ if and only if }i=j\}.
    $$
    Then the convex hull of $Q_n$ in $\RR^n$ is $I^n.$
\end{lem}
\begin{proof}
    We proceed by induction on $n.$  The $n=1$ case is obvious.  For $n>1,$ $Q_n$ is equal to
    $$
    \{(\alpha_1,\dots,\alpha_{n-1},\alpha_n):(\alpha_1,\dots,\alpha_{n-1}) \in Q_{n-1} \textrm{ and } \alpha_n \neq \alpha_{i} \textrm{ for }1 \leq i \leq n-1\}.
    $$
    It follows that the convex hull of $Q_n$ contains $Q_{n-1} \times I.$  We now apply the inductive hypothesis to conclude.
\end{proof}

\begin{cor} \label{cor:rational}
   The quotient
    \begin{align}\label{eq:gcd}
    \frac{Z(W,W',b^{\mathrm{nai}},\underline{s},s)}{
    L(s+\tfrac{1}{2},\pi^\vee, \otimes^3) \prod_{i=1}^3L(C_i(\underline{s},s)+\tfrac{1}{2},\pi_i \times \pi_i') L(s_i+\tfrac{1}{2}, \pi_i^\vee)}
    \end{align}
    is a polynomial in $q^{\pm (r_i-2)^{-1}s},q^{\pm (r_i-2)^{-1}s_i}.$  
\end{cor}

\begin{proof}
By Theorem \ref{thm:loc:comp} we have 
\begin{align} \label{holo} \begin{split}
\frac{Z(W,W',b^{\mathrm{nai}},\underline{s},s)}{
\prod_{i=1}^3L(C_i(\underline{s},s)+\tfrac{1}{2},\pi_i \times \pi_i')}=\sum_{\substack{\underline{k} \in \ZZ_{\geq 0}^3}}\sum_{\ell=0}^\infty q^{-\ell (s+1/2)}\prod_{i=1}^3      \sum_{n_i=0}^{r_i-2}\frac{(-1)^{n_i} \mathbb{S}_{0,\dots,0,n_i-\ell-k_i}(\alpha_i)\mathrm{Tr}(\wedge^{n_i}\alpha'_i)}{q^{n_i\left(C_i(\underline{s},s)+1/2\right)+k_i(s_i+1/2)}}. \end{split}
\end{align}
Let $[\omega](\varpi):=\prod_{i=1}^3\omega_i(\varpi)$ and $S_{\underline{r}}:=S_{r_1} \times S_{r_2} \times S_{r_3},$ where $S_{r_i}$ is the symmetric group on $r_i$ letters.
    Using the proof of \Cref{lem:rankin}, specifically \eqref{for:Schur}, the expression above is 
\begin{align*}
    & (-1)^{r_1+r_2+r_3+3}[\omega](\varpi) \sum_{\substack{\underline{k} \in \ZZ_{\geq 0}^3}}\sum_{\ell=0}^\infty q^{-\ell( s+1/2)} \prod_{i=1}^3 \sum_{j=1}^{r_i} \alpha_{ij}^{-\ell-k_i}\frac{\prod_{1\leq m\leq r_i-2} (1-\alpha_{ij}\alpha_{im}'q^{-C_i(\underline{s},s)-1/2})}{\alpha_{ij} \prod_{n\neq j}(\alpha_{ij}-\alpha_{in})q^{k_i(s_i+1/2)}} \\
    = &  (-1)^{r_1+r_2+r_3+3}[\omega](\varpi) \sum_{\substack{\underline{k} \in \ZZ_{\geq 0}^3}}\sum_{\ell=0}^\infty q^{-\ell (s+1/2)} \sum_{\underline{\sigma} \in S_{\underline{r}}} \alpha_{1\sigma_1(1)}^{-\ell-k_1}\alpha_{2\sigma_2(1)}^{-\ell-k_2}\alpha_{3\sigma_3(1)}^{-\ell-k_3} \\ 
    & \times \prod_{i=1}^3 \frac{q^{-k_i(s_i+1/2)}}{(r_i-1)!}  \frac{\prod_{1\leq m\leq r_i-2} (1-\alpha_{i\sigma_i(1)}\alpha_{im}'q^{-C_i(\underline{s},s)-1/2})}{\alpha_{i\sigma_i(1)} \prod_{n\neq \sigma_i(1)}(\alpha_{i\sigma_i(1)}-\alpha_{in})}\\
     = & (-1)^{r_1+r_2+r_3+3}[\omega](\varpi) \sum_{\underline{\sigma} \in S_{\underline{r}}} \left( 1-\frac{q^{-s-\tfrac{1}{2}}}{\alpha_{1\sigma_1(1)}\alpha_{2\sigma_2(1)}\alpha_{3\sigma_3(1)}} \right)^{-1} \prod_{i=1}^3 \frac{1}{(r_i-1)!}\left( 1 - \frac{q^{-s_i-\tfrac{1}{2}}}{\alpha_{i\sigma_i(1)}}\right)^{-1} \\
     & \times \frac{\prod_{1\leq m\leq r_i-2} (1-\alpha_{i\sigma_i(1)}\alpha_{im}'q^{-C_i(\underline{s},s)-1/2})}{\alpha_{i\sigma_i(1)} \prod_{n\neq \sigma_i(1)}(\alpha_{i\sigma_i(1)}-\alpha_{in})} \\
     = & (-1)^{r_1+r_2+r_3+3}[\omega](\varpi) L(s+\tfrac{1}{2}, \pi^{\vee}, \otimes^3)\left(\prod_{i=1}^3L(s_i+\tfrac{1}{2}, \pi_i^{\vee})\right) \\
     & \times \sum_{\underline{\sigma} \in S_{\underline{r}}} \left( \prod_{\substack{ 1\leq a_j \leq r_j \\ (a_1,a_2,a_3) \neq (1,1,1) }} 1- \frac{q^{-s-\tfrac{1}{2}}}{\alpha_{1\sigma_1(a_1)}\alpha_{2\sigma_2(a_2)}\alpha_{3\sigma_3(a_3)}} \right) \prod_{i=1}^3\frac{1}{(r_i-1)!} \left( \prod_{2 \leq b \leq r_i} 1-\frac{q^{-s_i-\tfrac{1}{2}}}{\alpha_{i\sigma_i(b)}} \right) \\
      & \times \frac{\prod_{1\leq m\leq r_i-2} (1-\alpha_{i\sigma_i(1)}\alpha_{im}'q^{-C_i(\underline{s},s)-1/2})}{\alpha_{i\sigma_i(1)} \prod_{n\neq \sigma_i(1)}(\alpha_{i\sigma_i(1)}-\alpha_{in})}.
\end{align*}
Thus it suffices to show that 
\begin{align} \label{to:show:rational} \begin{split}
&(-1)^{r_1+r_2+r_3+3}[\omega](\varpi)\sum_{\underline{\sigma} \in S_{\underline{r}}} \left( \prod_{\substack{ 1\leq a_j \leq r_j \\ (a_1,a_2,a_3) \neq (1,1,1) }} 1- \frac{q^{-s-\tfrac{1}{2}}}{\alpha_{1\sigma_1(a_1)}\alpha_{2\sigma_2(a_2)}\alpha_{3\sigma_3(a_3)}} \right) \\
      & \times \prod_{i=1}^3\frac{1}{(r_i-1)!} \left( \prod_{2 \leq b \leq r_i} 1-\frac{q^{-s_i-\tfrac{1}{2}}}{\alpha_{i\sigma_i(b)}} \right) \frac{\prod_{1\leq m\leq r_i-2} (1-\alpha_{i\sigma_i(1)}\alpha_{im}'q^{-C_i(\underline{s},s)-1/2})}{\alpha_{i\sigma_i(1)} \prod_{n\neq \sigma_i(1)}(\alpha_{i\sigma_i(1)}-\alpha_{in})} \end{split}
\end{align}
     is a polynomial in $q^{\pm (r_i-2)^{-1}s},q^{\pm (r_i-2)^{-1}s_i}.$

For any subset $U \subset \RR^{\underline{r}} \times \RR^{\underline{r-2}} \times \RR^{\underline{1}} \times \RR,$ let 
$$
\mathcal{T}_U:=\{((\underline{z},\underline{z}',\underline{s},s) \in \CC^{\underline{r}} \times \CC^{\underline{r-2}} \times \CC^{\underline{1}} \times \CC: (\mathrm{Re}(\underline{z}),\mathrm{Re}(\underline{z}'),\mathrm{Re}(\underline{s}),\mathrm{Re}(s)) \in U\}
$$ be the tube over $U.$  Let
\begin{align}
    Q \subset \{(\alpha_{ij}),(\alpha_{ij}'),\underline{s},s) \in (q^{-1/2},q^{1/2})^{\underline{2r-2}}  \times \RR^{\underline{1}} \times \RR \}
\end{align}
be the subset in which $\alpha_{ij}\neq \alpha_{in}$ for $j\neq n$ and all $i$. 
Using the equality of \eqref{eq:gcd} and \eqref{to:show:rational}, we can view \eqref{eq:gcd} as a meromorphic function $C(\alpha_{ij},\alpha_{ij}',\underline{s},s)$ on $\CC^{\underline{r}} \times \CC^{\underline{r-2}} \times \CC^{\underline{1}} \times \CC.$  The expression \eqref{to:show:rational} implies that $C$ is holomorphic on $\mathcal{T}_Q.$  On the other hand, \eqref{holo} implies that there is an open cone $Q'\subset \RR^4$ such that $C(\alpha_{ij},\alpha_{ij}',\underline{s},s)$ is holomorphic on $\mathcal{T}_{(q^{-1/2},q^{1/2})^{\underline{2r-2}} \times Q'}$. 

By Lemma \ref{lem:Qn} the convex hull of $Q$ is $(q^{-1/2},q^{1/2})^{\underline{2r-2}}\times \RR^{\underline{1}} \times \RR.$
It follows that the convex hull of $Q\cup \left((q^{-1/2},q^{1/2})^{\underline{2r-2}}\times Q'\right)$ is also $(q^{-1/2},q^{1/2})^{\underline{2r-2}}\times \RR^{\underline{1}} \times \RR.$  On the other hand $Q\cup \left((q^{-1/2},q^{1/2})^{\underline{2r-2}}\times Q'\right)$ is open and path connected, hence connected.  Therefore, by Bochner's tube theorem, $C$ is holomorphic on $\mathcal{T}_{(q^{-1/2},q^{1/2})^{\underline{2r-2}}\times \RR^{\underline{1}} \times \RR}$. Since $C(\alpha_{ij},\alpha_{ij}',\underline{s},s)\prod_{i=1}^3\prod_{j<n} (\alpha_{ij}-\alpha_{in})$ is a polynomial in $q^{\pm (r_i-2)^{-1}s},q^{\pm (r_i-2)^{-1}s_i},$ we deduce that  $C(\alpha_{ij},\alpha_{ij}',\underline{s},s)$ is a polynomial in $q^{\pm (r_i-2)^{-1}s},q^{\pm (r_i-2)^{-1}s_i}$ for any $(\alpha_{ij}),(\alpha_{ij}') \in (q^{-1/2},q^{1/2})^{\underline{2r-2}}.$
\end{proof}

 The following lemma asserts that $b^{\mathrm{nai}}$ satisfies the bound for the basic function hypothesized in \ref{rapidlydecreasing}. 
\begin{lem} \label{lem:L2:NA} 
Let $d>0.$ If the residual characteristic of $F$ is large enough in a sense depending on $d,$ then  one has 
$$
\frac{|b^{\mathrm{nai}}|\sigma_{X^\circ}^d}{\Xi_{X^\circ}}(x) \leq \one_{X(\OO_F)}(x).
$$
\end{lem}

\begin{proof}
We assume $x=(g,m,v)$ is given in coordinates as in \eqref{coord} where $b=\begin{psmatrix} 1 & \Delta(y) \\ & 1 \end{psmatrix} a$ with $a \in T^e(F).$  Then $b^{\mathrm{nai}}(x)=0$ unless $g \in \GL_{\underline{r}}(\OO_F).$  
Assuming $g \in \GL_{\underline{r}}(\OO_F)$ and using \eqref{prod:log}, \eqref{iso:inv}, \eqref{d0:bound} 
one has
\begin{align} \label{bnai:ineq} \begin{split}
&\frac{|b^{\mathrm{nai}}|\sigma_{X^\circ}^d}{\Xi_{X^\circ}}\left(g,k'^tm_0p_1(\begin{psmatrix} 1  & \Delta(y)\\ & 1 \end{psmatrix}ak_0)^{\iota}, v_0.p_2(\begin{psmatrix} 1 & \Delta(y)\\ & 1 \end{psmatrix}ak_0)\right)\\
&\leq  \one_{Y(\OO_F)}(m_0p_1(\begin{psmatrix} 1  & \Delta(y)\\ & 1 \end{psmatrix}a)^{\iota}, v_0.p_2(\begin{psmatrix} 1 & \Delta(y) \\ & 1 \end{psmatrix}a))\\& \times  \mathrm{Pl}^{\underline{(r-2)/4}}(v_0.p_2(\begin{psmatrix} 1 & \Delta(y) \\ & 1 \end{psmatrix}a))\sigma(a)^{dd_0}\sigma(y)^{dd_0}\Xi_{Y^\circ}\left(\begin{psmatrix} 1 & \Delta(y)\\ & 1 \end{psmatrix}a \right)^{-1} \end{split}
\end{align}
for some absolute constant $d_0>0.$

Using coordinates as in Lemma \ref{lem:Bisom} for $a,$ we have 
\begin{align*}
&\one_{Y(\OO_F)}(m_0x^{-1}, (c_1,c_2,c_3,[c]y))\\&
\times \max\left\{\frac{\max(1, \delta_{B^e}^{-1/2}(a')|y'|)^{1/2}}{|\det p_1(a')|^{\underline{(r-2)/4}}\delta_{B^e}(a')^{1/4}}: \begin{psmatrix} 1 & \Delta(y)  \\  & 1 \end{psmatrix} a \in N_0(F)\begin{psmatrix} 1 & \Delta(y')\\ & 1 \end{psmatrix}a'K_{H^e} \right\}\\
&=\frac{\one_{Y(\OO_F)}(m_0x^{-1}, (c_1,c_2,c_3,[c]y))}{|\det p_1(a)|^{\underline{(r-2)/4}}\delta_{B^e}(a)^{1/4}}
=\frac{\one_{Y(\OO_F)}(m_0x^{-1}, (c_1,c_2,c_3,[c]y))}{|x|^{\underline{(r-2)/2}}|\Delta([c])c^{-1}|^{\underline{(r-3)/4}}
}.
\end{align*}
Hence, recalling the definition of $\Xi_{Y^{\circ}}$ from \eqref{XiY}, we see that \eqref{bnai:ineq} is bounded by
\begin{align*}
&\one_{Y(\OO_F)}(m_0x^{-1}, (c_1,c_2,c_3,[c]y)) \mathrm{Pl}^{\underline{(r-2)/4}}(c_1,c_2,c_3,[c]y) \frac{\sigma(x)^{dd_0}\sigma(c)^{dd_0}\sigma(y)^{dd_0}}{\left(|x|^{\underline{(r-2)/2}}|\Delta([c])c^{-1}|^{\underline{(r-3)/4}}
\right)}\\
&=\one_{Y(\OO_F)}(m_0x^{-1}, (c_1,c_2,c_3,[c]y)) \sigma(x)^{dd_0}\sigma(c)^{dd_0}\sigma(y)^{dd_0}\frac{|c|^{\underline{1/2}}}{|x|^{\underline{(r-2)/2}}}.
\end{align*}
This bound suffices to establish the lemma.
\end{proof}

\section{Additional local desiderata} \label{sec:add:des}

Corollary \ref{cor:rational} suggests that the local zeta integals of functions in the Schwartz space should be controlled by certain local $L$-functions.  We make this precise in the current section.  These previously unexplained desiderata were assumed in  Conjecture \ref{PS:conj}. 

The work in this section is not used in the global analysis of the zeta integrals in \S\ref{sec:conv} and \S\ref{sec:unfold}.  It is however necessary in \S\ref{sec:apply} when we apply Conjecture \ref{PS:conj} to deduce the analytic properties of triple product $L$-functions.

Throughout we use Langlands' definition of the local $L$-functions $L(s,\pi^\vee,\otimes^3),$ $L(s,\pi^\vee_i),$ $L(s,\pi_i \times \pi_i')$ and the associated $\gamma$-factors.  Since the local Langlands correspondence for general linear groups is known, this poses no difficulties.  We remind the reader that the definition of $L(s,\pi)$
and $L(s,\pi_i \times \pi_i')$ via Rankin-Selberg theory agrees with the definition due to Langlands; this is essentially built into the local Langlands correspondence.

For any subset $U \subset \RR^{\underline{1}} \times \RR,$ let 
\begin{align}
\mathcal{T}_U:=\{(\underline{s},s) \in \CC^{\underline{1}} \times \CC: (\mathrm{Re}(\underline{s}),\mathrm{Re}(s)) \in U\}
\end{align}
be the tube over $U.$ 
Let $(W,W',f) \in \mathcal{W}(\pi,\psi) \times \mathcal{W}(\pi',\overline{\psi}) \times \mathcal{S}(X(F),\mathcal{H}).$  We expect that $Z(W,W',f,\underline{s},s)$ admits a meromorphic continuation to $\CC^{\underline{1}} \times \CC$ that is rapidly decreasing in tubes away from its poles.  More precisely, when $F$ is Archimedean we expect that:
\begin{enumerate}
\myitem[(7A)] \label{Arch:ratio}  The local zeta function $Z(W,W',f,\underline{s},s)$ admits a meromorphic continuation to $(\underline{s},s) \in \CC^{\underline{1}} \times \CC.$ Let $U \subset \RR^{\underline{1}} \times \RR$ be an open connected set with compact closure $\overline{U}.$ Let $p:\CC^{\underline{1}} \times \CC \to \CC$ be a polynomial such that 
$$
p(\underline{s},s)L(s+\tfrac{1}{2},\pi^\vee,\otimes^3)\prod_{i=1}^3L(C_i(\underline{s},s)+\tfrac{1}{2},\pi_i \times \pi_i')L(s_i+\tfrac{1}{2},\pi_i^\vee)
$$
is holomorphic in an open neighborhood of $\mathcal{T}_{\overline{U}}.$  Then
$$
\mathrm{sup}_{(\underline{s},s) \in \mathcal{T}_{\overline{U}}} |p(\underline{s},s)Z(W,W',f,\underline{s},s)| <\infty.
$$

\end{enumerate}

When $F$ is non-Archimedean, we expect the following: 
\begin{enumerate}
    \myitem[(7nA)] \label{nA:ratio}  The quotient $$\frac{Z(W,W',f,\underline{s},s)}{L(s+\tfrac{1}{2},\pi^\vee,\otimes^3)\prod_{i=1}^3L(C_i(\underline{s},s)+\tfrac{1}{2},\pi_i \times \pi_i')L(s_i+\tfrac{1}{2},\pi_i^\vee)}
    $$
    is a rational function of $q^{\pm (r_i-2)^{-1}s},q^{\pm (r_i-2)^{-1}s_i}.$
\end{enumerate}

Recall from \S \ref{sec:SS} that to define the global Schwartz spaces, we must choose basic sections $b$ in the hypothetical Schwartz space $\mathcal{S}(X(F),\mathcal{H}).$ 
Due to Corollary \ref{cor:rational}, in order for the Poisson summation conjecture to hold for these Schwartz spaces, we require that
\begin{enumerate}
    \myitem[(8)] \label{zeta:basic} In the unramified setting, 
    $$
Z(W,W',b,\underline{s},s)=L(s+\tfrac{1}{2},\pi^\vee,\otimes)\prod_{i=1}^3L(C_i(\underline{s},s)+\tfrac{1}{2},\pi_i \times \pi_i^\vee)L(s_i+\tfrac{1}{2},\pi_i^\vee).
$$
\end{enumerate}
Of course, we also expect a local functional equation, but we omit a discussion to avoid discussing how to normalize ``dual" Whittaker functionals.

For the last result of this section, we revert to global notation.  Thus $F$ is a number field.  For $N_1 \in \ZZ_{>0}$ let 
$$
\mathcal{T}_{N_1}:=\{\underline{s} \in \CC^{\underline{1}}: |\mathrm{Re}(s_i)| <N_1\}
$$
and let $\overline{\mathcal{T}}_{N_1}$ be the closure of $\mathcal{T}_{N_1}$ in $\CC^{\underline{1}}.$

\begin{lem} \label{lem:MI}
Let $\pi=\otimes_{i=1}^3 \pi_i$ and $\pi'$ be unitary irreducible generic representations of $A_{\GL_{\underline{r}}} \backslash \GL_{\underline{r}}(F_\infty)$ and $A_{\GL_{\underline{r-2}}} \backslash \GL_{\underline{r-2}}(F_\infty)$ respectively.  For $v|\infty,$ let $(W_v,W'_v,f_v) \in \mathcal{W}(\pi_v,\psi_v) \times \mathcal{W}(\pi'_v,\overline{\psi}_v) \times \mathcal{S}(X(F_v),\mathcal{H}).$
Assume \ref{rapidlydecreasing} and \ref{Arch:ratio}.  Let $N_1 \in \ZZ_{>0},$ and let $p \in \CC[\underline{s}]$ be chosen so that $p(\underline{s})\prod_{v|\infty}L(s_i+\tfrac{1}{2},\pi_{iv}^\vee)$ has no pole in $\overline{\mathcal{T}}_{N_1}.$
Let
\begin{align}
\phi_{p,\underline{s}}(a):=\frac{1}{2\pi i}\int_{\mathrm{Re}(s)=\sigma} a^{-s}\left(p(\underline{s})\prod_{v|\infty}Z(W_v,W_v',f_v,\underline{s},s)\right) ds
\end{align}
where we assume $\sigma=\mathrm{Re}(s)$ and $\mathrm{Re}(\underline{s})$ satisfy \eqref{bounds}.
Then $\phi_{p,\underline{s}}(a)$ admits an analytic continuation to all $\underline{s} \in \mathcal{T}_{N_1}$ that satisfies
\begin{align*}
    |\phi_{p,\underline{s}}(a)|\ll_{N_1,N_2,p,W_v,W_v',f_v} a^{-N_2}
\end{align*}
for all $N_2>1$ and  $a  \geq 1.$
\end{lem}
\begin{proof}
The zeta integral in the lemma is convergent for $\sigma=\mathrm{Re}(s)$ and $\mathrm{Re}(\underline{s})$ satisfying \eqref{bounds} by \ref{rapidlydecreasing} and Proposition \ref{prop:ram:conv}.

Using \ref{Arch:ratio} we can shift the contour to $\mathrm{Re}(\sigma) \gg 1$ without passing any poles of the integrand.  
Moreover, again by \ref{Arch:ratio} for $\sigma \gg_{N_1} 1$ sufficiently large $p(\underline{s})\prod_{v|\infty}Z(W_v,W_v',f_v,\underline{s},s)$ has no poles with $|\mathrm{Re}(s)-\sigma|<\tfrac{1}{2}$ and $\underline{s} \in \mathcal{T}_{N_1}$ and is rapidly decreasing as a function of $s$ for $|\mathrm{Re}(s)-\sigma|<\tfrac{1}{2}$ provided $\underline{s} \in \mathcal{T}_{N_1}.$  This implies the analytic continuation statement.

To prove the bound, 
we observe that 
\begin{align*}
 a^{N_2}\phi_{p,\underline{s}}(a)&=\frac{1}{2\pi i }\int_{\mathrm{Re}(s)=\sigma}a^{N_2-s}\left(p(\underline{s})\prod_{v|\infty}Z(W_v,W_v',f_v,\underline{s},s)\right) ds\\
&=\frac{1}{2\pi i }\int_{\mathrm{Re}(s)=\sigma-N_2}a^{-s}\left(p(\underline{s})\prod_{v|\infty}Z(W_v,W_v',f_v,\underline{s},s+N_2)\right) ds\\
&=\frac{1}{2\pi i }\int_{\mathrm{Re}(s)=\sigma}a^{-s}\left(p(\underline{s})\prod_{v|\infty}Z(W_v,W_v',f_v,\underline{s},s+N_2)\right) ds
\end{align*}
Here in the last line we have employed a contour shift as in the previous paragraph.
Assume that $\sigma+N_2 \gg_{N_1} 1.$  By \ref{Arch:ratio}, for $\underline{s} \in \mathcal{T}_{N_1},$ the expression in parenthesis is rapidly decreasing as a function of $s$ satisfying $|\mathrm{Re}(s)-\sigma|<\tfrac{1}{2}.$  Thus the integral above is $O_{N_1,N_2,W_v,W_{v'}f_v}(1).$  This implies the estimate in the lemma.
\end{proof}

\section{Convergence in a cone} \label{sec:conv}

For the moment we focus on a single factor $\GL_{r_i}$ of $\GL_{\underline{r}}.$  We omit $i$ from notation for simplicity.  
Let $K_{\infty}\GL_{r}(\widehat{\OO}_F)=\prod_vK_v$ be a maximal compact subgroup of $\GL_r(\A_F).$  
For places $v$ of $F$ and
$(p,k) \in P_{r-2,2}(F_v) \times K_{v},$
let 
\begin{align} \label{Xi1factor} \begin{split}
\Xi_{r-2,2,F_{v}}:X_r^\circ(F_v)&\lto \RR\\
N_{r-2,2}(F_v)pk &\longmapsto \delta^{1/2}_{P_{r-2,2}}(p). \end{split}
\end{align}
Moreover let
\begin{align}
\Xi_{r-2,2}:&=\prod_{v}\Xi_{r-2,2,F_v}.
\end{align}

Let $P_{r-2,1,1} \leq \GL_{r}$ be the standard parabolic subgroup of type $(r-2,1,1)$ and let $N_{r-2,1,1}$ be its unipotent radical. 

\begin{lem} \label{lem:factor:mod:growth}
Let $\Phi_1 \in \mathcal{S}_{\mathrm{ES}}(X_r(\A_F)).$ 
Provided $\mathrm{Re}(s) \gg_r 1$ and $(r-1)\mathrm{Re}(z)-\mathrm{Re}(s) \gg_r 1,$ the integral
\begin{align*}
&\int_{[\GG_m]^3 }\sum_{x \in N_{r-2,1,1}(F) \backslash \GL_r(F)} \Xi_{r-2,2}\left(xg\right)\left(\int_{\A_F}\Phi_1\left(\begin{psmatrix} (a'g'')^{-1} & &\\ & 1 & \\ & &\lambda\end{psmatrix}\begin{psmatrix} I_{r-2} & & \\ & 1 & y \\ & & 1\end{psmatrix}xag\right)dy\right) 
|\lambda|^{s}|\det a'|^{z} d^\times a d^\times a'd^\times \lambda  
\end{align*}
converges absolutely and is of moderate growth as a  function of $(g,g'') \in [\GL_r]\times [\GL_{r-2}].$
\end{lem}
\begin{proof} 
We assume without loss of generality that $s,z \in \RR.$  All implicit constants in this proof are allowed to depend on $r.$

Using Lemma \ref{lem:FL} we see that it suffices to prove that for all $\Phi \in \mathcal{S}_{\mathrm{ES}}(X_r(\A_F))$ the expression 
\begin{align} \label{deg:ES} \begin{split}
E(g,g''):=&\sum_{x \in P_{r-2,1,1}(F) \backslash \GL_r(F)} \Xi_{r-2,2}\left(xg\right)\\& \times \int_{\GL_{r-2}(\A_F) \times (\A_F^\times)^2 \times \A_F}\Phi\left(\begin{psmatrix} (g'g'')^{-1} & &\\ & 1 & \\ & &\lambda\end{psmatrix}\begin{psmatrix} I_{r-2} & & \\ & 1 & y \\ & & 1\end{psmatrix}xag\right)dy 
|\lambda|^{s}|\det g'|^{z} d^\times a dg'd^\times \lambda \end{split}
\end{align}
converges for $(s,z)$ as in the lemma and is of moderate growth as a function of $(g,g'') \in [\GL_r] \times [\GL_{r-2}].$  Provided \eqref{deg:ES} converges absolutely, a change of variables shows $E(g,g'')=|\det g''|^{-z}E(g,I_{r-2}),$ so it suffices to show that $E(g,I_{r-2})$ converges for $(s,z)$ as in the lemma and is of moderate growth as a function of $g \in [\GL_r].$

As the notation indicates, $E(g,I_{r-2})$ is a degenerate Eisenstein series on $\GL_r(\A_F)$.  It is well-known to converge and define a function of moderate growth for $s \gg 1$ and $z \gg 1$ provided that the section
\begin{align*}
g\longmapsto \Xi_{r-2,2}\left(g\right) \int_{\GL_{r-2}(\A_F) \times (\A_F^\times)^2 \times \A_F}\Phi\left(\begin{psmatrix} g'^{-1} & &\\ & 1 & \\ & &\lambda\end{psmatrix}\begin{psmatrix} I_{r-2} & & \\ & 1 & y \\ & & 1\end{psmatrix}ag\right)dy 
|\lambda|^{s}|\det g'|^{z} d^\times a dg'd^\times\lambda 
\end{align*}
converges \cite[\S II.1.5]{MW:Spectral:Decomp:ES}.

 Using the embedding $\det \times \mathrm{Pl}_{X_r}$ of Lemma \ref{lem:Pl:Iso}, we see that this is equivalent to the assertion that for all 
$$
(\Phi_1,\Phi_2,\Phi_3) \in \mathcal{S}(M_{r-2,r-2}(\A_F)) \times \mathcal{S}(M_{2,2}(\A_F)) \times \mathcal{S}(\A_F^\times)
$$
the integral
\begin{align*}
\int_{\GL_{r-2}(\A_F) \times (\A_F^\times)^2 \times \A_F}
\Phi_1(g'^{-1}a^3\lambda)\Phi_2\begin{psmatrix} a & ay\\ & \lambda a\end{psmatrix}\Phi_3(a^r\lambda\det g'^{-1}) dy |\lambda|^s|\det g'|^{z}d^\times adg'd^\times \lambda
\end{align*}
converges for $(s,z)$ as in the lemma.
We change variables $(g',y) \mapsto ((a^3\lambda)g'^{-1},a^{-1}y)$ and then execute the integral over $y$ to see that the above is 
\begin{align*}
\int_{\GL_{r-2}(\A_F) \times (\A_F^\times)^2}
\Phi_1(g')\Phi_2'(a,\lambda a)\Phi_3(a^{r-3(r-2)}\lambda^{3-r} \det g')  |\lambda|^{s+(r-2)z}|a|^{3(r-2)z}|\det g'|^{-z}\frac{d^\times adg'd^\times\lambda}{|a|}
\end{align*}
for some $\Phi_2' \in \mathcal{S}(\A_F^2).$  Changing variables $\lambda \mapsto \lambda a^{-1},$ this is 
\begin{align*}
\int_{\GL_{r-2}(\A_F) \times (\A_F^\times)^2}
\Phi_1(g')\Phi_2'(a,\lambda )\Phi_3(a^{3-r
}\lambda^{3-r} \det g')  |\lambda|^{s+(r-2)z}|a|^{-s+
2(r-2)z-1}\frac{d^\times adg'd^\times\lambda}{|\det g'|^z}.
\end{align*}
For any $\sigma >0,$ this is bounded by a constant depending on $\sigma$ times 
\begin{align*}
\int_{\GL_{r-2}(\A_F) \times (\A_F^\times)^2}
\Phi_1(g')\Phi_2'(a,\lambda ) |\lambda|^{s+(r-2)z+(3-r)\sigma}|a|^{-s+2(r-2)z+(3-r)\sigma-1}|\det g'|^{\sigma-z}d^\times adg'd^\times\lambda.
\end{align*}

Now for $\Phi \in \mathcal{S}(M_{n,n}(\A_F)),$ the integral
\begin{align*}
\int_{\GL_{n}(\A_F)}\Phi(g)|\det g|^{s_0}dg
\end{align*}
converges absolutely for $\mathrm{Re}(s_0)>n$ by \cite[Lemma 12.5]{GodementJacquetBook}.
Thus the integral above converges provided that $\sigma-z>(r-2),$ $s+z(r-2)+(3-r)\sigma>1,$ and $-s+2(r-2)z+(3-r)\sigma>2.$  It is possible to choose $\sigma$ satisfying these conditions for $(s,z)$ as in the lemma.
\end{proof}

We now return to considerations involving all three factors of $\GL_{\underline{r}}.$
Since we do not yet have a definition of $\mathcal{S}(X(\A_F),\mathcal{H}),$ we consider sections  satisfying fairly weak conditions that we expect to be valid for elements of the Schwartz space.  
Recall the definitions of $\Xi_{X^\circ}$ and $\sigma_{X^\circ}$ from  \eqref{Xi:glob}.
We consider sections $f$ satisfying the 
following conditions:
\begin{enumerate}
    \myitem[($f$1)]\label{ad1} 
    $
    f=f_\infty \otimes f^\infty \in C^\infty(X^\circ(F_\infty),\mathcal{H})  \otimes C^\infty(X^\circ(\A_F^\infty),\mathcal{H}).
    $
    \myitem[($f$2)]\label{ad2} $f^\infty$ is finite under a compact open subgroup of $G(\A_F^\infty).$
    \myitem[($f$3)] \label{ad3} For any $d>0,$ there is a function
    $
    \Phi \in \mathcal{S}_{\mathrm{ES}}(X(\A_F))
    $
    such that 
    $$
    \left|\frac{f\sigma_{X^\circ}^d}{\Xi_{X^\circ}}(x)\right| \leq \Phi(x).
    $$
\end{enumerate}
For  $g \in G(\A_F)$ and $f$ satisfying \ref{ad1}, \ref{ad2} and \ref{ad3}, we set 
\begin{align*} 
    \Theta_{f}
    (g)&= \sum_{x \in X^\circ(F)}\mathcal{R}_u(g)f(x).
    \end{align*}
    The sum converges absolutely by the proof of Lemma \ref{lem:abs:conv00}. We define \index{$H^e(\A_F)'$}
\begin{align} \label{He1} \begin{split}
    H^e(\A_F)':&=\mathrm{ker}(|\det \circ p_1|:H^e(\A_F) \lto \RR_{>0}^{\underline{1}}\},\\
[H^e]':&=H^e(F) \backslash H^e(\A_F)'. \end{split}
\end{align} 
Thus $H^e(\A_F)'/H^e(\A_F)^1 \cong \RR_{>0}.$  

\begin{lem} \label{lem:abs:conv0}
Let $\omega:A_{\GL_{\underline{1}}} \backslash [\GL_{\underline{1}}] \to \CC^\times$ be a unitary character.
For $f$ satisfying \ref{ad1}, \ref{ad2}, \ref{ad3} and $(\underline{s},s) \in \CC^{\underline{1}} \times \CC$ satisfying  
\begin{align}\label{bounds2}
\mathrm{Re}(s) \gg 1, \quad 
\mathrm{Re}(s_i)-\mathrm{Re}(s)/3  \gg 1,  \quad 
\mathrm{Re}(C_i(\underline{s},s))-\frac{\mathrm{Re}(s)}{3(r_i-1)} \gg 1,
 \end{align}
the integral
\begin{align*}
\int_{[\GL_{\underline{1}}] \times [\GL_{\underline{1}}] \times [H^e]'}\Theta_f(ag,a'g,h)\omega(a)\overline{\omega}(a')\eta_{\underline{s},s}(a'I_{\underline{r-2}},h)d^\times ad^\times a' dh
\end{align*}
converges absolutely and is of moderate growth as a function of $[\GL_{\underline{r}}] \times [\GL_{\underline{r-2}}].$  Here the implicit constants depend on $\underline{r}.$
\end{lem}

\begin{proof}
We assume without loss of generality that $\omega$ is trivial and $(\underline{s},s) \in \RR^{\underline{1}} \times \RR.$  All implicit constants in this proof are allowed to depend on $\underline{r}.$

Let $Y_{v_0}^\circ$ be the fiber of $Y^\circ \to V$ over $v_0.$
By \eqref{Htildestab} $H^e_{v_0}=Z_{\GL_{\underline{2}}} \wedge^\mathcal{E}\mathcal{E}T_H N_0\cong Z_{\GL_{\underline{2}}} \times T_H \times N_0.$  Thus by Hilbert's theorem 90 every element of $Y^{\circ}(F)$ is in the orbit of an element in $Y_{v_0}^\circ(F).$  We conclude that 
\begin{align*}
&\int_{[\GL_{\underline{1}}]^2  \times [H^e]'}\Theta_f(ag,a'g',h)\eta_{\underline{s},s}(a'I_{\underline{r-2}},h)d^\times ad^\times a'dh
\\
&=\int\sum_{(x,v) \in X^\circ(F) \times \prod_{i=1}^3(F^{r_i-2}-\{0\})} \mathcal{R}_u(ag,a
'g',h)f\left(x,(v,0),v_0\right)\eta_{\underline{s},s}(a'I_{\underline{r-2}},h)d^\times a d^\times a'dh.
\end{align*}
Here the integral is over $[\GL_{\underline{1}}]^2 \times (Z_{\GL_{\underline{2}}}T_HN_0)(F) \backslash H^e(\A_F)'.$ 

By the Iwasawa decomposition of $H^e(\A_F),$ this is 
\begin{align*}
\int_{}&\sum_{(x,v)} \mathcal{R}_u\left(ag,a'g',\begin{psmatrix} \Delta(\lambda) & \\ & 1 \end{psmatrix}b\right)f'\left(x,(v,0),v_0 \right) 
\eta_{\underline{s},s}(a'I_{\underline{r-2}},\begin{psmatrix} \Delta(\lambda) & \\ & 1 \end{psmatrix}b )|\lambda|^{-2}  d^\times ad^\times a'd^\times \lambda d_\ell b
\end{align*}
where $f'=\int_{K_{H^e}}\mathcal{R}_u(I_{\underline{r}},I_{\underline{r-2}},k)fdk$ 
and the integral is over 
$$
[\GL_{\underline{1}}]^2 \times [\GL_1] \times  (Z_{\GL_{\underline{2}}}N_0)(F) \backslash (B^e(\A_F) \cap H^e(\A_F)').
$$ 
Using Lemma \ref{lem:Bisom} we can write the above as $\mathrm{meas}([Z_{\GL_{\underline{2}}}]^1[N_0])$ times 
\begin{align*}
\int&\sum_{(x,v)} 
f'\left(\begin{psmatrix} a'g' & \\ & \begin{psmatrix} 1 & y \\ & 1 \end{psmatrix}\begin{psmatrix}\Delta(\lambda) & \\ & \Delta([c])c^{-1}\end{psmatrix} \end{psmatrix}^{-1}xag,((a'g')^tv\Delta(\lambda)^{-1},0),\left(c_1,c_2,c_3,\frac{[c]}{\lambda}\sum_{i=1}^3y_i\right)\right) \\& \times \eta_{\underline{s},s}(a'I_{\underline{r-2}},\begin{psmatrix} \Delta(\lambda) & \\ & 1 \end{psmatrix}\ell\begin{psmatrix} 1 & \\ & \Delta([c])c^{-1}\end{psmatrix}  ) \frac{|\det a'g'|^{\underline{3/2}}|c|^{\underline{2}}d^\times a d^\times a' d^\times\lambda dyd^\times c}{|\Delta(c)c^{-1}|^{\underline{3(r-2)/4}}|\Delta(\lambda)|^{\underline{r-2}}|\lambda|^2}
\end{align*}
where the integral is over $[\GL_{\underline{1}}]^2  \times [\GL_1] \times N_0 (F) \backslash N_{\underline{2}}(\A_F)\times \GL_{\underline{1}}(\A_F).$  Here we are writing elements of $V(F)$ as tuples $(a,b,c,d)$ with respect to the basis $e_{156},-e_{246},e_{345},e_{456}.$  
   
We now change variables $(a,a') \mapsto (\Delta(\lambda)a,\Delta(\lambda)a')$ and then $\lambda \mapsto \lambda[c],$ and use \eqref{eta:invariance} to see that the above is 
\begin{align*}
\int&\sum_{(x,v)} 
f'\left(\begin{psmatrix} a'g' & \\ &\begin{psmatrix} 1 & y \\ & 1 \end{psmatrix}\begin{psmatrix}1 & \\ & \Delta(\lambda)^{-1}c^{-1}\end{psmatrix} \end{psmatrix}^{-1}xag,(a'g'^tv,0),(c_1,c_2,c_3,\lambda^{-1}\sum_{i=1}^3y_i)\right) \\& \eta_{\underline{s},s}\left(a'I_{\underline{r-2}},\begin{psmatrix} \Delta(\lambda[c]) & \\ & 1 \end{psmatrix}\ell\begin{psmatrix} \Delta(\lambda[c])^{-1} & \\ & \Delta(\lambda^{-1})c^{-1}\end{psmatrix}  \right) \frac{|\Delta(\lambda[c])|^{\underline{(r-2)/2}}|\det a'g'|^{\underline{3/2}}d^\times a d^\times a' d^\times\lambda dyd^\times c}{|\Delta(c)c^{-1}|^{\underline{3(r-2)/4}}|\lambda|^{2}}.
\end{align*}

By our assumptions on $f'$ and Lemma \ref{lem:Xi:norm},
 there are  nonnegative  
$$
(\Phi_1, \Phi_2,\Phi_3)\in \mathcal{S}_{\mathrm{ES}}(X(\A_F))  \times \mathcal{S}(\A_F^{\underline{r-2}}) \times  \mathcal{S}(V(\A_F))
$$
such that the above is bounded by 
\begin{align*}
\int&\sum_{(x,v)} \Xi_{\underline{r-2},\underline{2}}\left(\begin{psmatrix} a'g' & \\ & \begin{psmatrix} 1 & y \\ & 1 \end{psmatrix}\begin{psmatrix}1 & \\ & \Delta(\lambda)^{-1}c^{-1}\end{psmatrix}\end{psmatrix}^{-1}xag\right)
\norm{a'g'^tv}^{\underline{(2-r)/2}}|\Delta(c)c^{-1}|^{\underline{(r-3)/4}}
\\& \times \Phi_1\left(\begin{psmatrix} a'g' & \\ & \begin{psmatrix} 1 & y \\ & 1 \end{psmatrix}\begin{psmatrix}1 & \\ & \Delta(\lambda)^{-1}c^{-1}\end{psmatrix}\end{psmatrix}^{-1}xag\right)
\Phi_2(a'g'^tv)\Phi_3\left(c_1,c_2,c_3,\lambda^{-1}\sum_{i=1}^3y_i\right)\\& \times
\eta_{\underline{s},s}\left(a'I_{\underline{r-2}},\begin{psmatrix} \Delta(\lambda[c]) & \\ & 1 \end{psmatrix}\ell\begin{psmatrix} \Delta(\lambda[c])^{-1} & \\ & \Delta(\lambda^{-1})c^{-1}\end{psmatrix}  \right) \frac{|\Delta(\lambda[c])|^{\underline{(r-2)/2}}|\det a'g'|^{\underline{3/2}}
d^\times a d^\times a' d^\times\lambda dyd^\times c}{|\lambda|^2|\Delta([c])c^{-1}|^{\underline{3(r-2)/4}}}.
\end{align*}
Here we have used \eqref{Pli2:comp}.  Simplifying the above using \eqref{small:comp},  we arrive at 
\begin{align*}
&\int\sum_{(x,v)} \Xi_{{\underline{r-2},\underline{2}}}\left(xg\right)
\Phi_1\left(\begin{psmatrix} a'g' & \\ & \begin{psmatrix} 1 & y \\ & 1 \end{psmatrix}\begin{psmatrix}1 & \\ & \Delta(\lambda)^{-1}c^{-1}\end{psmatrix}\end{psmatrix}^{-1}xag\right)\frac{|\det g'|^{\underline{1/2}}\Phi_2(a'g'^tv)}{\norm{g'^tv}^{\underline{(r-2)/2}}}\\& \times \Phi_3\left(c_1,c_2,c_3,\lambda^{-1}\sum_{i=1}^3y_i\right) 
|\lambda|^{s+\rho-3}|c|^{\underline{s+\rho}}\prod_{i=1}^3|a'_i|^{(r_i-2)C_i(\underline{s},s)+\rho_i'} da da' d^\times\lambda dyd^\times c\\
&=\int_{\GL_{\underline{1}}(\A_F) \times [\GL_{\underline{1}}]^2 \times [\GL_1] \times [N_{\underline{2}}]}\sum_{(x,v)} \Xi_{{\underline{r-2},\underline{2}}}\left(xg\right)
\Phi_1\left(\begin{psmatrix} a'g' & \\ & \begin{psmatrix} 1 &  y \\ & 1 \end{psmatrix}\begin{psmatrix}1 & \\ & \Delta(\lambda)^{-1}c^{-1}\end{psmatrix}\end{psmatrix}^{-1}xag\right) \frac{|\det g'|^{\underline{1/2}}\Phi_2(a'g'^tv)}{\norm{g'^tv}^{\underline{(r-2)/2}}}\\& \times\sum_{\alpha \in F}\Phi_3\left(c_1,c_2,c_3,\lambda^{-1}(\alpha+\sum_{i=1}^3y_i)\right) |\lambda|^{s+\rho-3}|c|^{\underline{s+\rho}}\prod_{i=1}^3|a'_i|^{(r_i-2)C_i(\underline{s},s)+\rho_i'} da da' d^\times\lambda dyd^\times c
\end{align*}
for some $(\rho,\underline{\rho},\underline{\rho}') \in \RR \times \RR^{\underline{1}} \times \RR^{\underline{1}}.$
Replacing the sum over $\alpha$ by an integral using Lemma \ref{lem:FL} and changing variables, to bound the above it suffices to bound 
\begin{align*}
&\int_{\GL_{\underline{1}}(\A_F) \times [\GL_{\underline{1}}]^2 \times [\GL_1] \times [N_{\underline{2}}]}\sum_{(x,v)}\Xi_{{\underline{r-2},\underline{2}}}\left(xg\right)
 \Phi_1\left(\begin{psmatrix} a'g' & \\ & \begin{psmatrix} 1 & y \\ & 1 \end{psmatrix}\begin{psmatrix}1 & \\ & \Delta(\lambda)^{-1}c^{-1}\end{psmatrix}\end{psmatrix}^{-1}xag\right)
 \\& \times\frac{|\det g'|^{\underline{1/2}}\Phi_2(a'g'^tv)\Phi_{3}'(c_1,c_2,c_3)  }{\norm{g'^tv}^{\underline{(r-2)/2}}}
|\lambda|^{s+\rho}|c|^{\underline{s+\rho}}\prod_{i=1}^3| a'_i|^{(r_i-2)C_i(\underline{s},s)+\rho_i'} da da' d^\times\lambda dyd^\times c
\end{align*}
for any $\Phi_3' \in \mathcal{S}(\A_F^3).$ 

Using Lemma \ref{lem:FL} we can bound the integral over $\lambda \in \A_F^\times$ by an integral over $\lambda \in (\A_F^\times)^3$ and see that it suffices to bound
\begin{align*}
&\int_{\GL_{\underline{1}}(\A_F) \times [\GL_{\underline{1}}]^3  \times [N_{\underline{2}}]}\sum_{(x,v)} \Xi_{{\underline{r-2},\underline{2}}}\left(xg\right)
 \Phi_1'\left(\begin{psmatrix} a'g' & \\ & \begin{psmatrix} 1 & y \\ & 1 \end{psmatrix}\begin{psmatrix}1 & \\ & \lambda^{-1}c^{-1}\end{psmatrix}\end{psmatrix}^{-1}xag\right)\\& \times\frac{|\det g'|^{\underline{1/2}}\Phi_2(a'g'^tv)\Phi_{3}'(c_1,c_2,c_3) }{\norm{g'^tv}^{\underline{(r-2)/2}}}
|c|^{\underline{s+\rho}}\prod_{i=1}^3|\lambda_i|^{(s+\rho)/3}|a'_i|^{(r_i-2)C_i(\underline{s},s)+\rho_i'} d^\times a d^\times a' d^\times\lambda dy d^\times c
\end{align*}
for any $\Phi_1' \in \mathcal{S}_{\mathrm{ES}}(X(\A_F)).$
Provided that $s_i-s/3 \gg_{\underline{\rho},\rho} 1,$ we can change variables $\lambda \mapsto c^{-1}\lambda$ and execute the integral over $c$ to see that the above is dominated by
\begin{align*}
\int_{[\GL_{\underline{1}}]^3 \times [N_{\underline{2}}]}&\sum_{(x,v)}\Xi_{{\underline{r-2},\underline{2}}}\left(xg\right)
\Phi_1'\left(\begin{psmatrix} a'g' & \\ & \begin{psmatrix} 1 & y \\ & 1 \end{psmatrix}\begin{psmatrix}1 & \\ & \lambda^{-1}\end{psmatrix}\end{psmatrix}^{-1}xag\right)\\& \times\frac{|\det g'|^{\underline{1/2}}\Phi_2(a'g'^tv)\prod_{i=1}^3|\lambda_i|^{(s+\rho)/3}|a'_i|^{(r_i-2)C_i(\underline{s},s)+\rho_i'}}{\norm{g'^tv}^{\underline{(r-2)/2}}} d^\times a d^\times a' dy d^\times \lambda.
\end{align*}

Using Lemma \ref{lem:FL} again, we see that to bound this sum it suffices to bound
\begin{align*}
&\int\sum_{x} \Xi_{{\underline{r-2},\underline{2}}}\left(xg\right)
 \Phi_1'\left(\begin{psmatrix} a'g' & \\ & \begin{psmatrix} 1 & y \\ & 1 \end{psmatrix}\begin{psmatrix}1 & \\ & \lambda^{-1}\end{psmatrix}\end{psmatrix}^{-1}xag\right)\\& \hskip.2in \times \sum_{v \in \mathbb{P}^{\underline{r-3}}(F)}\left(\int_{\GL_{\underline{1}}(\A_F)}|a_1|^{\underline{z}}
 \frac{|\det g'|^{\underline{1/2}}\Phi_2'(a'a_1g'^tv)\prod_{i=1}^3|\lambda_i|^{(s+\rho)/3}|a'_i|^{(r_i-2)C_i(\underline{s},s)+\rho_i'}}{\norm{g'^tv}^{\underline{(r-2)/2}}}
 d^\times a_1\right)  
 d^\times a d^\times a' d^\times\lambda dy\\
&=\int \sum_{x} \Xi_{{\underline{r},\underline{r-2}}}\left(xg\right)\Phi_1\left(\begin{psmatrix} a'g' & \\ & \begin{psmatrix} 1 & y \\ & 1 \end{psmatrix}\begin{psmatrix}1 & \\ & \lambda^{-1}\end{psmatrix}\end{psmatrix}^{-1}xag\right)
\\& \hskip.2in \times \sum_{v \in \mathbb{P}^{\underline{r-3}}(F)}\left(\int_{\GL_{\underline{1}}(\A_F)}|a_1|^{\underline{z}}|\det g'|^{\underline{1/2}}\norm{g'^tv}^{\underline{(2-r)/2}}\Phi_2'(a_1g'^tv)d^\times a_1\right) \\& \hskip.2in \times
\prod_{i=1}^3|\lambda_i|^{(s+\rho)/3}|a'_i|^{(r_i-2)C_i(\underline{s},s)+\rho_i'-z_i} d^\times a d^\times a' d^\times\lambda dy
\end{align*}
for $\underline{z} \in \RR^{\underline{1}}$ sufficiently large and all $\Phi_2' \in \mathcal{S}(\A_F^{\underline{r-2}}).$ 
 Here the unmarked integrals are over $[\GL_{\underline{1}}]^3  \times [N_{\underline{2}}].$

  If we view the quantity in parenthesis as a function of $g',$ it is a mirabolic Eisenstein series.  It is well-known to converge for $z_i \gg_{r_i} 1$  and defines a function of moderate growth when it converges \cite[Lemma 4.2]{JacquetShalikaEPI}.
Thus we are left with proving that 
\begin{align*} 
\int&\sum_{x} \Xi_{{\underline{r-2},\underline{2}}}\left(xg\right)
\Phi_1\left(\begin{psmatrix} a'g' & \\ & \begin{psmatrix} 1 & y \\ & 1 \end{psmatrix}\begin{psmatrix}1 & \\ & \lambda^{-1}\end{psmatrix}\end{psmatrix}^{-1}xag\right)
\prod_{i=1}^3|\lambda_i|^{(s+\rho)/3}|a'_i|^{(r_i-2)C_i(\underline{s},s)+\rho_i'-z_i} d^\times a d^\times a' d^\times\lambda dy
\end{align*}
converges and defines a function of moderate growth in $(g',g) \in \GL_{\underline{r-2}}(\A_F) \times \GL_{\underline{r}}(\A_F)$ provided  $z_i \gg_{r_i} 1$ and $\underline{s},s$ satisfy the conditions in the statement of the lemma.  Here the integral is over $[\GL_{\underline{1}}]^3 \times [N_{\underline{2}}].$

We rewrite the sum as 
\begin{align*}
\int_{[\GL_{\underline{1}}]^3 }&\sum_{x' } \Xi_{{\underline{r-2},\underline{2}}}\left(x'g\right)
\int_{\A_F}\Phi_1\left(\begin{psmatrix} (a'g')^{-1} & \\ & 1 & \\ & &\lambda^{-1}\end{psmatrix}\begin{psmatrix} I_{\underline{r-2}} & & \\ & 1 & -y \\ & & 1 \end{psmatrix}x'ag\right)dy \\& \times
\prod_{i=1}^3|\lambda_i|^{(s+\rho)/3}|a'_i|^{(r_i-2)C_i(\underline{s},s)+\rho_i'-z_i} d^\times a d^\times a' d^\times\lambda 
\end{align*}
where the sum is now over $x'\in (N_{\underline{r-2},\underline{1},\underline{1}} \backslash \GL_{\underline{r}})(F).$  
At this point we apply Lemma \ref{lem:factor:mod:growth} to see that this converges and defines a function of moderate growth provided that $s \gg  1,$ $(r_i-1)C_i(\underline{s},s)-\mathrm{Re}(s)/3 \gg_{z_i} 1.$  
\end{proof}

Applying Lemma \ref{lem:Sobolev} to Lemma \ref{lem:abs:conv0}, we deduce the following lemma:
\begin{lem} \label{lem:mod:nos}
The integral $\int_{[\GL_{\underline{1}}] \times [H^e]^1}\Theta_{f}(ag,g',h)d^\times adh$ is absolutely convergent and is of moderate growth as a function of  $(g,g') \in \GL_{r}(\A_F) \times \GL_{r-2}(\A_F).$  
\qed
\end{lem}

\section{Unfolding the global zeta integrals} \label{sec:unfold}

Let $\pi=\pi_1 \otimes \pi_2\otimes \pi_3$ (resp. $\pi'=\pi_1' \otimes \pi_2' \otimes \pi_3'$) be a cuspidal automorphic representation of $\GL_{\underline{r}}(\A_F)$ (resp. $\GL_{\underline{r-2}}(\A_F)).$  We let $\omega,\omega_i,\omega',\omega_{i}'$ be the central characters of 
$\pi$, $\pi_i$, $\pi'$ and $\pi_i'$, respectively.  We assume that $\omega=\omega'^{-1}.$   Let $\varphi$ be an automorphic form in the space of $\pi$ and let $\varphi'$ be an automorphic form in the space of $\pi'.$ 

In \eqref{glob:zeta}, we defined
\begin{align} \begin{split}    &Z(\varphi,\varphi',f,\underline{s},s):=\int_{[G]'}\Theta_f(g)\varphi(g^{(1)})\varphi'(g^{(2)})\eta_{\underline{s},s}(g)dg.
    \end{split}
\end{align} 
Conditions ensuring absolute convergence are given in Lemma \ref{lem:abs:conv0} above.  In this section we prove our unfolding result, Theorem \ref{thm:intro1}, restated as Theorem \ref{thm:unfold} below.  
  
Let 
\begin{align} \label{mira}
    \mathcal{P}_{r_i}(R):=\{ \begin{psmatrix} g & x\\ & 1\end{psmatrix}: g \in \GL_{r_i-1}(R) \times R^{r_i-1}\}
\end{align}
be the usual mirabolic subgroup and let $\mathcal{P}_{\underline{r}}:=\prod_{i=1}^3\mathcal{P}_{r_i}.$  By convention $\mathcal{P}_1$ is the trivial group.

\begin{lem}\label{lem:unfold} 
We have that
\begin{align*}
&\int_{[\GG_m]^1 \times [\mathcal{P}_{\underline{r-2}}]^1}\int_{[N_{\underline{r-2},\underline{2}}] \times [N_{0}]} \overline{\psi}\left( \langle m_0,(u_1,u_2)\rangle \right)\varphi\begin{psmatrix} \Delta(\lambda)p & \Delta(\lambda) u_1 & u_2\\ & \Delta(\lambda) & t \\ & & 1\end{psmatrix} du_1du_2dt \varphi'(p) dp d^\times \lambda 
\\&=\int_{(\A_F^\times)^1 \times U_{\underline{r-2}}(\A_F) \backslash \mathcal{P}_{\underline{r-2}}(\A_F)^1} W^\varphi_{\psi}\begin{psmatrix}  \Delta(\lambda) p & & \\ & \Delta(\lambda) & \\ & & 1\end{psmatrix} W^{\varphi'}_{\overline{\psi}}(p)d^\times \lambda d\dot{p}.
\end{align*}
\end{lem}

\begin{proof} Let $N_2 \leq \SL_2$ be the unipotent radical of the Borel subgroup of upper triangular matrices.
By Fourier inversion, we have
\begin{align*}
    &\int_{[\GG_a^{\underline{r-2}}] \times [N_0]}
    \varphi\begin{psmatrix} \Delta(\lambda)p & \Delta(\lambda) u_1 & u_2\\ & \Delta(\lambda) & t \\ & & 1\end{psmatrix} du_2dt\\
    &=\sum_{\alpha\in F}\int_{[\GG_a^{\underline{r-2}}] \times [N_{\underline{2}}]}\varphi\left(\begin{psmatrix} I_{\underline{r-2}} & & u_2 \\ & 1 & t \\ & & 1\end{psmatrix} \begin{psmatrix} \Delta(\lambda)p & \Delta(\lambda) u_1 & \\ & \Delta(\lambda) &  \\ & & 1\end{psmatrix}\right) \psi\left(-\alpha \sum_{i=1}^3t_i\right)du_2 dt\\
    &=\sum_{\alpha\in F^\times }\int_{[\GG_a^{\underline{r-2}}] \times [N_{\underline{2}}]}\varphi\left( \begin{psmatrix} I_{\underline{r-2}} & & u_2\\ & 1 & t\\ & & 1\end{psmatrix} \begin{psmatrix} \Delta(\lambda \alpha )p & \Delta(\lambda \alpha ) u_1 & \\ & \Delta(\lambda \alpha) &  \\ & & 1\end{psmatrix}\right)\psi\left(-\sum_{i=1}^3 t\right)du_2 dt.
\end{align*}
Here in the last equality we have used the fact that the contribution of $\alpha=0$ is zero because $\varphi$ is a cusp form, a change of variables $(u_2,t) \mapsto (\Delta(\alpha)^{-1}u_2,\Delta(\alpha)^{-1}t),$ and the fact that $\varphi$ is left $\GL_{\underline{r}}(F)$-invariant.    The sum is absolutely convergent because $\varphi$ is smooth and the integral is absolutely convergent because it has compact domain.   The last integral is nothing but the projection of $\varphi$ to  $\mathcal{P}_{\underline{r-1}}(\A_F)$ in the sense of \cite[\S 5.2]{CogdellFields}.  By applying the usual Whittaker expansion (see \cite[Proposition 5.2]{CogdellFields}) to the last expression in the displayed equation above, we see that 
\begin{align} \label{first:Id} \begin{split}
\int_{[\GG_a^{\underline{r-2}}] \times [N_0]}
    \varphi\begin{psmatrix} \Delta(\lambda)p & \Delta(\lambda) u_1 & u_2\\ & \Delta(\lambda) & t \\ & & 1\end{psmatrix} du_2dt
    &=\sum_{\alpha\in F^\times }\sum_{\gamma \in U_{\underline{r-2}}(F) \backslash \GL_{\underline{r-2}}(F)}W^{\varphi}_{\psi}\begin{psmatrix} \gamma \Delta(\alpha\lambda )p & \gamma \Delta(\alpha\lambda ) u_1 & \\ & \Delta(\alpha \lambda) &  \\ & & 1\end{psmatrix}. \end{split}
\end{align}
The double sum here converges absolutely because the full Whittaker expansion of a cusp form converges absolutely. 

Using the fact that $\psi\left( \langle m_0,(u_1,u_2)\rangle \right)=\psi\left( \langle m_{0},(u_1,0)\rangle \right),$ the identity \eqref{first:Id}
implies that
\begin{align} \label{before:abs} \begin{split}
    &\int_{[N_{\underline{r-2},\underline{2}}] \times [N_{0}]} \overline{\psi}\left( \langle m_0,(u_1,u_2)\rangle \right)\varphi\begin{psmatrix} \Delta(\lambda)p & \Delta(\lambda) u_1 & u_2\\ & \Delta(\lambda) & t \\ & & 1\end{psmatrix} du_1du_2dt \\
   &=\int_{[\GG_a^{\underline{r-2}}]}\overline{\psi}\left( \langle m_0,(u_1,0)\rangle \right)\left(\int_{[\GG_a^{\underline{r-2}}] \times [N_{0}]} \varphi\begin{psmatrix} \Delta(\lambda)p & \Delta(\lambda) u_1 & u_2\\ & \Delta(\lambda) & t \\ & & 1\end{psmatrix} dtdu_2\right)du_1\\
    &=\int_{[\GG_a^{\underline{r-2}}]}\overline{\psi}\left( \langle m_{0},(u_1,0)\rangle \right)\sum_{(\alpha,\gamma)\in F^\times \times U_{\underline{r-2}}(F)\backslash \GL_{\underline{r-2}}(F)} W^\varphi_\psi\begin{psmatrix} \gamma  \Delta(\alpha \lambda)p & \gamma\Delta(\alpha \lambda) u_1 & \\ & \Delta(\alpha \lambda) &  \\ & & 1\end{psmatrix}du_1. \end{split}
\end{align}

 Consider the function
\begin{align*}
\GL_{\underline{r-1}}(\A_F)^1 \times \GL_{\underline{r-2}}(\A_F)^1 &\lto \CC\\
(g,g') &\longmapsto \sum_{\gamma \in U_{\underline{r-1}}(F) \backslash \GL_{\underline{r-1}}(F)}\left|W^{\varphi}_{\psi}\begin{psmatrix} \gamma g& \\  & 1 \end{psmatrix} \right||\varphi'(g')|.
\end{align*}
It is rapidly decreasing by \cite[\S 4]{Cogdell:PS:PPS}. Hence by reduction theory
\begin{align*}
    \int_{[\GG_m]^1 \times [\mathcal{P}_{\underline{r-2}}]^1} \int_{[\GG_{a}^{\underline{r-2}}]}\sum_{\alpha\in F^\times }\sum_{\gamma \in U_{\underline{r-2}}(F) \backslash \GL_{\underline{r-2}}(F)}\left|W^{\varphi}_{\psi}\begin{psmatrix} \gamma \Delta(\alpha\lambda )p & \gamma\Delta(\alpha\lambda ) u_1 & \\ & \Delta(\alpha \lambda) &  \\ & & 1\end{psmatrix} \varphi'(p)\right|du_1dpd^\times \lambda<\infty.
\end{align*}
Thus we can unfold the sum over $\alpha$ and bring the integral over $[\GG_a^{\underline{r-2}}]$ inside the sum over $\gamma$ to see that
\begin{align*}
&\int_{[\GG_m]^1 \times [\mathcal{P}_{\underline{r-2}}]^1}\int_{[\GG_a^{\underline{r-2}}]}\overline{\psi}\left( \langle m_{0},(u_1,0)\rangle \right)\sum_{(\alpha,\gamma)\in F^\times \times U_{\underline{r-2}}(F)\backslash \GL_{\underline{r-2}}(F)} W^\varphi_\psi\begin{psmatrix} \gamma \Delta(\alpha \lambda)p & \gamma\Delta(\alpha \lambda) u_1 & \\ & \Delta(\alpha \lambda) &  \\ & & 1\end{psmatrix}du_1 \varphi'(p)dpd^\times \lambda\\
&=\int_{(\A_F^\times)^1}\int_{[\mathcal{P}_{\underline{r-2}}]^1}\sum_{\gamma\in U_{\underline{r-2}}(F)\backslash \mathcal{P}_{\underline{r-2}}(F)} W^\varphi_\psi\begin{psmatrix} \gamma\Delta(\lambda)p &  & \\ & \Delta(\lambda) &  \\ & & 1\end{psmatrix}\varphi'(p)dpd^\times \lambda \\
    &=\int_{(\A_F^\times)^1}\int_{U_{\underline{r-2}}(\mathbb{A}_F)\backslash\mathcal{P}_{\underline{r-2}}(\mathbb{A}_F)^1} \int_{[U_{\underline{r-2}}]}W^\varphi_\psi \begin{psmatrix}  \Delta(\lambda) p & & \\ & \Delta(\lambda) & \\ & & 1\end{psmatrix}\varphi'(np)\psi(n) dn d\dot{p}d^\times \lambda\\
    &= \int_{(\A_F^\times)^1}\int_{U_{\underline{r-2}}(\A_F) \backslash \mathcal{P}_{\underline{r-2}}(\A_F)^1} W^\varphi_{\psi}\begin{psmatrix}  \Delta(\lambda)p & & \\ & \Delta(\lambda) & \\ & & 1\end{psmatrix} W^{\varphi'}_{\overline{\psi}}(p) d\dot{p}d^\times \lambda.
\end{align*}
Combining this with \eqref{before:abs} we deduce the result.
\end{proof}

\begin{thm}\label{thm:unfold}
Let $f$ satisfy conditions \ref{ad1}, \ref{ad2} and \ref{ad3} of \S \ref{sec:conv}. The integral $Z(\varphi,\varphi',f,\underline{s},s)$ is absolutely convergent for $\mathrm{Re}(\underline{s}),\mathrm{Re}(s)$ satisfying \eqref{bounds2}
and is equal to $\mathrm{meas}(Z_{\GL_{\underline{2}}}(F) \backslash Z_{\GL_{\underline{2}}}(\A_F)^1)$ times the absolutely convergent integral 
\begin{align*}
    &\int_{Z_{\GL_{\underline{2}}}(\A_F)N_0(\A_F) \backslash H^e(\A_F)}\int_{U_{\underline{r-2}}(\A_F) \backslash \GL_{\underline{r-2}}(\A_F)}\int W^{\varphi}_{\psi}(g)W^{\varphi'}_{\overline{\psi}}(g')\mathcal{R}_u\left(g,g',h\right)f(x_0) \eta_{\underline{s},s}(g',h)dg dg'dh,
\end{align*}
where the inner integral is over
$N_{\underline{r-2},\underline{2}}(\A_F) \backslash \GL_{\underline{r}}(\A_F).$
\end{thm}

\begin{proof} The integral $Z(\varphi,\varphi',f,\underline{s},s)$ converges absolutely for $\mathrm{Re}(\underline{s}), \mathrm{Re}(s)$ as in the theorem by Lemma \ref{lem:abs:conv0}.  Consider the following subgroups of $G$
\begin{align*}
\widetilde{Z}(R):&=\{(aI_{\underline{r}},aI_{\underline{r-2}},aI_{\underline{2}}):a \in (R^\times)^{\underline{1}}\},\\
    M(R):&=\left\{\left(\begin{psmatrix}  \Delta(\lambda)p & & \\ &  \Delta(\lambda) &\\ &  & 1 \end{psmatrix},   \Delta(\lambda) p
    ,\begin{psmatrix}\Delta(\lambda) & \\ & 1 \end{psmatrix}\right): (p,\lambda) \in   \mathcal{P}_{\underline{r-2}}(R)  \times R^\times\right\},\\
    N(R):&=\left\{\left(\begin{psmatrix} I_{\underline{r-2}} & m\\ & n_0\end{psmatrix}, I_{\underline{r-2}},n_0\right): (m,n_0) \in M_{\underline{r-2},\underline{2}}(R) \times N_0(R)\right\}.
\end{align*}
Lemma \ref{lem:orb:stab2} asserts that $G_{x_0}=\widetilde{Z}MN.$ In particular, by Hilbert's theorem 90, $X^\circ(F)$ is the $G(F)$-orbit of $x_0.$
We have 
\begin{align*}
&Z(\varphi,\varphi',f,\underline{s},s)\\&=\int_{[G]' }\sum_{x\in X^{\circ}(F)}\mathcal{R}_u(g)f(x) \varphi(g^{(1)})\varphi'(g^{(2)})\eta_{\underline{s},s}(g^{(2)},g^{(3)})dg \\
&=\int_{G_{x_0}(F) \backslash G(\A_F)'}\mathcal{R}_u(g)f(x_0)\varphi(g^{(1)})\varphi'(g^{(2)})\eta_{\underline{s},s}(g^{(2)},g^{(3)})dg\\
&=\int \mathcal{R}_u(g)f(x_0) \eta_{\underline{s},s}(g^{(2)},g^{(3)})\\& \times\iint \psi\left( \langle m_0,(u_1,u_2)\rangle \right)\varphi\left(\begin{psmatrix} \Delta(\lambda)p & \Delta(\lambda)u_1 & u_2\\ & \Delta(\lambda) & t \\ & & 1\end{psmatrix}mg^{(1)}\right)du_1du_2dt \varphi'(\Delta(\lambda)pg^{(2)}) d^\times \lambda dp d\dot{g}
\end{align*}
where the integrals, from left to right, are over 
$$
\left(g,\lambda,p,\begin{psmatrix} I_{\underline{r-2}} & u_1 & u_2 \\ & 1& \\ & & 1 \end{psmatrix},\begin{psmatrix} 1 & t \\ & 1\end{psmatrix}\right)
\in \widetilde{Z}(F)M(\A_F)^1 N(\A_F)\backslash G(\A_F)' \times [\GG_m]^1 \times [\mathcal{P}_{\underline{r-2}}]^1 \times [N_{\underline{r-2},\underline{2}}] \times [N_0].
$$

Applying Lemma \ref{lem:unfold}, the above is
\begin{align*}
    \int \mathcal{R}_u(g)f(x_0)\eta_{\underline{s},s}(g^{(2)},g^{(3)})\int W^\varphi_{\psi}\left(\begin{psmatrix}  \Delta(\lambda) p & & \\ & \Delta(\lambda) & \\ & & 1\end{psmatrix} g^{(1)}\right)W^{\varphi'}_{\overline{\psi}}(\Delta(\lambda)pg^{(2)})d^\times \lambda d\dot{p} d\dot{g}
\end{align*}
where the outer integral is as above and the inner integral is now over $(\A_F^\times)^1 \times U_{\underline{r-2}}(\A_F) \backslash \mathcal{P}_{\underline{r-2}}(\A_F)^1.$  
We rearrange the order of integration to obtain the expression in the theorem.  This is justified since
the integral in the theorem
converges absolutely for $\mathrm{Re}(\underline{s}),\mathrm{Re}(s)$ satisfying \eqref{bounds2} by the local bounds in Proposition \ref{prop:ram:conv} and Lemma \ref{lem:unr:conv}.
\end{proof}

\section{Applying the Poisson summation conjecture}
 \label{sec:apply} 
In this section we explain how the Poisson summation formula of Conjecture \ref{PS:conj} implies the analytic continuation and functional equation of $Z(\varphi,\varphi',f,\underline{s},s),$ and in turn the analytic continuation of $L^S(s,\pi,\otimes^3).$ 
The idea of the argument is familiar and goes back to Tate's thesis (or even Riemann).  Nevertheless, our setting is nonstandard, and in particular involves a multivariable zeta function.  Thus it warrants a full discussion.  This also serves as a check that various normalizations are consistent.  

 Let $\pi$ and $\pi'$ be cuspidal automorphic representations of $A_{\GL_{\underline{r}}} \backslash \GL_{\underline{r}}(\A_F)$ and $A_{\GL_{\underline{r-2}}} \backslash \GL_{\underline{r-2}}(\A_F),$ respectively, and let $\varphi$ (resp.~$\varphi'$) be a cuspidal form in the space of $\pi$ (resp.~$\pi'$).  We assume that the central characters $\omega$ of $\pi$ and $\omega'$ of $\pi'$ satisfy $\omega \omega'=1.$

\begin{thm} \label{thm:func:eqn}
    Assume Conjecture \ref{PS:conj}.  Let $v_1$ and $v_2$ be two finite places of $F$, and let $f=f_{v_1}f_{v_2}f^{v_1v_2} \in \mathcal{S}(
    X(\A_F),\mathcal{H})$ where $f_{v_1} \in \mathcal{S}(X^\circ(F_{v_1}),\mathcal{H})$ and $\mathcal{F}(f_{v_2}) \in \mathcal{S}(X^\circ(F_{v_2}),\mathcal{H}).$   Assume that $\varphi,$ $\varphi',$ and $f$ are pure tensors.  The zeta function $Z(\varphi, \varphi',f,\underline{s},s)$ admits a holomorphic continuation to $ \CC^{\underline{1}} \times \CC$.  Moreover,
$Z(\varphi,\varphi',f,\underline{s},s)=Z(\varphi^\vee,\varphi'^\vee,\mathcal{F}(f),-\underline{s},-s),$
where $\varphi^\vee(g):=\varphi(g^{-t})$ and $\varphi'^\vee(g):=\varphi'(g^{-t}).$
\end{thm}
\noindent We point out that for general $f \in \mathcal{S}(
    X(\A_F),\mathcal{H}),$ the function $Z(\varphi,\varphi',f,\underline{s},s)$ will have poles.  The local conditions on $f$ in the theorem eliminate these poles.  

Let $\delta_{ij}$ be the Kronecker $\delta$-function and let 
\begin{align*}
\kappa'_i:\RR_{>0} &\lto  A_{\GL_{\underline{r-2}}} \times A_{H}\\
a &\longmapsto \left( \left( a^{(\delta_{i1}+1)/2}I_{r_1-2},a^{(\delta_{i2}+1)/2}I_{r_2-2},a^{(\delta_{i3}+1)/2}I_{r_3-2}\right),\Delta(a^{1/2})I_H  \right),\\
\kappa:\RR_{>0} &\lto   A_{\GL_{\underline{r-2}}} \times A_{H}\\
a &\longmapsto (\Delta(a)I_{\underline{r-2}},\Delta(a^{1/2})I_{H}).
\end{align*}
We moreover define $\kappa':=\prod_{i=1}^3\kappa_i':\RR_{>0}^{\underline{1}} \to A_{\GL_{\underline{r-2}}} \times A_{H}.$  Then by \eqref{CA},
\begin{align} \label{etass'}
\eta_{\underline{s},s}(\kappa'(a')\kappa(a))=a'^{\underline{s}}a^s
\end{align}
where $a'^{\underline{s}}:=\prod_{i=1}^3a'^{s_i}_i.$  In particular
\begin{align} \label{is:central:iso}
    \kappa'\kappa:\RR^1_{\underline{>0}} \times \RR_{>0} \lto A_{\GL_{\underline{r-2}}} \times A_{H}
\end{align}
is an isomorphism.

For $(a',a) \in \RR_{>0}^{\underline{1}} \times \RR_{>0},$ let
\begin{align} \label{Zaa} \begin{split}
&Z_{a',a}(\varphi,\varphi',f)=\int_{[\GL_{\underline{r}}] \times [\GL_{\underline{r-2}} \times H^e]^1} \Theta_f(\kappa'(a')\kappa(a)(g,g',h))\varphi(g)\varphi'(g')dgdg'dh. \end{split}
\end{align}
This integral is absolutely convergent by Lemma \ref{lem:mod:nos}.  Since \eqref{is:central:iso} is an isomorphism, \eqref{etass'} implies that for $(\underline{s},s)$ as in \eqref{bounds2},
we have
\begin{align} \label{MT}
\int_{\RR_{>0}^{\underline{1}} \times \RR_{>0}}Z_{a',a}(\varphi,\varphi',f)a^{s} a'^{\underline{s}}d^\times a' d^\times a=Z(\varphi,\varphi',f,\underline{s},s)
\end{align}
and the integral converges absolutely.  

\begin{prop} \label{prop:first:fe} Under the assumptions of Theorem \ref{thm:func:eqn}, we have
\begin{align*}
    Z_{a',a}(\varphi,\varphi',f)=Z_{a'^{-1},a^{-1}}(\varphi^\vee,\varphi'^\vee,\mathcal{F}(f)).
\end{align*}
\end{prop}

\begin{proof}
We have 
\begin{align*}
&Z_{a',a}(\varphi,\varphi',f)=\int_{}\sum_{x \in X^\circ(F)}\mathcal{R}_u(\kappa'(a')\kappa(a)(g,g',h))f(x)\varphi(g)\varphi'(g')dg dg'dh,
\end{align*}
where the integral is over $[\GL_{\underline{r}}] \times [\GL_{\underline{r-2}} \times H^e]^1.$
By Conjecture \ref{PS:conj} this is
\begin{align*}
    \int \sum_{x \in X^\circ(F)}\mathcal{F}(\mathcal{R}_u(\kappa'(a')\kappa(a)(g,g',h))f)(x)\varphi(g)\varphi'(g)dg dg'dh.
\end{align*}
Applying \ref{FT2} this is
\begin{align*}
    \int \sum_{x \in X^\circ(F)}\mathcal{R}_u(\kappa'(a'^{-1})\kappa(a^{-1})(g^{-t},g'^{-t},h^{\iota}))\mathcal{F}(f)(x)\varphi(g)\varphi'(g)dg dg'dh.
\end{align*}
We change variables $(g,g',h) \longmapsto (g^{-t},g'^{-t},h^{\iota})$ to deduce the proposition.
\end{proof}

The following lemma is a consequence of Lemma \ref{lem:Sobolev} and Lemma \ref{lem:abs:conv0}:

\begin{lem} \label{lem:conv:fixed:a}
Let $a \in \RR_{>0}$ and $\underline{s} \in \CC^{\underline{1}}.$  Assume that there exists an $s \in \CC$ such that $(\underline{s},s)$ satisfies \eqref{bounds2}.  Then the integral
\begin{align*}
Z_a(\varphi,\varphi',f,\underline{s}):=\int_{\RR_{>0}^{\underline{1}}} Z_{a',a}(\varphi,\varphi',f)\prod_{i=1}^3a_i'^{s_i}d^\times a'
\end{align*}
converges absolutely. \qed
\end{lem}

Let $N_1 \in \ZZ_{\geq 1}$ and let 
$$
\mathcal{T}_{N_1}:=\{\underline{s} \in \CC^{\underline{1}}:-N_1 <\mathrm{Re}(s_i)<N_1\}.
$$

\begin{prop} \label{prop:first:ac} Assume the hypotheses of Theorem \ref{thm:func:eqn}. For each $a \in \RR_{>0},$ the function $Z_a(\varphi,\varphi',f,\underline{s})$  admits a holomorphic continuation to $\CC^{\underline{1}}$.   Moreover, for each $p \in \CC[\underline{s}]$  one has that
\begin{align} \label{p:bound}
\sup_{\underline{s} \in \mathcal{T}_{N_1}}|p(\underline{s})Z_{a}(\varphi,\varphi',f,\underline{s})| \ll_{N_1,\varphi,\varphi',f,a,p} 1.
\end{align}
Let $N_2 \in \ZZ_{ \geq 0}.$  
If $a \geq 1,$ then we have the refined estimate
\begin{align} \label{p:bound2}
\sup_{\underline{s} \in \mathcal{T}_{N_1}}|p(\underline{s})Z_{a}(\varphi,\varphi',f,\underline{s})| \ll_{N_1,N_2,\varphi,\varphi',f,p} a^{-N_2}.
\end{align}
\end{prop}

\begin{proof}
Let $S$ be a finite set of places of $F$ including the infinite places, the places dividing $6$, the places lying above places of $\QQ$ ramified in $F/\QQ,$ and the places where $\psi$ is ramified.  Upon enlarging $S$ if necessary, we may assume $\varphi \otimes \varphi'$ is fixed by $\GL_{\underline{r}}(\widehat{\OO}_F^S) \times \GL_{\underline{r-2}}(\widehat{\OO}_F^S).$  If $(\underline{s},s)$ satisfy \eqref{bounds2}, Theorem \ref{thm:unfold} and \ref{zeta:basic} imply that 
\begin{align*} 
&Z(\varphi,\varphi',f,\underline{s},s)\\&=\left(\prod_{v \in S}Z(W_v,W_v',f_v,\underline{s},s)\right) L^S(s+\tfrac{1}{2},\pi^\vee,\otimes^3)\prod_{i=1}^3L^S(C_i(\underline{s},s)+\tfrac{1}{2},\pi_i \times \pi_i')L^S(s_i+\tfrac{1}{2},\pi_i^\vee)
\end{align*}
for suitable $(W_{v},W_v') \in \mathcal{W}(\pi_v,\psi) \times \mathcal{W}(\pi'_v,\overline{\psi}).$   

By Mellin inversion, we have 
\begin{align} \label{by:Mell}
Z_a(\varphi,\varphi',f,\underline{s})=\frac{1}{2\pi i}\int_{\mathrm{Re}(s)=\sigma} a^{-s}Z(\varphi,\varphi',f,\underline{s},s)ds
\end{align}
provided that $\sigma=\mathrm{Re}(s)$ and $\mathrm{Re}(s_i)$ satisfy \eqref{bounds2}.  Here we are using the fact that in this range $Z(\varphi,\varphi',f,\underline{s},s)$ decays rapidly as a function of $s$ in vertical strips by \ref{nA:ratio}, \ref{Arch:ratio} and \ref{zeta:basic}.

Let $\overline{\mathcal{T}}_{N_1}$ be the closure of $\mathcal{T}_{N_1}$ in $\CC^{\underline{1}}.$ Let $p \in \CC[\underline{s}]$ be chosen so that $p(\underline{s})\prod_{i=1}^3 L(s_i+\frac{1}{2},\pi_{i\infty}^\vee)$ has no poles in $\overline{\mathcal{T}}_{N_1}.$  Under the assumptions on $\sigma$ and $\mathrm{Re}(\underline{s})$ above, \eqref{by:Mell} implies that
\begin{align}\label{Mellin} \begin{split}
\frac{Z_a(\varphi,\varphi',f,\underline{s})p(\underline{s})}{\prod_{i=1}^3L(s_i+\tfrac{1}{2},\pi_i^{\infty \vee})}&=\frac{1}{2\pi i}\int_{\mathrm{Re}(s)=\sigma}\frac{p(\underline{s})}{a^s}\left(\prod_{v|\infty}Z(W_v,W_v',f_v,\underline{s},s)\right)\left(\prod_{v \in S-\infty }\frac{Z(W_v,W_v',f_v,\underline{s},s)}{\prod_{i=1}^3L(s_i+\tfrac{1}{2},\pi_{iv}^\vee)}\right)\\& \hskip1in \times L^S(s+\tfrac{1}{2},\pi^\vee,\otimes^3)\prod_{i=1}^3L^S(C_i(\underline{s},s)+\tfrac{1}{2},\pi_i \times \pi_i')ds. \end{split}
\end{align}
We claim this admits a holomorphic continuation to $\underline{s} \in \mathcal{T}_{N_1}$ that is bounded by $O_{N_1,\varphi,\varphi',f,p,a}(1)$.  We claim moreover that for $a \geq 1$ and $\underline{s} \in \mathcal{T}_{N_1}$, it is bounded by $O_{N_1,N_2,\varphi,\varphi,f,p}(a^{-N_2})$ for any $N_2 \in \ZZ_{\geq 0}.$

First we point out that \ref{nA:ratio}, \ref{Arch:ratio} and \ref{zeta:basic} imply that \eqref{Mellin} converges absolutely for $\mathrm{Re}(s) \gg 1$ and $\mathrm{Re}(C_i(\underline{s},s)) \gg 1.$  In this range it is equal to a finite sum of expressions of the form 
\begin{align} \label{for:MI} \begin{split}
\frac{1}{2\pi i}\int_{\mathrm{Re}(s)=\sigma}& a^{-s}\left(p(\underline{s})\prod_{v|\infty}Z(W_v,W_v',f_v,\underline{s},s)\right)\\& \gamma^{s}\gamma'^{\underline{s}}\sum_{n_0=1}^\infty \frac{b_0(n_0)}{n_0^{s+1/2}} \prod_{i=1}^3\left(\sum_{n_i=1}^\infty \frac{b_i(n_i)}{n_i^{(r_i-2)^{-1}(s+s_i-s_{i+1}-s_{i+2})+1/2}}\right) ds
\end{split}
\end{align}
for some $\gamma \in \QQ^\times$ and $\gamma' \in (\QQ^\times)^{\underline{1}}$ and $b_0(n_0), b_i(n_i) \in \CC.$  Here the indices over $i$ are taken modulo $3$ in the obvious sense.  
If we define $\phi_{p,\underline{s}}$ as in Lemma \ref{lem:MI} then
 \eqref{for:MI} is
\begin{align} \label{Mell3}
&\gamma'^{\underline{s}}\sum_{n_0=1}^\infty \frac{b_0(n_0)}{n_0^{1/2}}\prod_{i=1}^3\left(\sum_{n_i=1}^\infty \frac{b_i(n_i)}{n_i^{(r_i-2)^{-1}(s_i-s_{i+1}-s_{i+2})+1/2}}\right)\phi_{p,\underline{s}}\left(\frac{an_0n_1^{(r_i-2)^{-1}}n_2^{(r_i-2)^{-1}}n_3^{(r_i-2)^{-1}}}{\gamma}\right).
\end{align}
By Lemma \ref{lem:MI} the function $\phi_{p,\underline{s}}(a)$ is holomorphic for $\underline{s} \in \mathcal{T}_{N_1}$ and rapidly decreasing as a function of $a \geq 1.$  The claim follows.

We now explain how to deduce the proposition from the claim.  Let $p \in \CC[\underline{s}]$ be a polynomial such that $p(\underline{s})\prod_{i=1}^3L(s_i+\tfrac{1}{2},\pi_{i\infty}^\vee)$ has no poles or zeros in $\overline{\mathcal{T}}_{N_1}.$  
Since $L(s_i+\tfrac{1}{2},\pi_i^\vee)$ is holomorphic we deduce that $p(\underline{s})^{-1}\prod_{i=1}^3L(s_i+\tfrac{1}{2},\pi^{\infty \vee}_i)$ is holomorphic on $\mathcal{T}_{N_1}.$  Hence the claim implies that $Z_a(\varphi,\varphi',f,\underline{s})$ is also holomorphic on $\mathcal{T}_{N_1}.$  Since $N_1$ is arbitrary, we deduce that $Z_a(\varphi,\varphi',f,\underline{s})$ is holomorphic.   To deduce the bounds in the proposition, we multiply \eqref{Mellin} by $\prod_{i=1}^3L(s_i+\tfrac{1}{2},\pi^{\infty \vee}_i )$ and use the fact that $L(s_i+\tfrac{1}{2},\pi_i^{\infty \vee})$ has at worst polynomial growth in vertical strips \cite[\S 1]{BrumleyNarrow}.
\end{proof}

\begin{prop} \label{prop:second:fe}
Under the assumptions of Theorem \ref{thm:func:eqn}, one has
\begin{align*}
Z_a(\varphi,\varphi',f,\underline{s})=Z_{a^{-1}}(\varphi^\vee,\varphi'^\vee,\mathcal{F}(f),-\underline{s}).
\end{align*}
\end{prop}

\begin{proof}
Choose $(\underline{s},s) \in \CC^{\underline{1}} \times \CC$ that satisfy \eqref{bounds2}.  Let $\sigma'=\mathrm{Re}(\underline{s}).$
For any $a' \in \RR_{>0}^{\underline{1}},$ we have
\begin{align*}
    \frac{1}{2\pi i}\int_{\mathrm{Re}(\underline{s})=\sigma'}a'^{-\underline{s}}Z_a(\varphi,\varphi',f,\underline{s})ds=Z_{a',a}(\varphi,\varphi',f)
\end{align*}
by the definition \eqref{Zaa} of $Z_{a',a}(\varphi,\varphi',f)$ and Mellin inversion. The application of Mellin inversion is justified by 
 Lemma \ref{lem:conv:fixed:a} and Proposition \ref{prop:first:ac}.  Using Proposition \ref{prop:first:fe} we deduce that 
\begin{align*}
&\frac{1}{2\pi i}\int_{\mathrm{Re}(\underline{s})=\sigma'}a'^{-\underline{s}}Z_a(\varphi,\varphi',f,\underline{s})ds=Z_{a'^{-1},a^{-1}}(\varphi^\vee,\varphi'^\vee,\mathcal{F}(f)).
\end{align*}

By the same reasoning we have
\begin{align*}
\frac{1}{2\pi i} \int_{\mathrm{Re}(\underline{s})=-\sigma'}a'^{-\underline{s}} Z_{a^{-1}}(\varphi^\vee,\varphi'^\vee,\mathcal{F}(f),-\underline{s})ds=Z_{a'^{-1},a^{-1}}(\varphi^\vee,\varphi'^\vee,\mathcal{F}(f)).
\end{align*}
In this latter expression we can shift the contour to $\mathrm{Re}(\underline{s})=\sigma'$ using Proposition \ref{prop:first:ac}.  We then conclude the identity in the proposition by Mellin inversion.
\end{proof}

With this preparation in place, we can now execute the classic argument:

\begin{proof}[Proof of Theorem \ref{thm:func:eqn}]
By \eqref{MT} and Lemma \ref{lem:conv:fixed:a}, for $(\underline{s},s)$ satisfying \eqref{bounds2} we have
\begin{align*}
Z(\varphi,\varphi',f,\underline{s},s)&=\int_{\RR_{>0}}a^{s}Z_a(\varphi,\varphi',f,\underline{s})d^\times a\\
&=\int_{1}^\infty a^{s}Z_a(\varphi,\varphi',f,\underline{s})d^\times a+
\int_{0}^1 a^{s}Z_a(\varphi,\varphi',f,\underline{s})d^\times a.
\end{align*}
By Proposition \ref{prop:first:ac}, the integral from $1$ to $\infty$ admits a holomorphic continuation to $(\underline{s},s) \in \CC^{\underline{1}} \times \CC$.  On the other hand, using Proposition \ref{prop:second:fe} we have
\begin{align*}
\int_{0}^1a^{s}Z_a(\varphi,\varphi',f,\underline{s})d^\times a&=\int_{0}^1a^{s}Z_{a^{-1}}(\varphi^\vee,\varphi'^\vee,\mathcal{F}(f),-\underline{s})d^\times a\\
&=\int_{1}^\infty a^{-s}Z_a(\varphi^\vee,\varphi'^\vee,\mathcal{F}(f),-\underline{s})d^\times a.
\end{align*}
The latter integral admits a holomorphic continuation to the plane by Proposition \ref{prop:first:ac}.   Thus we have proved the analytic continuation.  To prove the functional equation we observe that
\begin{align*}
Z(\varphi,\varphi',f,\underline{s},s)=\int_{1}^\infty a^{s}Z_a(\varphi,\varphi',f,\underline{s})d^\times a+\int_{1}^\infty a^{-s}Z_a(\varphi^\vee,\varphi'^\vee,\mathcal{F}(f),-\underline{s})d^\times a
\end{align*}
is unchanged if we replace $(\varphi,\varphi',f,\underline{s},s)$ by $(\varphi^\vee,\varphi'^\vee,\mathcal{F}(f),-\underline{s},-s)$ by \ref{FTsquared}.
\end{proof}

To deduce the meromorphic continuation of $L(s,\pi,\otimes^3)$ from Conjecture \ref{PS:conj}, we require the following nonvanishing hypothesis:
\begin{enumerate}
    \myitem[(NV)] \label{nonvanish} If $v$ is non-Archimedean and $\pi_v,\pi_v'$ are unitary generic unramified representations of $\GL_{\underline{r}}(F_v)$ and $\GL_{\underline{r-2}}(F_v),$ then there exists $(f_v,W_v,W'_v) \in \mathcal{S}(X(F_v),\mathcal{H}) \times \mathcal{W}(\pi_v,\psi) \times \mathcal{W}(\pi'_v,\overline{\psi})$ such that $\mathcal{F}(f_v) \in \mathcal{S}(X^\circ(F_v) ,\mathcal{H})$ and $Z(W_v,W_v',f_v,\underline{s},s)$ is not identically zero.
\end{enumerate}
We expect that one can prove \ref{nonvanish} using our assumptions that the Fourier transform is anti-equivariant \ref{FT2} and preserves the basic function \ref{FT3}.

\begin{cor} \label{cor:mero}
 Assume \ref{nonvanish} and Conjecture \ref{PS:conj}.  Let $\pi$ be a cuspidal automorphic representation of $A_{\GL_{\underline{r}}} \backslash \GL_{\underline{r}}(\A_F).$  Then $L(s,\pi,\otimes^3)$ admits a meromorphic continuation to the plane.  
\end{cor}
\begin{proof}
Let $W^{\varphi}_{\psi}:=\prod_vW_v$ and $W^{\varphi'}_{\overline{\psi}}:=\prod_v W_v'$ be pure tensors and assume that $f=\otimes f_v \in \mathcal{S}(X(\A_F) ,\mathcal{H})$ is a pure tensor.  
Assume for the moment that $(\underline{s},s)$ satisfies \eqref{bounds2}.  By \ref{zeta:basic}, for a sufficiently large finite set of places $S$ one has that
\begin{align*}
Z(\varphi,\varphi',f,\underline{s},s)
=&\left(\prod_{v \in S}\frac{Z(W_v,W_v',f_v,\underline{s},s)}{L(s+\tfrac{1}{2},\pi_v^\vee,\otimes^3)\prod_{i=1}^3L(C_{i}(\underline{s},s)+\tfrac{1}{2},\pi_{iv}\times \pi_{iv}')L(s_i+\tfrac{1}{2},\pi_{iv}^\lor)}\right) \\& \times L(s+\tfrac{1}{2},\pi^\vee,\otimes^3)\prod_{i=1}^3L(C_i(\underline{s},s)+\tfrac{1}{2},\pi_i \times \pi_i')L(s_i+\tfrac{1}{2},\pi_i^\vee).
\end{align*}
The product $\prod_{i=1}^3L(C_i(\underline{s},s)+\tfrac{1}{2},\pi_i \times \pi_i')L(s_i+\tfrac{1}{2},\pi_i^\vee)$ is known to be entire. On the other hand, for each $v$
$$
\frac{Z(W_v,W_v',f_v,\underline{s},s)}{L(s+\tfrac{1}{2},\pi_v^\vee,\otimes^3)\prod_{i=1}^3L(C_{i}(\underline{s},s)+\tfrac{1}{2},\pi_{iv}\times \pi_{iv}')L(s_i+\tfrac{1}{2},\pi_{iv}^\lor)}
$$
is entire by  \ref{nA:ratio}, \ref{Arch:ratio}, and we can choose $W_v,W_v',f_v$ so that it is not identically zero by Lemma \ref{lem:non:van}.  We may additionally assume $\mathcal{F}(f_v) \in \mathcal{S}(X^{\circ}(F_v),\mathcal{H})$ by \ref{nonvanish}. Thus the corollary follows from Theorem \ref{thm:func:eqn}.
\end{proof}
\noindent 
We point out that the argument above will give a precise functional equation for $L(s,\pi,\otimes^3)$ once the corresponding local theory is complete.

\section{The fiber bundle method} \label{sec:fb}

In the previous section we proved that the Poisson summation formula of Conjecture \ref{PS:conj} implies the meromorphic continuation of $L(s,\pi,\otimes^3).$ It remains to prove Conjecture \ref{PS:conj}.  In this section we propose a method.

Since $Y \subset \mathcal{M} \times V_3$ we have a diagram
\begin{equation} \label{hat}
\begin{tikzcd}
X_{\underline{r}} \times \mathcal{M}& \arrow[l,swap,"\phi_1"] X \arrow[r,"\phi_2"] & V_3
\end{tikzcd}
\end{equation}
where $\phi_1$ and $\phi_2$ are the canonical maps.
Our strategy for proving the Poisson summation conjecture for $X$ is  to prove the Poisson summation conjecture for the fibers of the two maps and thereby deduce it for $X.$ After a bit of motivation from the case of a cone in \S \ref{ssec:cone}, we compute the fibers of $\phi_1$ and $\phi_2$ in \S\ref{sec:linear:FT}.   In both cases, the Poisson summation formula is available.

\subsection{The case of the cone} \label{ssec:cone}
The strategy just mentioned is motivated by the Poisson summation formula for the affine closure of  $U_3 \backslash \SL_3,$ where $U_3 <\SL_3$ is the unipotent radical of a Borel subgroup.   Let us describe it briefly to highlight some conceptual consequences.  We refer to \cite{BK:basic:affine,GurK:Cone}  for proofs and precise statements. 

Let $W=\GG_a^3$ and let $W^\vee=\mathrm{Hom}(\GG_a^3,\GG_a).$   The affine closure $\overline{U_3 \backslash \SL_3}$ may be identified with the following quadric cone:
$$
C(R):=\{(w^\vee,w) \in W^\vee(R) \times W(R): w^\vee(w)=0\}. 
$$
The relevant diagram in this case is
\begin{equation}
    \begin{tikzcd}
W^\vee & \arrow[l,swap,"\Psi_1"] C \arrow[r,"\Psi_2"] & W,
\end{tikzcd}
\end{equation}
where the arrows are the canonical projections.  The generic fibers of $\Psi_1$ and $\Psi_2$ are isomorphic to $\GG_a^2$ and carry canonical $\SL_2$-invariant inner products.  Thus we have Fourier transforms $\mathcal{F}_i$ along the fibers of $\Psi_i$ and corresponding Poisson summation formulae.  The Fourier transforms $\mathcal{F}_i$ generate a group isomorphic to the Weyl group of $\SL_3$ with respect to a split maximal torus.  Under this isomorphism, each $\mathcal{F}_i$ corresponds to a simple reflection.  Now $C$ also admits an action of the orthogonal group $\mathrm{O}_{W^\vee \oplus W}$ of the canonical quadratic form $W^\vee \otimes W \to \GG_a.$  The individual Fourier transforms $\mathcal{F}_i$ are not $\mathrm{O}_{W^\vee \oplus W}(\A_F)$-invariant, but the Fourier transform $\mathcal{F}_1 \circ \mathcal{F}_2 \circ \mathcal{F}_1$ (which corresponds to the long Weyl element) is invariant under $\mathrm{O}_{W^\vee \oplus W}(\A_F).$  The references above do not contain this precise statement, but the relevant computation in the finite field case is contained in \cite{Slipper}.

We draw two conclusions from this example.  First, suppose a  scheme is the total space of a fiber bundle such that fibers admit compatible families of Poisson summation formulae.  Then the whole scheme inherits a Poisson summation formula.  Second, if one constructs a group of Fourier transforms, certain compositions of the Fourier transforms may admit equivariance properties that the individual Fourier transforms do not enjoy.

\subsection{Poisson summation formulae on the fibers}
\label{sec:linear:FT} 

\subsubsection{The fibers of $\phi_1$} We first consider the fibers of $\phi_1.$   
Fix a place $v$ of the number field $F$ and omit it from notation, writing $F:=F_v.$  For $m \in \mathcal{M}^\circ(F),$ the fiber of $\phi_1^{-1}(x,m)$ may be identified with $Y_m(F),$ which is a vector space.  Therefore there is a corresponding Poisson summation formula.  We point out that the vector bundle $Y_{\mathcal{M}^\circ}(F)$ is not self-dual; this will likely cause additional difficulties.

 \subsubsection{Geometry of the fibers of $\phi_2$}

Recalling again that $Y \subset \mathcal{M} \times V_3,$ we have a canonical map $\phi:Y \to V_3.$  
If $v \in V_3(F),$ then 
$$
\phi_2^{-1}(v)=X_{\underline{r}} \times \phi^{-1}(v).
$$
Thus to compute the fibers of $\phi_2,$ it suffices to compute the fibers of $\phi.$
Using notation from \eqref{Yell}, let
\begin{align*}
\mathrm{pr}_Y:\mathbb{P}_1(Y) \lto \mathbb{P}^{\underline{1}}
\end{align*}
be the vector bundle whose fiber over $\ell$ is $Y_{\ell}.$ 

Recall the open $H$-orbit $O(v_0)\subset V_3$.
We have a commutative diagram
\begin{equation} \label{fib:comp}
\begin{tikzcd}
Y^\circ \ar[d,"p"]\ar[r,hook]&\mathcal{M}^\circ \times O(v_0)\arrow[r] \ar[d,"{\mathrm{pr}}"] &O(v_0) \arrow[d]\\
\mathbb{P}_1(Y)\ar[r,hook]&\mathbb{P}^{\underline{1}} \times V_3 \arrow[r] & V_3
\end{tikzcd}
\end{equation}
where $p$ is the canonical quotient map and the unmarked arrows are canonical.  All of the morphisms are $\GL_{\underline{r-2}} \times H^e$-equivariant, where $\GL_{\underline{r-2}} \times H^e$ acts through $H^e$ on $O(v_0),$ $V_3,$ and the bottom row of the diagram.  In particular, $p$ is the quotient by the action of $\GL_{\underline{r-2}}.$

\begin{lem} \label{lem:bir}
    The composite of the bottom two horizontal arrows of \eqref{fib:comp}
    admits an $H^e$-equivariant section $v \mapsto (v^\perp,v)$ over $O(v_0).$  The section is an isomorphism from $O(v_0)$ onto an open subscheme of $\mathbb{P}_1(Y).$
\end{lem}

\begin{proof} Using \eqref{Y:tilde} we see that
there is a unique $v^\perp_0 \in \mathbb{P}^1$ such that $(v_0^\perp,v_0) \in \mathbb{P}_1(Y).$  Moreover, using \eqref{Htildestab} we have $H_{v_0}^e=Z_{\GL_{\underline{2}}} \wedge^{\mathcal{E}} \mathcal{E}T_H N_0.$  This group fixes $v_0^\perp,$ hence $H_{v_0}^e=H_{(v_0^\perp,v_0)}^e.$  Thus the desired section is induced by action of $H^e$ and the choices of basepoints $v_0 \in O(v_0)(F)$ and $(v_0^\perp,v_0) \in \mathbb{P}_1(Y)(F).$
\end{proof}

\begin{prop} \label{prop:fibers}
For $v \in O(v_0)(F)$ with $v^{\perp}=(v_1^\perp,v_2^\perp,v_3^\perp),$
one has a canonical $\GL_{\underline{r-2}} \times H^e_{v}$-equivariant isomorphism
\begin{align*}
    \phi^{-1}(v) \tilde{\lto} \mathcal{M}_{1v_1^\perp} \times \mathcal{M}_{2v_2^{\perp}} \times \mathcal{M}_{3v_2^\perp}.
\end{align*}
\end{prop}
\begin{proof} In view of Lemma \ref{lem:bir}, the fibers of $\phi$ over $\OO(v_0)$ may be canonically identified with the fibers of $p$ over $\mathbb{P}_1(Y)_{O(v_0)}.$  These fibers may in turn be identified with the fibers of $\mathrm{pr}$ over the image of $\mathbb{P}_1(Y)_{O(v_0)}$ by the commutativity of the diagram.  The proposition follows.
\end{proof}

 Given Proposition \ref{prop:fibers}, to construct a Poisson summation formula for the fibers of $\phi_2$, it suffices to prove a Poisson summation formula for $X_{r_i}\times \mathcal{M}_{i\ell}$ for each $i$.  This Poisson summation formula is contained in \cite{GGH}.

\subsection{Remaining work}

To prove a Poisson summation formula for $X,$ we plan to combine the two Poisson summation formulae mentioned earlier.  This will involve constructing the Schwartz space $\mathcal{S}(X(\A_F),\mathcal{H})$ and piecing the Fourier transforms described above together to form a Fourier transform 
$$
\mathcal{F}:\mathcal{S}(X(\A_F),\mathcal{H}) \lto \mathcal{S}(X(\A_F),\mathcal{H})
$$
The main difficulties are to construct $\mathcal{F}$ so that
\begin{enumerate}
\item One can deduce the Poisson summation formula of Conjecture \ref{PS:conj} from the known Poisson summation formulae on the fibers of $\phi_1$ and $\phi_2.$
\item The Fourier transform $\mathcal{F}$ satisfies the equivariance property \ref{FT2}.
\end{enumerate}
These may be difficult problems, but we emphasize that, a priori, they may be addressed locally.

\bibliography{refs}{}
\bibliographystyle{alpha}

\printindex
\end{document}